\pdfoutput=1
\RequirePackage{ifpdf}
\ifpdf 
\documentclass[pdftex]{sigma}
\else
\documentclass{sigma}
\fi

\usepackage{subfigure}

\numberwithin{equation}{section}

\newtheorem{Theorem}{Theorem}[section]
\newtheorem*{Theorem*}{Theorem}

\newtheorem{Lemma}[Theorem]{Lemma}

 { \theoremstyle{definition}

\newtheorem{Remark}[Theorem]{Remark} }

\newcommand{\R}		{\mathbb{R}}
\newcommand{\C}		{\mathbb{C}}

\newcommand{\Z}		{\mathbb{Z}}

\newcommand{\dist}{\operatorname{dist}}

\newcommand{\diag}{\operatorname{diag}}
\newcommand{\re}{\operatorname{Re} }

\renewcommand{\arg}{\operatorname{arg}}
\renewcommand{\det}{\operatorname{det}}

\newcommand{\n}{{\vec n}}
\newcommand{\ic}{\mathrm{i}}
\newcommand{\rhy} {\textnormal{RHP}-${\boldsymbol Y}$}
\newcommand{\rhx} {\textnormal{RHP}-${\boldsymbol X}$}
\newcommand{\rhn} {\textnormal{RHP}-${\boldsymbol N}$}

\newcommand{\rhp} {\textnormal{RHP}-${\boldsymbol P}$}

\newcommand{\rhz} {\textnormal{RHP}-${\boldsymbol Z}$}

\newcommand{\rhpsi} {\textnormal{RHP}-${\boldsymbol \Psi}$}
\newcommand{\rhwpsi} {\textnormal{RHP}-${\boldsymbol \Psi_*}$}
\newcommand{\rhphi} {\textnormal{RHP}-${\boldsymbol \Phi}$}

\newcommand{\rhwphi} {\textnormal{RHP}-${\widetilde{\boldsymbol \Phi}}$}
\newcommand{\rhhphi} {\textnormal{RHP}-${\widehat{\boldsymbol \Phi}}$}
\newcommand{\rhth} {\textnormal{RHP}-${\boldsymbol \Theta}$}

\newcommand{\RS}{\boldsymbol{\mathfrak{S}}}
\newcommand{\z}	{{\boldsymbol z}}
\newcommand{\x}	{{\boldsymbol x}}

\begin{document}

\newcommand{\arXivNumber}{2411.04206}

\renewcommand{\PaperNumber}{033}

\FirstPageHeading

\ShortArticleName{Uniformity of Strong Asymptotics in Angelesco Systems}
\ArticleName{Uniformity of Strong Asymptotics\\ in Angelesco Systems}

\Author{Maxim L.~YATTSELEV}
\AuthorNameForHeading{M.L.~Yattselev}

\Address{Department of Mathematical Sciences, Indiana University Indianapolis,\\ 402~North Blackford Street, Indianapolis, IN 46202, USA}
\Email{\href{mailto:maxyatts@iu.edu}{maxyatts@iu.edu}}
\URLaddress{\url{https://math.indianapolis.iu.edu/~maxyatts/}}

\ArticleDates{Received November 08, 2024, in final form April 28, 2025; Published online May 08, 2025}

\Abstract{Let $ \mu_1 $ and $ \mu_2 $ be two complex-valued Borel measures on the real line such that $ \operatorname{supp} \mu_1 =[\alpha_1,\beta_1] < \operatorname{supp} \mu_2 =[\alpha_2,\beta_2] $ and $ {\rm d}\mu_i(x) = -\rho_i(x){\rm d}x/2\pi\ic $, where $ \rho_i(x) $ is the restriction to $ [\alpha_i,\beta_i] $ of a function non-vanishing and holomorphic in some neighborhood of~$ [\alpha_i,\beta_i] $. Strong asymptotics of multiple orthogonal polynomials is considered as their multi-indices $ (n_1,n_2) $ tend to infinity in both coordinates. The main goal of this work is to show that the error terms in the asymptotic formulae are uniform with respect to $ \min\{n_1,n_2\} $.}


\Classification{42C05; 41A20; 41A25}

\section{Main results}

\subsection{Multiple orthogonal polynomials}

Let $ \mu_1 $ and $ \mu_2 $ be two complex-valued Borel measures on the real line and $ \n=(n_1,n_2) $ be a~multi-index, where $ n_1$, $n_2 $ are non-negative integers. A non-identically zero polynomial $ P_{\n}(x) $ of degree at most $ |\n| := n_1+n_2 $ is called a type II multiple orthogonal polynomial with respect to a system of measures $ (\mu_1,\mu_2) $ if it satisfies
\begin{equation}\label{typeII}
\int P_\n(x)x^l{\rm d}\mu_i(x) =0, \qquad l\in\{0,\dots,n_i-1\}, \quad i\in\{1,2\}.
\end{equation}
In what follows, we take $ P_\n(x) $ to be the monic polynomial of minimal degree satisfying \eqref{typeII}, which makes it unique. Type~I multiple orthogonal polynomials are not identically zero polynomial coefficients of the linear form
\begin{gather}
Q_\n(x):= A_\n^{(1)}(x){\rm d}\mu_1(x) + A_\n^{(2)}(x){\rm d}\mu_2(x), \qquad \deg A_\n^{(i)} < n_i, \nonumber \\
\int x^lQ_\n(x) = 0, \qquad l<|\n|-1, \qquad A_{(0,1)}^{(1)} = A_{(1,0)}^{(2)} \equiv 0.\label{typeI}
\end{gather}
It is known \cite[Section~23.1]{Ismail} that when $\deg P_\n=|\n|$, in which case the multi-index $ \n $ is called \emph{normal}, the linear form $ Q_\n(x) $ is defined uniquely up to multiplication by a constant. In this case, it is customary to normalize it by requiring
\begin{equation}
\label{Qn_norm}
\int x^{|\n|-1}Q_{\n}(x)=1 .
\end{equation}
The polynomials $ A_\n^{(i)}(x) $ are no longer monic and their leading coefficients, which we denote by~$ 1/h_{\n-\vec e_i,i} $, are closely related to the type~II polynomials $ P_\n(x) $, where $ \vec e_1=(1,0) $ and $ \vec e_2 = (0,1) $. Indeed, it holds that
\begin{align}
\label{connection}
h_{\n-\vec e_i,i} &= h_{\n-\vec e_i,i} \int x^{|\n|-1}Q_\n(x) = h_{\n-\vec e_i,i}\int P_{\n-\vec e_i}(x)Q_\n(x) \nonumber \\
& = h_{\n-\vec e_i,i}\int P_{\n-\vec e_i}(x) A_\n^{(i)}(x){\rm d}\mu_i(x) = \int P_{\n-\vec e_i}(x) x^{n_i-1}{\rm d}\mu_i(x)
\end{align}
as follows from \eqref{typeII}--\eqref{Qn_norm}. It is known \cite[Theorem~23.1.11]{Ismail} that if indices $ \n $ and $ \n+\vec e_i $ are normal, multiple orthogonal polynomials satisfy nearest-neighbor recurrence relations of the form
\begin{gather}
xP_\n(x) = P_{\n+\vec e_i}(x) + b_{\n,i}P_\n(x) + a_{\n,1}P_{\n-\vec e_1}(x) + a_{\n,2}P_{\n-\vec e_2}(x), \nonumber \\
xQ_\n(x) = Q_{\n-\vec e_i}(x) + b_{\n-\vec e_i,i}Q_\n(x) + a_{\n,1}Q_{\n+\vec e_1}(x) + a_{\n,2}Q_{\n+\vec e_2}(x).\label{rec_rel}
\end{gather}
All these definitions can be formulated for a collection of more than two measures, however, we shall not pursue such an extension here.

\subsection{Angelesco systems}

Assume that each measure $ \mu_i $ is compactly supported and let $ \Delta_i $ be the smallest closed interval containing the support of $ \mu_i $. If $ \Delta_1\cap \Delta_2 = \varnothing $, then the pair $ (\mu_1,\mu_2) $ is said to form an Angelesco system. Angelesco himself considered the case of non-negative measures and had shown that such systems are always \emph{perfect} (all multi-indices are normal) \cite{Ang19}. In what follows, we are only interested in the case where
\begin{equation}
\label{mus}
\operatorname{supp} \mu_i = \Delta_i =: [\alpha_i,\beta_i] \qquad \text{and}\qquad  {\rm d}\mu_i(x) = -\rho_i(x)\frac{{\rm d}x}{2\pi\ic},
\end{equation}
i.e., each $ \mu_i $ is supported by an interval and is absolutely continuous with respect to the Lebesgue measure. We allow densities $ \rho_i(x) $ to be complex-valued\footnote{The specific choice of the normalization in \eqref{mus} is there for two reasons. First, under such a normalization the Markov function of $ \mu_i $ becomes the Cauchy transform of $ \rho_i $, i.e., $ \int\frac{{\rm d}\mu_i(x)}{z-x}=\int\frac{\rho_i(x)}{x-z}\frac{{\rm d}x}{2\pi\ic} $. Second, when $ \mu_i $ is a~positive measure, the density $ \log(\rho_iw_{c,i+}) $, appearing in \eqref{Szego}, is a real-valued function.} and assume for definiteness that~${ \beta_1<\alpha_2} $.

When each $ \ic\rho_i(x) $ is positive almost everywhere on the corresponding interval $ \Delta_i $, the zero distribution of the polynomials $ P_\n(x) $ and their $ |\n| $-th root (i.e., weak) asymptotic behavior were studied in \cite{GRakh81} along subsequences of indices that satisfy
\begin{equation}
\label{subseq}
\lim_{|\n|\to\infty}n_1/|\n| \qquad \text{exists and belongs to} \ (0,1).
\end{equation}
When the functions $ \ic\rho_i(x) $ are non-negative with integrable logarithms (Szeg\H{o} class), strong asymptotics of these polynomials along the diagonal sequence $ n_1=n_2 $ was obtained in \cite{Ap88}, see also \cite{Kal79,NutTr87}, and more generally under assumption \eqref{subseq} in \cite{uApDenYa}. When each function $ \rho_i(x) $ is the product of a restriction to $ \Delta_i $ of a non-vanishing (complex-valued) holomorphic function and a Fisher--Hartwig weight, the strong asymptotics of type~II polynomials along subsequences satisfying \eqref{subseq} was derived in \cite{Ya16}. Asymptotics of type~I polynomials as well as of the recurrence coefficients was later deduced in \cite{ApDenYa20}, but just for weights that are restrictions of holomorphic functions only. Moreover, assuming more stringently that each density $ \rho_i(x) $ is a restriction of a~holomorphic function and is positive on $ \Delta_i $ while allowing the limit in \eqref{subseq} to be $ 0 $ or $ 1 $ under the additional assumption
\begin{equation}
\label{vareps}
\varepsilon_\n := 1/\min\{n_1,n_2\} \to 0 \qquad \text{as} \ |\n|\to\infty,
\end{equation}
strong asymptotics of the polynomials of both types and the asymptotics of their recurrence coefficients were derived in \cite{ApDenYa21}. The goal of this work is to show that the error rates obtained in \cite{ApDenYa21} can be made uniform in $ \n $; that is, the asymptotic formulae can be derived solely under condition \eqref{vareps}, no assumption on the existence of the limit in \eqref{subseq} is needed. We shall assume that the densities $ \rho_i(x) $ are restrictions of analytic, non-vanishing, and in general complex-valued functions. Analyticity assumption can be relaxed, but this will be addressed in a subsequent publication.

\subsection{Riemann surfaces}

The functions describing strong asymptotics of multiple orthogonal polynomials naturally live on a sequences of Riemann surfaces. To describe these surfaces, we need to start with the already mentioned work by Gonchar and Rakhmanov \cite{GRakh81}. There, assuming \eqref{mus} and that each $ \ic\rho_i(x) $ is almost everywhere positive function on $ \Delta_i $, it was shown that if a subsequence of multi-indices satisfies \eqref{subseq} and $ c $ is the limit, then the normalized counting measures of the zeros of $ P_\n(x) $ converge weak$^*$ to a certain measure $ \omega_c $ and
\[
|\n|^{-1} \log |P_\n(z)| = - (1+o(1))V^{\omega_c}(z)
\]
locally uniformly in $ \overline\C\setminus(\Delta_1\cup\Delta_2) $ along this subsequence of multi-indices, where $ V^\omega(z):=-\int\log|z-x|{\rm d}\omega(x) $ is the logarithmic potential of a measure $ \omega $ and $ \omega_c := \omega_{c,1} + \omega_{c,2} $ with~$ \omega_{c,1}$,~$\omega_{c,2} $ being the unique pair of measures such that $ |\omega_{c,1}| = c $, $ |\omega_{c,2}|=1-c $ (here, $ |\omega| $ is the total mass~$ \omega $), $ \operatorname{supp} \omega_{c,i}\subseteq \Delta_i $, and
\[
\begin{cases}
\ell_{c,i} - V^{\omega_c+\omega_{c,i}}(x) \equiv 0, & x\in \operatorname{supp} \omega_{c,i}, \\
\ell_{c,i} - V^{\omega_c+\omega_{c,i}}(x) < 0, & x\in\Delta_i\setminus \operatorname{supp} \omega_{c,i},
\end{cases}
\]
 $ i\in\{1,2\} $, for some constants $ \ell_{c,1}$, $\ell_{c,2} $ (the measures $ \omega_{c,1}$, $\omega_{c,2} $ can also be defined through a~certain vector energy minimization problem).

 What is of main importance to us from the above results are the supports of the vector equilibrium measures $ \omega_{c,1} $ and $ \omega_{c,2} $. It was explained in \cite{GRakh81} that
\[
\Delta_{c,i} := \operatorname{supp} \omega_{c,i} = [\alpha_{c,i},\beta_{c,i}] \subseteq \Delta_i=[\alpha_i,\beta_i]
\]
 $ i\in\{1,2\} $, where $ \alpha_{c,1} := \alpha_1 $ and $ \beta_{c,2} := \beta_2 $ for any $ c\in(0,1) $. However, it is possible that $ \beta_{c,1}<\beta_1 $ or $ \alpha_{c,2}>\alpha_2 $ (this is so-called \emph{pushing effect}). In fact, it is known \cite[Propositions~4.1 and 4.2]{ApDenYa21} that there exist $ 0<c^*<c^{**}<1 $ such that
\begin{equation}
\label{cstars}
\begin{cases}
\beta_{c,1}<\beta_1, & c<c^*, \\
\beta_{c,1}=\beta_1, & c \geq c^*,
\end{cases}
\qquad \text{and}\qquad
\begin{cases}
\alpha_{c,2}=\alpha_2, & c \leq c^{**}, \\
\alpha_{c,2}>\alpha_2, & c > c^{**}.
\end{cases}
\end{equation}
Moreover, $ \beta_{c,1} $ is a continuous strictly increasing function of $ c $ on $ [0,c^*] $ with $ \beta_{0,1}:= \alpha_1 $ while~$ \alpha_{c,2} $ is a continuous strictly increasing function of $ c $ on $ [c^{**},1] $ with $ \alpha_{1,2}:= \beta_2 $. It is also known that the constants $ \ell_{c,1} $ and $ \ell_{c,2} $ are continuous functions of $ c $ and so are the measures $ \omega_{c,1} $ and~$ \omega_{c,2} $ in the sense of weak$^*$ convergence of measures.

\begin{figure}[ht!]\centering
\includegraphics[scale=.5]{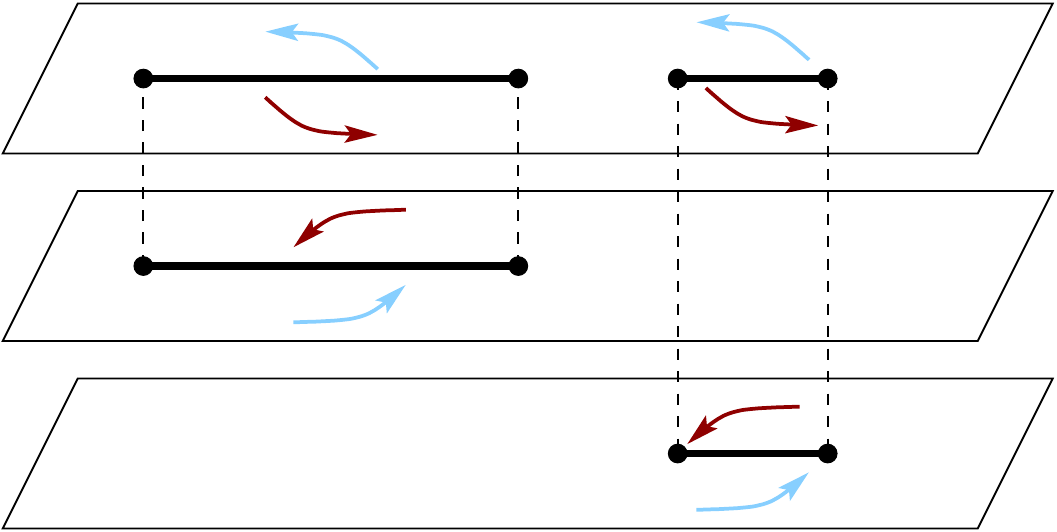}
\begin{picture}(0,0)
\put(-238,110){$\alpha_1$}
\put(-128,110){$\beta_{c,1}$}
\put(-109,110){$\alpha_2$}
\put(-53,110){$\beta_2$}
\put(-30,115){$\RS_c^{(0)}$}
\put(-30,69){$\RS_c^{(1)}$}
\put(-30,24){$\RS_c^{(2)}$}
\end{picture}
\caption{Realization of the surface $ \RS_c $ (the notation is meant to emphasize that pushing effect happens on one interval only, i.e., $ \beta_{c,1}<\beta_1 $ and therefore $\alpha_{c,2}=\alpha_2$ for this value of $ c $).}
\label{fig:surfaceA}
\end{figure}

Given $ c\in(0,1) $, let $\RS_c$ be a 3-sheeted Riemann surface realized as follows. Define
\[
\Delta_c := \Delta_{c,1}\cup \Delta_{c,2}, \qquad E_c := E_{c,1} \cup E_{c,2}, \qquad \text{where} \quad E_{c,i} := \{\alpha_{c,i},\beta_{c,i}\}.
\]
Denote by $\RS_c^{(0)}$, $\RS_c^{(1)}$, and $\RS_c^{(2)}$, three copies of $ \overline \C $ cut along $ \Delta_c $, $ \Delta_{c,1} $, and $ \Delta_{c,2} $, respectively. These copies are then glued to each other crosswise along the corresponding cuts, see Figure~\ref{fig:surfaceA}. It can be easily seen from the Riemann–-Hurwitz formula that $ \RS_c $ has genus $ 0 $. We denote by~$\pi$ the natural projection from $\RS_c$ to $\overline\C$ and employ the notation $ \z $ for a generic point on $ \RS_c $ with~${ \pi(\z)=z }$. We let $ \boldsymbol\alpha_1= \boldsymbol\alpha_{c,1}, \boldsymbol\beta_{c,1}, \boldsymbol\alpha_{c,2},\boldsymbol\beta_2 = \boldsymbol\beta_{c,2} $ to stand for the ramification points of $ \RS_c $ with natural projections $ \alpha_1$, $\beta_{c,1}$, $\alpha_{c,2}$, $\beta_2 $, respectively. We set
\begin{gather*}
\boldsymbol E_c := \boldsymbol E_{c,1} \cup \boldsymbol E_{c,2},\qquad \boldsymbol E_{c,i} = \{ \boldsymbol\alpha_{c,i}, \boldsymbol\beta_{c,i}\}, \\
\boldsymbol \Delta_c:= \boldsymbol \Delta_{c,1} \cup \boldsymbol \Delta_{c,2}, \qquad \boldsymbol \Delta_{c,i} := \RS_c^{(0)}\cap \RS_c^{(i)} = \partial\RS_c^{(i)}.
\end{gather*}
Notice that $ \boldsymbol \Delta_{c,i}\setminus \boldsymbol E_{c,i} $ is a two-to-one cover of $ \Delta_{c,i}^\circ := (\alpha_{c,i},\beta_{c,i}) $.

We call a linear combination $ \sum m_i\z_i$, $ m_i\in\Z $, a divisor. Its degree is defined as $ \sum m_i $. We say that $ \sum m_i\z_i$ is a zero/pole divisor of a meromorphic function if this function has a zero at $ \z_i $ of multiplicity $ m_i $ when $ m_i>0 $, a pole at $ \z_i $ of order $ -m_i $ when $ m_i<0 $, and has no other zeros or poles in the domain of its definition. Zero/pole divisors of rational functions on~$ \RS_c $ necessarily have degree zero. Conversely, since $\RS_c$ has genus $ 0 $, any degree zero divisor is a~zero/pole divisor of a~rational function, which is unique up to multiplication by a constant. For any function $ G(\z) $, defined on $ \RS_c\setminus\boldsymbol \Delta_c $, we denote by
\[
G^{(k)}(z) := G(\z), \qquad \z\in \RS_c^{(k)}\setminus \partial\RS_c^{(k)}, \qquad k\in\{0,1,2\},
\]
its pull-back under the natural projection from the $ k $-th sheet to the cut complex plane.

\subsection{Conformal maps}

To proceed, we introduce a certain conformal map of $ \RS_c $ onto $ \overline \C $, say $ \chi_c(\z) $, which is a rational function with one pole and one zero, defined by the relation
\begin{equation}
\label{chi}
\chi_c^{(0)}( z ) = z + \mathcal O\bigl(z^{-1} \bigr) \qquad \text{as} \ z\to\infty.
\end{equation}
Since prescribing the absence of a constant term around $ \infty^{(0)} $ ($ \infty^{(k)} $ is a point on $ \RS_c^{(k)} $ whose natural projection is the point at infinity) is equivalent to prescribing a zero, the function $ \chi_c(\z) $ is uniquely determined by \eqref{chi}. Further, let the numbers $ A_{c,1}$, $A_{c,2}$, $B_{c,1}$, $B_{c,2} $ be defined by
\begin{equation}
\label{AngPar1}
\chi_c^{(i)}(z) =: B_{c,i} + A_{c,i}z^{-1} + \mathcal O\bigl(z^{-2} \bigr) \qquad \text{as} \ z\to\infty, \ i\in\{1,2\}.
\end{equation}
Since $ \chi_c(\z) $ is a conformal map, the numbers $ A_{c,i} $ are necessarily non-zero (otherwise $ \chi_c(\z)- B_{c,i} $ would have had a double zero). Moreover, by tracing the image of $ \pi^{-1}(\R)\setminus\boldsymbol\Delta_c $ under $ \chi_c(\z) $, which is necessarily equal to the real line, one gets that each $ A_{c,i}>0 $ and $ B_{c,2}>B_{c,1} $. Set
\begin{gather}
\varphi_i(z) := (z-(\beta_i+\alpha_i)/2 + w_i(z))/2, \qquad
w_i(z) :=\sqrt{(z-\alpha_i)(z-\beta_i)} = z + \mathcal O(1)\label{phi_w}
\end{gather}
to be the branches holomorphic off $ \Delta_i $, $ i\in\{1,2\} $ ($ \varphi_i(z) $ is the conformal map of the complement of $ \Delta_i $ onto the complement of $ \{|z|>(\beta_i-\alpha_i)/4 \} $ that behaves like $ z $ at infinity). It was also shown in \cite[Proposition~2.1]{ApDenYa21} that $ A_{c,i} $ and $ B_{c,i} $ are (real-valued) continuous functions of the parameter $ c\in (0,1) $ such that
\begin{equation}
\lim_{c\to0}
\begin{cases}
A_{c,2} = \big[(\beta_2-\alpha_2)/4\big]^2 =: A_{0,2}, \\
B_{c,2} = (\beta_2+\alpha_2)/2 =: B_{0,2}, \\
A_{c,1} = 0 =: A_{0,1}, \\
B_{c,1} = B_{0,2}+\varphi_2(\alpha_1) =: B_{0,1},
\end{cases}\label{AngPar2}
\end{equation}
and the analogous limits hold when $ c\to1 $.

In the limit $ c\to0 $ (similar considerations apply to the case $ c \to 1 $), the Riemann surface $ \RS_c $ becomes disconnected: one connected component is a copy of $ \overline\C $ and the other one is two copies of $ \overline \C $ glued together crosswise across $ \Delta_2 $. As expected, in this case the conformal maps $ \chi_c(\z) $ converge to a conformal map of the second connected component onto the Riemann sphere. More precisely, it was shown in \cite[equation~(5.2)]{ApDenYa21} that
\[
\chi_c(\z) - B_{c,2} = (1+o(1))
\begin{cases}
\varphi_2(z), & \z\in \RS_c^{(0)}, \\
\varphi_2(\alpha_1), & \z\in \RS_c^{(1)}, \\
\left(\dfrac{\beta_2-\alpha_2}{4}\right)^2/\varphi_2(z), & \z\in \RS_c^{(2)},
\end{cases}
\]
uniformly on $ \RS_c $ as $ c\to0 $ (the top two limits in \eqref{AngPar2} follow immediately from this formula).

The vector equilibrium measures discussed before admit the following explicit formulae. Let~$ \Pi_c(\z) $ be the derivative of $ \chi_c(\z) $, that is, the derivative of $ \chi_c^{(k)}(z) $ is equal to $ \Pi_c^{(k)}(z) $, for each~${k\in\{0,1,2\}}$. Equivalently, $ \Pi_c(\z) $ is a rational function on $ \RS_c $ with the zero/pole divisor and normalization given by
\begin{equation}
\label{Pin}
2\bigl(\infty^{(1)} + \infty^{(2)} \bigr) - \boldsymbol\alpha_1 - \boldsymbol\beta_{c,1} - \boldsymbol\alpha_{c,2} - \boldsymbol\beta_2 \qquad \text{and}\qquad  \Pi_c^{(0)}(\infty) =1.
\end{equation}
Observe that $ \lim_{z\to\infty} z^2\Pi_c^{(i)}(z) = -A_{c,i} $, $ i\in\{1,2\} $. Set
\begin{equation}
\label{h_chi}
h_c(\z) := \Pi_c(\z) \frac{\chi_c(\z) - (1-c)B_{c,1} - c B_{c,2}}{(\chi_c(\z)-B_{c,1})(\chi_c(\z)-B_{c,2})}.
\end{equation}
Then $ h_c(\z) $ is a rational function on $ \RS_c $ with the zero/pole divisor given by
\begin{equation}
\label{hc1}
\infty^{(0)} + \infty^{(1)} + \infty^{(2)} + \z_c - \boldsymbol\alpha_1 - \boldsymbol\beta_{c,1} - \boldsymbol\alpha_{c,2} - \boldsymbol\beta_2,
\end{equation}
where $ \z_c $ is some point on $ \RS_c $ whose existence is guaranteed by the fact that zero/pole divisors of rational functions on compact Riemann surfaces must have degree zero. Moreover, it holds that
\begin{equation}
\label{hc2}
\lim_{z\to\infty} zh_c^{(0)}(z) = 1, \qquad \lim_{z\to\infty} zh_c^{(1)}(z) = -c, \qquad \text{and}\qquad  \lim_{z\to\infty} zh_c^{(2)}(z) = c-1
\end{equation}
(this, in particular, means that all three branches of $ h_c(\z) $ add up to the identically zero function because their sum must be an entire function that is equal to zero at infinity). The function~$ h_c(\z) $ is uniquely determined by \eqref{hc1} and \eqref{hc2}. Indeed, the ratio of any two such functions (that is, functions corresponding to possibly different points $ \z_c $) minus $ 1 $ must have at most one necessarily simple pole by \eqref{hc1} and at least three zeros (at the points on top of infinity) by \eqref{hc2}, which is only possible for the identically zero function. Now, it follows from \cite[Proposition~2.3]{Ya16}, where the functions $ h_c(\z) $ were introduced differently, that
\begin{equation}
\label{point_zc}
\z_c\in \RS_c^{(0)} \qquad \text{and}\qquad
\begin{cases}
\z_c=\boldsymbol\beta_{c,1}, & c\in(0,c^*], \\
\pi(\z_c) \in (\beta_1,\alpha_2), & c\in(c^*,c^{**}), \\
\z_c=\boldsymbol\alpha_{c,2}, & c\in[c^*,1)
\end{cases}
\end{equation}
(if either $ \z_c=\boldsymbol\beta_{c,1} $ or $ \z_c=\boldsymbol\alpha_{c,2} $, then these points cancel each other out in \eqref{hc1} and $ h_c(\z) $ has only three poles and three zeros). Finally, the following claim can be found in \cite[Section~4]{ApDenYa21}. It holds that
\begin{equation}
\label{eq_meas}
{\rm d}\omega_{c,i}(x) = \bigl( h_{c+}^{(i)}(x)-h_{c-}^{(i)}(x) \bigr)\frac{{\rm d}x}{2\pi\ic}, \qquad x\in\Delta_{c,i}, \ i\in\{1,2\}.
\end{equation}

Notice that $ h_{c\pm}^{(i)}(x) = h^{(0)}_{c\mp}(x) $ for $ x\in\Delta_{c,i} $. Using the limiting behavior of $ \chi_c^{(0)}(z) $ as $ c\to0 $ discussed above, one can show that $ h_c(\z) $ converges to $ 1/w_2(z) $ on the zero-th sheet as $ c\to0 $. This, in particular, implies that $ \omega_{c,2} $ converge to the arcsine distribution on $ \Delta_2 $ as $ c\to 0 $ (of course, similar limits take place as $ c\to1 $).

\subsection{Main asymptotic terms}

The following construction was carried out in \cite[Section~6]{Ya16}. Similarly to \eqref{phi_w}, for each $ c\in(0,1) $, set
\begin{equation}
\label{wci}
w_{c,i}(z) := \sqrt{(z-\alpha_{c,i})(z-\beta_{c,i})}, \qquad z\in \C\setminus\Delta_{c,i}, \ i\in\{1,2\},
\end{equation}
to be the branch normalized so that $ w_{c,i}(z) = z + \mathcal O(1) $ as $ z\to\infty $. Denote by $ \mathcal C_\z $ the discontinuous Cauchy kernel on $ \RS_c $, that is, $ \mathcal C_\z $ is the third kind differential with three simple poles, located at $ \z $ and the other two points with the same natural projection, and the residues $ 2 $ at $ \z $ and $ -1 $ at the other points. Put
\begin{equation}
\label{Szego}
S_c(\z) := \exp\left\{ \frac1{6\pi\ic} \sum_{i=1}^2 \int_{\boldsymbol\Delta_{c,i}} \log(\rho_iw_{c,i+})\mathcal C_\z\right\},
\end{equation}
where we choose a continuous determination of $ \log\rho_i(x) $ and set
\[
\log w_{c,i+}(x) := \log |w_{c,i}(x)| + \pi\ic/2
\]
(when $ \mu_i $ is a~positive measure, we can take $ \log\rho_i(x) = \log|\rho_i(x)| - \pi\ic/2$, see \eqref{mus}, so that $ (\rho_i w_{c,i+})(x)$ is a positive function on $ (\alpha_{c,i},\beta_{c,i}) $). It is known \cite[Proposition~2.4]{Ya16} that $ S_c(\z) $ is holomorphic in $ \RS_c\setminus\boldsymbol\Delta_c $ and has continuous traces on $ \boldsymbol\Delta_{c,i}\setminus \boldsymbol E_{c,i} $ that satisfy
\begin{gather}
S_{c\pm}^{(i)}(x) = S_{c\mp}^{(0)}(x)(\rho_iw_{c,i+})(x), \qquad x\in \Delta_{c,i}^\circ, \nonumber \\
\big|S_c^{(0)}(z)\big| \sim |z-e|^{-1/4} \qquad \text{as} \ z\to e\in E_c.\label{szego-pts2}
\end{gather}
Moreover, $ \bigl(S_c^{(0)}S_c^{(1)}S_c^{(2)}\bigr)(z) \equiv 1 $. These functions continuously depend on the parameter $ c\in(0,1) $ and possess limits as $ c\to 0 $ and $ c\to 1 $. Namely, we have, see \cite[Proposition~3.1]{ApDenYa21}, that
\[
\frac{S_c^{(k)}(z)}{S_c^{(k)}(\infty)} = ( 1 + o(1) )
\begin{cases}
S_{\rho_2}(z)/S_{\rho_2}(\infty), & k=0, \\
1, & k=1, \\
S_{\rho_2}(\infty)/S_{\rho_2}(z), & k=2,
\end{cases}
\]
as $ c\to0 $, where $ o(1) $ holds locally uniformly in $ \C\setminus\{\alpha_1\} $ when $ k\in\{0,1\} $ and uniformly in $ \overline\C $ when $ k=2 $, and $ S_{\rho_2}(z) $ is the classical Szeg\H{o} function of $ \rho_2(x) $, that is,
\begin{equation}
\label{classical_szego}
S_{\rho_2}(z) := \exp\left\{ \frac{w_2(z)}{2\pi\ic}\int_{\Delta_2}\frac{\log(\rho_2w_{2+})(x)}{z-x}\frac{{\rm d}x}{w_{2+}(x)} \right\}.
\end{equation}
Moreover, it holds that the limits of $ S_c^{(0)}(\infty)c^{1/3} $, $ S_c^{(1)}(\infty)c^{-2/3} $, and $ S_c^{(2)}(\infty)c^{1/3} $ exist and are non-zero as $ c\to0 $. As usual, the above results have their counterparts when $ c\to1 $.

The terms describing the geometric growth of multiple orthogonal polynomials can be most conveniently defined as rational functions on the surfaces corresponding to rational values of the parameter $ c $. Namely, let $ c(\n) := n_1/|\n| $ and set
\begin{gather}
\Phi_\n(\z) := \tau_\n(\chi_{c(\n)}(\z)-B_{c(\n),1})^{n_1}(\chi_{c(\n)}(\z)-B_{c(\n),2})^{n_2}, \nonumber \\
\tau_\n^3 := (-1)^{n_2}A_{c(\n),1}^{-n_1}A_{c(\n),2}^{-n_2}(B_{c(\n),2}-B_{c(\n),1})^{-|\n|},\label{Phin}
\end{gather}
where we arbitrarily fix a cubic root of $ \tau_\n $. Thus defined, the function $ \Phi_\n(\z) $ is rational on $ \RS_{c(\n)} $ with the zero/pole divisor and the normalization given by
\begin{equation}
\label{Phin1}
n_1\infty^{(1)} + n_2\infty^{(2)} - |\n|\infty^{(0)} \qquad \text{and}\qquad  \Phi_\n^{(0)}(z)\Phi_\n^{(1)}(z)\Phi_\n^{(2)}(z) \equiv 1
\end{equation}
(to see that $ \Phi_\n(\z) $ is normalized this way is enough to observe that the product of all three branches is necessarily an entire function that assumes value $ 1 $ at infinity by \eqref{AngPar1}). Equivalently, it holds that
\begin{equation}
\label{Phin2}
\Phi_\n(\z) = \exp\left\{ |\n| \left(\int_{\boldsymbol \beta_2}^\z h_{c(\n)}(\x) {\rm d}x - \frac13\int_{\boldsymbol \beta_2}^{\beta_2^{(1)}} h_{c(\n)}(\x) {\rm d}x\right) \right\},
\end{equation}
where \smash{$ \beta_2^{(1)}\in\RS_{c(\n)}^{(1)} $} with $ \pi\bigl(\beta_2^{(1)}\bigr) = \beta_2 $, since the right-hand side of \eqref{Phin2} is a well-defined rational function on $ \RS_{c(\n)} $ with the divisor and normalization given by \eqref{Phin1} due to \eqref{hc2}. Let us point out that it is not hard to argue using \eqref{eq_meas}, see \cite[Proposition~2.1]{Ya16}, that
\begin{equation}
\label{Phin3}
\log|\Phi_\n(\z)| = |\n|\begin{cases}
-V^{\omega_{c(\n),1}+\omega_{c(\n),2}}(z)+ (\ell_{\n,1}+\ell_{\n,2})/3, & \z\in\RS_{c(\n)}^{(0)}, \\
V^{\omega_{c(\n),i}}(z)+ (\ell_{\n,3-i}-2\ell_{\n,i})/3, & \z\in\RS_{c(\n)}^{(i)}, \ i\in\{1,2\}.
\end{cases}
\end{equation}

\subsection{Main results}

Given Szeg\H{o} functions $ S_c(\z) $ as well as functions $ \Phi_\n(\z) $ and $ \Pi_c(\z) $ introduced in \eqref{Phin} and~\eqref{Pin}, respectively, we are ready to state our main results. Recall \eqref{vareps}. We start by describing the asymptotic behavior of type~II polynomials.

\begin{Theorem}
\label{thm:1}
Let $ \mu_1 $ and $ \mu_2 $ be as in \eqref{mus}, where $ \rho_1(x) $ and $ \rho_2(x) $ are the restrictions to~$ \Delta_1 $ and~$ \Delta_2 $, respectively, of non-vanishing functions analytic in some neighborhood of the corresponding interval. Set
\[
\mathcal P_{\n,i}(z) := \gamma_{\n,i} / \bigl(S_{c(\n)}\Phi_\n\bigr)^{(i)}(z), \qquad \lim_{z\to\infty} \mathcal P_{\n,i}(z)z^{-n_i} = 1,
\]
where $ c(\n) = n_1/|\n| $ and $ \gamma_{\n,i} $ are the normalizing constants, $ i\in\{1,2\} $. Let $ P_\n(z) $ be the type II multiple orthogonal polynomial defined via \eqref{typeII}. Then for all $ \varepsilon_\n $ small enough we can write~${P_\n(z) = P_{\n,1}(z)P_{\n,2}(z)}$ with the monic polynomials $ P_{\n,i}(z) $ satisfying
\begin{equation}
P_{\n,i}(z) = (1+ o(1)) \mathcal P_{\n,i}(z), \qquad
P_{\n,i}(x) = (1+ o(1)) \mathcal P_{\n,i+}(x) + (1+ o(1)) \mathcal P_{\n,i-}(x),\label{main1}
\end{equation}
uniformly for $ \dist(z,\Delta_{c(\n),i})\geq d $, $ z\in\overline\C $, and $ \dist(x,E_{c(\n),i})\geq d $, $ x\in\Delta_{c(\n),i} $, respectively, for any $ d>0 $ fixed, $ i\in \{1,2\} $. The error terms in the above formulae depend on $ d $ and satisfy
\[
\text{either} \qquad o(1) = \mathcal O\bigl( \varepsilon_\n^{1/3} \bigr) \qquad \text{or} \qquad o(1) = \mathcal O\bigl( \varepsilon_\n \bigr)
\]
uniformly for all $ \varepsilon_\n $ sufficiently small, where the second estimate holds if we additionally assume that $ c(\n) $ is uniformly separated from $ c^*$, $c^{**} $.
\end{Theorem}

Theorem~\ref{thm:1} generalizes \cite[Theorem~3.2]{ApDenYa21} in a sense that it provides asymptotic formulae for individual factors $ P_{\n,i}(z) $ rather the whole product $ P_\n(z) $ while the formulae themselves are uniform in $ \n $, that is, do not depend on ray sequences \eqref{subseq}.

\begin{Remark}
It easily follows from \eqref{main1} that each polynomial $ P_{\n,i}(z) $ has exactly $ n_i $ zeros and they all belong to $ \dist(z,\Delta_{c(\n),i})<d $ for any $ d>0 $ and all $ \varepsilon_\n $ small enough. Of course, if the measures $ \mu_1 $ and $ \mu_2 $ are arbitrary positive, it is straightforward to show that each $ P_{\n,i}(x) $ has exactly $ n_i $ zeros and they all belong to $ \Delta_i $, as initially has been observed in \cite{Ang19}.
\end{Remark}

\begin{Remark}
It readily follows from \eqref{Phin} that the constants $ \gamma_{\n,i} $ can be expressed as
\[
\gamma_{\n,i} = \tau_\n A_{c(\n),i}^{n_i} (B_{c(\n),i} - B_{c(\n),3-i})^{n_{3-i}} S_{c(\n)}^{(i)}(\infty), \qquad i\in\{1,2\}.
\]
\end{Remark}

\begin{Remark}
Since the products of all the branches of $ S_c(\z) $ as well as of $ \Phi_\n(\z) $ are identically equal to $ 1 $, we immediately deduce from \eqref{main1} that
\[
P_\n(z) = (1+ o(1)) \mathcal P_\n(z), \qquad
P_\n(x) = (1+ o(1)) \mathcal P_{\n+}(x) + (1+ o(1)) \mathcal P_{\n-}(x),
\]
uniformly for $ \dist(z,\Delta_{c(\n)})\geq d $, $ z\in\overline\C $, and $ \dist(x,E_{c(\n)})\geq d $, $ x\in\Delta_{c(\n)} $, respectively, for any~${ d>0} $ fixed, where
\begin{equation}
\label{Fn}
\mathcal P_\n(z) := \frac{(S_{c(\n)}\Phi_\n)^{(0)}(z)}{\tau_\n S_{c(\n)}^{(0)}(\infty)} = \mathcal P_{\n,1}(z)\mathcal P_{\n,2}(z) .
\end{equation}
\end{Remark}

\begin{Remark}
It might seem that asymptotic formulae \eqref{main1} do not significantly reduce complexity as the functions on both sides of the equalities depend on $ \n $. In this regard we would like to stress that the Szeg\H{o} functions only depend on the one-dimensional parameter $ c(\n)=n_1/|\n| $ rather than on the two-dimensional multi-index $ \n $, while the geometric factors must depend on~$ \n $ to properly match the behavior at infinity of $ P_{\n,1}(z) $ and $ P_{\n,2}(z) $, yet their absolute values satisfy \eqref{Phin3}, where the measures $ \omega_{c(\n),1} $ and $ \omega_{c(\n),2} $ again depend only on~$ c(\n) $.
\end{Remark}

\begin{Theorem}
\label{thm:2}
Under the conditions of Theorem~{\rm\ref{thm:1}}, let $ A_\n^{(1)}(z) $ and $ A_\n^{(2)}(z) $ be the type~I multiple orthogonal polynomials defined via \eqref{typeI} and \eqref{Qn_norm}. Define
\[
\mathcal A_{\n,i}(z) := \tau_\n S_{c(\n)}^{(0)}(\infty) w_{c(\n),i}(z) \left(\frac{-\Pi_{c(\n)}}{S_{c(\n)}\Phi_\n}\right)^{(i)}(z).
\]
Understanding the error terms $ o(1) $ exactly as in Theorem~{\rm\ref{thm:1}}, we have that
\begin{gather}
A_\n^{(1)}(z) = \left(1+ \frac{o(1)}{c(\n)}\right) \mathcal A_{\n,1}(z), \qquad
A_\n^{(2)}(z) = \left(1+ \frac{o(1)}{1-c(\n)}\right) \mathcal A_{\n,2}(z),\nonumber \\
A_\n^{(i)}(x) = (1+o(1)) \mathcal A_{\n,i+}(x) + (1+o(1)) \mathcal A_{\n,i-}(x),\label{main2}
\end{gather}
uniformly for $ \dist(z,\Delta_{c(\n),1})\geq d $, $ z\in\overline\C $, $ \dist(z,\Delta_{c(\n),2})\geq d $, $ z\in\overline\C $, and $ \dist(x,E_{c(\n),i})\geq d $, $ x\in\Delta_{c(\n),i} $, $ i\in\{1,2\} $, respectively, for any $ d>0 $ fixed.
\end{Theorem}

Theorem~\ref{thm:2} subsumes \cite[Theorem~3.3]{ApDenYa21} where the error terms were not uniform.

\begin{Remark}
It is known, see \eqref{c-rate} further below, that the length of $ \Delta_{\n,1} $ is proportional to~$ c(\n) $ while the length of $ \Delta_{\n,2} $ is proportional to $ 1 - c(\n) $. Therefore, the bottom formula in \eqref{main2} is meaningful only when $ c(\n) $ is separated from both $ 0 $ and $ 1 $ and thus there is no need to divide the error factors by $ c(\n) $ or $ 1- c(\n) $.
\end{Remark}

\begin{Remark}
It readily follows from the top formulae in \eqref{main2} that $ A_\n^{(i)}(z) $ has exactly $ n_i $ zeros which belong to $ \dist(z,\Delta_{c(\n),i})\leq d $ for all $ \varepsilon_\n $ small enough and $ c(\n) $ separated from $ 0 $ and~$ 1 $.
\end{Remark}

\begin{Remark}
Recall that we denoted by $ 1/h_{\n-\vec e_i,i} $ the leading coefficient of $ A_\n^{(i)}(z) $, $ i\in\{1,2\} $. We readily get from the sentence after \eqref{Pin}, \eqref{Phin}, and \eqref{main2} that
\[
h_{\n-\vec e_1,1} = \left(1+\frac{o(1)}{c(\n)}\right)A_{c(\n),1}^{n_1-1} \bigl(B_{c(\n),1} - B_{c(\n),2}\bigr)^{n_2} S_{c(\n)}^{(1)}(\infty)/S_{c(\n)}^{(0)}(\infty)
\]
and an analogous formula holds for $ h_{\n-\vec e_2,2} $. In fact, we can deduce that these constants are non-zero for all $ \varepsilon_\n $ small regardless of $ c(\n) $ being close to $ 0 $ or $ 1 $ because $ h_{\n,i}=\int P_\n P_{\n,i}{\rm d}\mu_i $ by~\eqref{connection}, which must be non-zero by Theorem~\ref{thm:1} ($ \deg P_\n = |\n| $ for all $ \varepsilon_\n $ small means there are no extra orthogonality conditions and hence $ h_{\n,i}\neq 0 $).
\end{Remark}

\begin{Theorem}
\label{thm:3}
In the setting of Theorem~{\rm\ref{thm:1}}, let the coefficients $ a_{\n,i},b_{\n,i} $ be as in \eqref{rec_rel} and the numbers $ A_{c(\n),i},B_{c(\n),i} $ as in \eqref{AngPar1}. Then it holds that
\begin{equation}
\label{main3}
a_{\n,i} = A_{c(\n),i} + o(1) \qquad \text{and}\qquad  b_{\n,i} = B_{c(\n+\vec e_i),i} + o(1),
\end{equation}
$ i\in\{1,2\} $, where the error terms $ o(1) $ should be understood exactly as in Theorem~{\rm\ref{thm:1}}.
\end{Theorem}

Again, Theorem~\ref{thm:3} generalizes \cite[Theorem~1.2]{ApDenYa21} by proving uniformity of the error terms.

The map $ \chi_c(\z) $ takes $ \boldsymbol\Delta_c $ onto two Jordan curves. It was shown in \cite[Lemma~4.1.2]{DenYa22} that these curves can be parametrized as
\begin{equation}
\label{chi_curves}
\chi_c(\boldsymbol\Delta_c) = \left\{ \chi\in\C : \frac{A_{c,1}}{|\chi-B_{c,1}|^2} + \frac{A_{c,2}}{|\chi-B_{c,2}|^2} = 1 \right\}.
\end{equation}
It follows from \eqref{hc1} and \eqref{point_zc} that $ \chi(\z_c) $ must belong to $ \chi_c(\boldsymbol\Delta_c) $ when $ c\in(0,c^*]\cup[c^{**},1) $. It must also hold that $ \chi(\z_c)=(1-c)B_{c,1} + cB_{c,2} $ as one can see from \eqref{Pin}, \eqref{h_chi}, and \eqref{hc1}. Therefore, it necessarily holds that
\begin{equation}
\label{AcBc}
c^{-2}A_{c,1} + (1-c)^{-2}A_{c,2} = B_c^2,
\end{equation}
$ c\in(0,c^*]\cup[c^{**},1) $, where we set $ B_c := B_{c,2}-B_{c,1} $. In particular, relations \eqref{AcBc} together with~\eqref{AngPar2} yield that
\begin{gather}
\lim_{c\to0} c^{-2} A_{c,1} = \big[(\beta_2-\alpha_2)/4\big]^2 + \varphi_2^2(\alpha_1),\nonumber \\
\lim_{c\to1} (1-c)^{-2} A_{c,2} = \big[(\beta_1-\alpha_1)/4\big]^2 + \varphi_1^2(\beta_2).\label{c_vanish}
\end{gather}

The recurrence coefficients $ a_{\n,i}$, $b_{\n,i} $ must satisfy what is known as \emph{compatibility conditions}, which are a system of discrete difference equations, see \cite[Theorem~3.2]{VA11}. This suggests that~$ A_{c,i}$,~$B_{c,i} $, as functions of the parameter $ c $, must satisfy a system of differential equation. This was conditionally confirmed in \cite[Theorem~3]{ApKoz20} for Angelesco systems of any size. The conditional part came from the requirement on the speed of convergence of the recurrence coefficients to their limits which is beyond of what is currently has been demonstrated including Theorem~\ref{thm:3} above. As it happens, we are able to show the validity of these differential equation using only our asymptotic analysis.

\begin{Theorem}
\label{thm:4}
Set $ R(c) := (c/(1-c))^2 (A_{c,2}/A_{c,1}) $, which is a continuous function on $ [0,1] $. It holds that
\begin{equation}
\label{main5}
R^\prime(c) = \frac{6R(c)(1+R(c))}{1-c^2+c(2-c)R(c)}
\end{equation}
on $ (0,c^*)\cup(c^{**},1) $. Moreover, we have that
\begin{gather}
 \frac{B_c^\prime}{B_c} = -\frac c{1-c} \frac{A_{c,1}^\prime}{A_{c,1}} = -\frac{1-c}c \frac{A_{c,2}^\prime}{A_{c,2}} = -2\frac{1-c-cR(c)}{1-c^2+c(2-c)R(c)},\nonumber \\
c B_{c,1}^\prime + (1-c) B_{c,2}^\prime = 0 \ \Leftrightarrow \ B_{c,2}^\prime = c B_c^\prime \ \Leftrightarrow \ B_{c,1}^\prime = -(1-c) B_c^\prime,\label{main4}
\end{gather}
on $ (0,c^*)\cup(c^{**},1) $, where $ ^\prime $ indicates the derivative with respect to the parameter $ c $.
\end{Theorem}

\begin{Remark}
Equation \eqref{main5} together with the initial conditions coming from \eqref{AngPar2} and \eqref{c_vanish} allows us to reconstruct $ R(c) $ uniquely on $ [0,c^*]\cup[c^{**},1] $. Since~$ R(c)>0 $, \eqref{main5} also shows that~$ R(c) $ is infinitely differentiable on $ (0,c^*)\cup(c^{**},1) $ and all the derivatives extend continuously to $ [0,c^*]\cup[c^{**},1] $. The first line of \eqref{main4} now allows one to recover $ A_{c,1} $, $ A_{c,2} $ \big(to remove singularities, it is better two rewrite these equations for $ c^{-2}A_{c,1} $ and $ (1-c)^{-2} A_{c,2} $\big), and~$ B_c $ as well as to draw the same conclusions about infinite differentiability and continuity; $ B_{c,1} $ and $ B_{c,2} $ are then recovered via the second line of \eqref{main4}.
\end{Remark}

To prove Theorems~\ref{thm:1}--\ref{thm:4} we use the extension to multiple orthogonal polynomials \cite{GerKVA01} of by now classical approach of Fokas, Its, and Kitaev \cite{FIK91,FIK92} connecting orthogonal polynomials to matrix Riemann--Hilbert problems. The RH problem is then analyzed via the non-linear steepest descent method of Deift and Zhou \cite{DZ93}.

\section{Model local parametrices}
\label{sec:MLP}

In this section, we formulate several Riemann--Hilbert problems with constant jumps for $ 2\times2 $ matrices that will be used in the main part of the proof. In what follows, the symbol $ \boldsymbol I $ stands for the identity matrix of any size, $ \sigma_3 := \diag(1,-1) $ is the third Pauli matrix, and we let
\begin{equation}
\label{K}
\boldsymbol K(\zeta) := \frac{\zeta^{-\sigma_3/4}}{\sqrt2}\left(\begin{matrix} 1 & \ic \\ \ic & 1 \end{matrix}\right),
\end{equation}
where the root is principal, i.e., $\arg(\zeta)\in(-\pi,\pi)$ (a convention we follow for all the power functions unless explicitly specified otherwise). Further, for brevity, we denote the rays $ \{\arg(z)=\pm2\pi/3\} $ by $ I_\pm $ and orient them towards the origin.

 \subsection{Hard edge}

 Let $\boldsymbol\Psi(\zeta)$ be a matrix-valued function such that
\begin{itemize}\itemsep=0pt
\label{rhpsi}
\item[(a)] $\boldsymbol\Psi(\zeta)$ is holomorphic in $\C\setminus\bigl(I_+\cup I_-\cup(-\infty,0]\bigr)$;
\item[(b)] $\boldsymbol\Psi(\zeta)$ has continuous boundary values on $I_+\cup I_-\cup(-\infty,0)$ that satisfy
\[
\boldsymbol\Psi_+(\zeta) = \boldsymbol\Psi_-(\zeta)
\begin{cases}
\left(\begin{matrix} 0 & 1 \\ -1 & 0 \end{matrix}\right), & \zeta\in(-\infty,0), \vspace{1mm}\\
\left(\begin{matrix} 1 & 0 \\ 1 & 1 \end{matrix}\right), & \zeta\in I_\pm;
\end{cases}
\]
\item[(c)] $ \boldsymbol\Psi(\zeta) = \boldsymbol{\mathcal O}(\log\zeta) $ as $\zeta\to0$, where $ \boldsymbol{\mathcal O}(\cdot) $ is understood entrywise;
\item[(d)] it holds uniformly for $ |\zeta| $ large that
\[
\boldsymbol\Psi(\zeta) = \boldsymbol K(\zeta)\bigl(\boldsymbol I + \boldsymbol{\mathcal O} \bigl(\zeta^{-1/2}\bigr)\bigr)\exp\bigl\{2\zeta^{1/2}\sigma_3\bigr\}.
\]
\end{itemize}
The solution of \hyperref[rhpsi]{\rhpsi} was constructed explicitly in \cite[Section~6]{KMcLVAV04} with the help of the modified Bessel and Hankel functions. Since the jump matrices in \hyperref[rhpsi]{\rhpsi}(b) have determinant one, $ \det \boldsymbol\Psi(\zeta) $ is analytic in $ \C\setminus\{0\} $. It then follows from \hyperref[rhpsi]{\rhpsi}(c,d) and \eqref{K} that $ \det \boldsymbol\Psi(\zeta) \equiv \sqrt 2 $.

 Set $ \boldsymbol \Psi_*(\zeta) := \sigma_3\boldsymbol\Psi(\zeta)\sigma_3 $. Then $ \boldsymbol \Psi_*(\zeta) $ solves the following Riemann--Hilbert problem:
 \begin{itemize}\itemsep=0pt\setlength{\leftskip}{0.1cm}
\label{rhwpsi}
\item[(a--d)] $\boldsymbol\Psi_*(\zeta)$ satisfies \hyperref[rhpsi]{\rhpsi}(a--d), but with the reverse orientation of the rays in \hyperref[rhpsi]{\rhpsi}(b) and $ \boldsymbol K(\zeta) $ replaced by $ \sigma_3\boldsymbol K(\zeta)\sigma_3 $ in \hyperref[rhpsi]{\rhpsi}(d).
\end{itemize}

\subsection{Sliding soft edge}
\label{ss:sse}

Let $ \tau\in[\tau_*,\infty] $ for some $ \tau_*>1 $ to be fixed later. If $ \tau<\infty $, define $ U_\tau $ to be the disk of unit radius centered at $ \tau $ and orient $ \partial U_\tau $ clockwise. Denote by $\boldsymbol\Theta(\zeta;\tau)$ the solution, if it exists, of the following Riemann--Hilbert problem (\rhth):
\begin{itemize}\itemsep=0pt
\label{rhth}
\item[(a)] $\boldsymbol\Theta(\zeta;\tau)$ is holomorphic in $\C\setminus\bigl(I_+\cup I_-\cup(-\infty,\tau]\bigr)$;
\item[(b)] $\boldsymbol\Theta(\zeta;\tau)$ has continuous boundary values on $I_+\cup I_-\cup(-\infty,0)\cup(0,\tau)$ that satisfy \hyperref[rhpsi]{\rhpsi}(b) and
\[
\boldsymbol\Theta_+(\zeta;\tau) = \boldsymbol\Theta_-(\zeta;\tau) \left(\begin{matrix} 1 & 1 \\ 0 & 1 \end{matrix}\right), \qquad \zeta\in(0,\tau);
\]
\item[(c)] $ \boldsymbol\Theta(\zeta;\tau)=\boldsymbol{\mathcal O}(1) $ as $\zeta\to0$ and $ \boldsymbol\Theta(\zeta;\tau)=\boldsymbol{\mathcal O}(\log|\zeta-\tau|) $ as $\zeta\to\tau$ when $ \tau $ is finite;
\item[(d)] it holds uniformly for $ \tau\geq\tau_* $ and $ |\zeta| $ large and such that $ \zeta\not\in U_\tau $ that
\[
\boldsymbol\Theta(\zeta;\tau) = \boldsymbol K(\zeta)\left(\boldsymbol I+\boldsymbol{\mathcal O}\left(\frac1{\sqrt{|\zeta}|\min\{\tau,|\zeta|\}}\right)\right) \exp\left\{-\frac23\zeta^{3/2}\sigma_3\right\}.
\]
\end{itemize}

The solution $ \boldsymbol\Theta_{\rm Ai}(\zeta) := \boldsymbol\Theta(\zeta;\infty) $ of \hyperref[rhth]{\rhth} for $ \tau=\infty $ is well known \cite{DKMLVZ99b} and is explicitly constructed using Airy functions. As in the previous subsection, $ \det \boldsymbol\Theta(\zeta;\tau) \equiv \sqrt 2 $.

To show that \hyperref[rhth]{\rhth} is also solvable for finite $ \tau\geq\tau_* $ and some $ \tau_*>1 $, define
\[
\boldsymbol\Theta_\tau(\zeta) := \boldsymbol K(\zeta) {\rm e}^{-(2/3)\zeta^{3/2}\sigma_3} \begin{pmatrix} 1 & l_\tau(\zeta) \\ 0 & 1 \end{pmatrix}, \qquad l_\tau(\zeta) := \frac1{2\pi\ic} \log(\zeta-\tau),
\]
for $ \zeta\in U_\tau\setminus(\tau-1,\tau) $, where we take the principal branch of the logarithm. Since $ U_\tau $ belongs to the right half-plane, the matrix $ \boldsymbol\Theta_\tau(\zeta) $ is analytic in the domain of its definition and satisfies
\[
\boldsymbol\Theta_{\tau+}(\zeta) = \boldsymbol\Theta_{\tau-}(\zeta) \begin{pmatrix} 1 & 1 \\ 0 & 1 \end{pmatrix}, \qquad \zeta\in(\tau-1,\tau).
\]
That is, $ \boldsymbol\Theta_\tau(\zeta) $ solves \hyperref[rhth]{\rhth} locally in $ U_\tau $. Consider the following Riemann--Hilbert problem: find a matrix function $ \boldsymbol R(\zeta;\tau) $ such that
\begin{itemize}\itemsep=0pt
\item[(a)] $\boldsymbol R(\zeta;\tau)$ is holomorphic in $ \C\setminus \bigl(\partial U_\tau \cup (\tau+1,\infty)\bigr) $ and $ \boldsymbol R(\zeta;\tau) = \boldsymbol I + \boldsymbol{\mathcal O}\bigl(\zeta^{-1}\bigr) $ as $ \zeta\to\infty $;
\item[(b)] $\boldsymbol R(\zeta;\tau)$ has continuous and bounded boundary values on $\partial U_\tau\setminus\{\tau+1\}$ and $ (\tau+1,\infty) $ that satisfy
\[
\boldsymbol R_+(\zeta;\tau) = \boldsymbol R_-(\zeta;\tau) \begin{cases}
\boldsymbol\Theta_\tau(\zeta)\boldsymbol\Theta_{\rm Ai}^{-1}(\zeta), & \zeta\in \partial U_\tau\setminus\{\tau+1\}, \\
\boldsymbol\Theta_{{\rm Ai}-}(\zeta)\boldsymbol\Theta_{{\rm Ai}+}^{-1}(\zeta), & \zeta\in(\tau+1,\infty).
\end{cases}
\]
\end{itemize}

By using the definition of $ \boldsymbol\Theta_\tau(\zeta) $ as well as \hyperref[rhth]{\rhth}(d) with $ \tau=\infty $, one can readily check that the jump of $ \boldsymbol R(\zeta;\tau) $ on $ \partial U_\tau $ can be estimated as
\[
\left(\boldsymbol I + \boldsymbol K(\zeta) \begin{pmatrix} 0 & l_\tau(\zeta){\rm e}^{-(4/3)\zeta^{3/2}} \\ 0 & 0 \end{pmatrix} \boldsymbol K^{-1}(\zeta) \right) \bigl( \boldsymbol I + \boldsymbol{\mathcal O}\bigl(\zeta^{-1}\bigr) \bigr) = \boldsymbol I + \boldsymbol{\mathcal O}\bigl(\tau^{-1}\bigr),
\]
where the error term is uniform in $ \tau $. Similarly, by using \hyperref[rhth]{\rhth}(b,d) with $ \tau=\infty $ we get that the jump of $ \boldsymbol R(\zeta;\tau) $ on $ (\tau+1,\infty) $ can be estimated as
\[
\boldsymbol I - \boldsymbol\Theta_{{\rm Ai}-}(\zeta) \begin{pmatrix} 0 & 1 \\ 0 & 0 \end{pmatrix} \boldsymbol\Theta_{{\rm Ai}-}^{-1}(\zeta) = \boldsymbol I + \boldsymbol{\mathcal O}\bigl(\sqrt\zeta {\rm e}^{-(4/3)\zeta^{3/2}} \bigr) = \boldsymbol I + \boldsymbol{\mathcal O}\bigl(\tau^{-1}\bigr),
\]
where again the estimate is uniform in $ \tau $. Therefore, we can conclude from \cite[Theorem~8.1]{FokasItsKapaevNovokshenov} that $ \boldsymbol R(\zeta;\tau) $ exists for all $ \tau\geq\tau_* $ and some $ \tau_*>1 $ and satisfies
\[
\boldsymbol R(\zeta;\tau) = \boldsymbol I + \boldsymbol{\mathcal O}\left( \tau^{-1}(1+|\zeta|)^{-1}\right)
\]
uniformly for all $ \zeta\in \C $ and $ \tau\geq\tau_* $, that is, including the boundary values (uniformity of the estimate in $ \zeta $ is achieved by varying the jump contour slightly, which is possible due to analyticity of the jump matrices). Now, it only remains to observe that \hyperref[rhth]{\rhth} is solved by
\[
\boldsymbol \Theta(\zeta;\tau) = \boldsymbol R(\zeta;\tau)
\begin{cases}
\boldsymbol\Theta_\tau(\zeta), & \zeta\in U_\tau, \\
\boldsymbol\Theta_{\rm Ai}(\zeta), & \zeta\in\C\setminus\overline U_\tau.
\end{cases}
\]

\subsection{Critical soft edge}

Given $ s\in(-\infty,\infty) $, let $\boldsymbol\Phi(\zeta;s)$ be such that
\begin{itemize}\itemsep=0pt\setlength{\leftskip}{0.1cm}
\label{rhphi}
\item[(a--c)] $\boldsymbol\Phi(\zeta;s)$ satisfies \hyperref[rhpsi]{\rhpsi}(a--c);
\item[(d)] it holds uniformly for $ |\zeta| $ large and locally uniformly in $ s $ that
\[
\boldsymbol\Phi(\zeta;s) = \boldsymbol K(\zeta) \bigl(\boldsymbol I + \boldsymbol{\mathcal O}\bigl(\zeta^{-1/2}\bigl) \bigl) \exp\left\{-\frac23(\zeta+s)^{3/2}\sigma_3\right\}.
\]
\end{itemize}

 The solvability of this problem was obtained in \cite{XuZh11}. The fact that $ \boldsymbol{\mathcal O}(\cdot)$ is locally uniform in~$ s $ was pointed out in \cite{IKOs08}. Again, observe that $ \det \boldsymbol\Phi(\zeta;s) \equiv \sqrt 2 $.

Further, given $ s\leq0 $, consider a similar Riemann--Hilbert problem \big(\rhwphi\big):
\begin{itemize}\itemsep=0pt\setlength{\leftskip}{0.1cm}
\label{rhwphi}
\item[(a--c)] $\widetilde{\boldsymbol\Phi}(\zeta;s)$ satisfies \hyperref[rhpsi]{\rhpsi}(a--c);
\item[(d)] it holds uniformly for all $ s \leq 0 $ and locally uniformly for $ \zeta/(1-s) \in\overline\C\setminus\{0\} $ that
\[
\widetilde{\boldsymbol\Phi}(\zeta;s) = \boldsymbol K(\zeta) \bigl(\boldsymbol I + \boldsymbol{\mathcal O}\bigl((1-s)^{-1}\zeta^{-1/2}\bigr) \bigr) \exp\left\{-\frac23\bigl(\zeta^{3/2}+s\zeta^{1/2}\bigr)\sigma_3\right\}.
\]
\end{itemize}

It was observed in \cite{Ya16} that solvability of \hyperref[rhphi]{\rhphi} for $ s\leq0 $ is equivalent to solvability of \hyperref[rhwphi]{\rhwphi} for $ s\leq0 $. It was also stated in \cite[equation~(4.3)]{Ya16} that the error term in \hyperref[rhwphi]{\rhwphi}(d) behaves like $ \boldsymbol{\mathcal O}\bigl((1-s)^{1/2}\zeta^{-1/2}\bigr) $ uniformly for all $ s\leq0 $. Below, we show how the bound from~\cite[equation~(4.3)]{Ya16} can be improved to the one stated in \smash{\hyperref[rhwphi]{\rhwphi}(d)}.

It has been already mentioned that the desired bound in \hyperref[rhwphi]{\rhwphi}(d) must hold locally uniformly for $ s\leq0 $, see \cite{IKOs08}. Therefore, we are only interested in what happens for $ -s $ large enough. To this end, let
$
 g(\xi) := (2/3)(\xi-1)\xi^{1/2}
$
be the principal branch holomorphic in $ \C\setminus(-\infty,0] $. Given $ \kappa\geq\kappa_0>0 $ for some $ \kappa_0 $ large enough, consider the following Riemann--Hilbert problem:
\begin{itemize}\itemsep=0pt\setlength{\leftskip}{0.1cm}
\label{rhhphi}
\item[(a--c)] $\widehat{\boldsymbol\Phi}(\xi;\kappa)$ satisfies \hyperref[rhpsi]{\rhpsi}(a--c);
\item[(d)] it holds uniformly for $ \kappa\geq\kappa_0 $ and locally uniformly for $ \xi\in\overline\C\setminus\{0\} $ that
\[
\widehat{\boldsymbol\Phi}(\xi;\kappa) = \boldsymbol K(\xi)\bigl(\boldsymbol I + \boldsymbol{\mathcal O}\bigl(\kappa^{-1}\xi^{-1/2}\bigr) \bigr) {\rm e}^{-\kappa g(\xi)\sigma_3}.
\]
\end{itemize}

We shall show that there exists $ \kappa_0 >0 $ such that \hyperref[rhhphi]{\rhhphi} is uniquely solvable. In this case it can be readily verified that the solution of \hyperref[rhhphi]{\rhhphi} yields the solution of \hyperref[rhwphi]{\rhwphi} via
\[
\widetilde{\boldsymbol\Phi}(\zeta;s) = (-s)^{-\sigma_3/4}\widehat{\boldsymbol\Phi}(-\zeta/s;\kappa), \qquad \kappa = (-s)^{3/2}.
\]

Let $ U_0 $ be a disk centered at the origin of any radius $ r_0<1 $ small enough so that $ g^2(\xi) $ is conformal in it. Set $ I_\pm^* := \bigl(g^2\bigr)^{-1}(I_\pm)\cap U_0 $ and orient these arcs towards the origin. Notice that the principal square root branch of $ g^2(\xi) $ is equal to $ -g(\xi) $ and that
\[
-g(I_\pm\cap U_0) \subset \bigl\{ \sqrt x\bigl(2x+1 \pm \sqrt 3\ic\bigr)/3 : x>0\bigr\},
\]
which are arcs that lie within the sector $ |\arg(\zeta)|<\pi/3 $. As $ g^2(\xi) $ preserves the negative and positive reals, the curve $ I_+^* $ lies between the rays $ I_+ $ and $ (-\infty,0) $ (since $ g^2(I_+^*) = I_+ $ lies between~$ g^2(I_+) $ and the negative reals) and the curve $ I_-^* $ lies between the rays $ I_- $ and $ (-\infty,0) $.

Let $ \boldsymbol\Psi(\zeta) $ be the solution of \hyperref[rhpsi]{\rhpsi}. For $ \xi \in U_0 $, define
\[
 \widehat{\boldsymbol\Phi}_0(\xi;\kappa) := \boldsymbol K(\xi) \bigl(\boldsymbol K^{-1} \boldsymbol\Psi\bigr)\bigl( (\kappa g(\xi)/2 )^2 \bigr)
\begin{cases}
\begin{pmatrix} 1 & 0 \\ \pm 1 & 1 \end{pmatrix}, & \xi\in S_\pm^*, \\
\boldsymbol I, & \text{otherwise},
\end{cases}
\]
where $ S_+^* $ and $ S_-^* $ are the sectors within $ U_0 $ delimited by $ I_+ $ and $ I_+^* $ in the second quadrant, and $ I_- $ and $ I_-^* $ in the third quadrant, respectively (notice that the sectors $ S_\pm^* $ do not depend on $ \kappa $ because the preimages of $ I_\pm $ under $ g^2(\xi) $ and $ (\kappa g(\xi)/2)^2 $ must coincide). One can readily verify that the matrix $ \widehat{\boldsymbol\Phi}_0(\xi;\kappa) $ satisfies \hyperref[rhhphi]{\rhhphi}(a--c) within $ U_0 $. Moreover, since the domains~$ - \kappa g(S_\pm^*)/2 $ lie with the sector $ |\arg(\zeta)|\leq\pi/3 $, it holds that
\[
\big| {\rm e}^{2\kappa g(\xi)} \big| = {\rm e}^{-2\kappa|\re g(\xi)|} < \frac1{2\kappa |\re g(\xi)|} \leq \frac1{\kappa|g(\xi)|}, \qquad \xi\in S_\pm^*,
\]
and therefore
\[
{\rm e}^{-\kappa g(\xi)\sigma_3} \begin{pmatrix} 1 & 0 \\ \pm1 & 1 \end{pmatrix} = \left( \boldsymbol I + \boldsymbol{\mathcal O}\left(\frac1{\kappa g(\xi)}\right) \right) {\rm e}^{-\kappa g(\xi)\sigma_3},
\]
uniformly for $ \xi\in S_\pm^* $ and all $ \kappa>0 $. Since $ |g(\xi)| $ is uniformly bounded away from $ 0 $ on $ \partial U_0 $, the last estimate and \hyperref[rhpsi]{\rhpsi}(d) yield that
\[
\widehat{\boldsymbol\Phi}_0(\xi;\kappa) = \boldsymbol K(\xi) \bigl(\boldsymbol I + \boldsymbol{\mathcal O}\bigl(\kappa^{-1}\bigr) \bigr) {\rm e}^{-\kappa g(\xi)\sigma_3}
\]
uniformly for $ \xi\in\partial U_0 $.

Orient $ \partial U_0 $ clockwise. Consider the following Riemann--Hilbert problem: find $ \widehat{\boldsymbol R}(\xi;\kappa) $ such that
\begin{itemize}\itemsep=0pt
\item[(a)] $\widehat{\boldsymbol R}(\xi;\kappa)$ is holomorphic in $ \overline \C $ away from $ \partial U_0 \cup ( (I_+\cup I_-) \setminus U_0) $ and $ \widehat{\boldsymbol R}(\xi;\kappa) = \boldsymbol I + \boldsymbol{\mathcal O}\bigl(\xi^{-1}\bigr) $ as $ \xi\to\infty $;
\item[(b)] $\widehat{\boldsymbol R}(\xi;\kappa)$ has continuous and bounded boundary values on $\partial U_0\setminus(I_+\cup I_-)$ and $ (I_+\cup I_-) \setminus \overline U_0 $ that satisfy
\[
\widehat{\boldsymbol R}_+(\xi;\kappa) = \widehat{\boldsymbol R}_-(\xi;\kappa)
\begin{cases}
\widehat{\boldsymbol\Phi}_0(\xi;\kappa) {\rm e}^{\kappa g(\xi)\sigma_3} \boldsymbol K^{-1}(\xi), & \xi\in \partial U_0\setminus(I_+\cup I_-), \\
\boldsymbol K(\xi) \begin{pmatrix} 1 & 0 \\ {\rm e}^{2\kappa g(\xi)} & 1 \end{pmatrix} \boldsymbol K^{-1}(\xi), & \xi\in(I_+\cup I_-) \setminus \overline U_0.
\end{cases}
\]
\end{itemize}

As we have already observed, $ 3\re(g(\xi))=-2|\xi|^{3/2}-|\xi|^{1/2} $ for $ \xi\in I_\pm $. Thus,
\[
\widehat{\boldsymbol R}_+(\xi;\kappa) = \widehat{\boldsymbol R}_-(\xi;\kappa) \bigl( \boldsymbol I + \boldsymbol{\mathcal O} \bigl(\kappa^{-1} \bigr) \bigr)
\]
uniformly on $ \partial U_0 \cup ( (I_+\cup I_-) \setminus U_0) $ (with respect to both $ \xi $ and $ \kappa $). Hence, as in the previous subsection, we can conclude on the basis of \cite[Theorem~8.1]{FokasItsKapaevNovokshenov} and the deformation of the contour technique that $ \widehat{\boldsymbol R}(\xi;\kappa) $ does indeed uniquely exist for all $ \kappa\geq\kappa_0 $ and some $ \kappa_0 $ large enough and satisfies
\smash{$
\widehat{\boldsymbol R}(\xi;\kappa) = \boldsymbol I + \boldsymbol{\mathcal O} \bigl( \kappa^{-1}(1+|\xi|)^{-1} \bigr)
$}
uniformly for all $ \xi\in \C $ and $ \kappa\geq\kappa_0 $. It remains to~observe that the solution of \hyperref[rhhphi]{\rhhphi} is given~by
\[
\widehat{\boldsymbol\Phi}(\xi;\kappa) := \widehat{\boldsymbol R}(\xi;\kappa)
\begin{cases}
\widehat{\boldsymbol\Phi}_0(\xi;\kappa), & \xi\in U_0, \\
\boldsymbol K(\xi){\rm e}^{-\kappa g(\xi)\sigma_3}, & \xi\in\C\setminus \overline U_0.
\end{cases}
\]

\section{Conformal maps}

The framework of the Riemann--Hilbert analysis, which we use, consists in formulating a multiplicative Riemann--Hilbert problem for $ 3\times3 $ matrices, whose jump relations are then factorized and partially moved into the complex plane onto the so-called ``lens''. The construction of this lens depends on the value of the parameter $ c $ via properties of local conformal maps around each point in $ E_c $. In this section, we define these maps and discuss some of their properties.

In what follows, it will be sometimes useful for us to recall that
\begin{equation}
\label{c-rate}
\lim_{c\to0} \frac{|\Delta_{c,1}|}c = 4|w_2(\alpha_1)| \qquad \text{and}\qquad  \lim_{c\to1} \frac{|\Delta_{c,2}|}{1-c} = 4|w_1(\beta_2)|,
\end{equation}
which was shown in \cite[equation~(4.8)]{ApDenYa21}, where the roots $ w_i(z) $ were introduced in \eqref{phi_w}. We shall also use the following well known fact: according to Koebe's 1/4 theorem, if $ \zeta(z) $ is conformal in $ \{|z-e|\leq r_*\} $ with $ \zeta(e)=0 $, then
\begin{equation}
\label{koebe}
(r/4)|\zeta^\prime(e)| \leq |\zeta(z)|, \qquad |z-e|=r, \qquad r\leq r_*.
\end{equation}

\subsection{Conformal maps}

The material of this section is taken from \cite[Section~7.4]{ApDenYa21}. We work only with the interval $ \Delta_{c,1} $, the maps around $ \Delta_{c,2} $ are constructed similarly. Given $ c\in(0,1) $, define
\begin{equation}
\label{4.1.1}
\zeta_{c,\alpha_1}(z) := \left(\frac14\int_{\alpha_1}^z\bigl(h_c^{(0)}-h_c^{(1)}\bigr)(s) {\rm d}s\right)^2, \qquad \re z<\beta_{c,1},
\end{equation}
where $ h_c(\z) $ was defined in \eqref{h_chi}. Then the following lemma holds, see \cite[Lemma~7.4]{ApDenYa21}.

\begin{Lemma}
\label{lem:4.1}
For each $ c\in(0,1) $, the map $ \zeta_{c,\alpha_1}(z) $ is positive on~$ (-\infty,\alpha_1) $ and negative on~$ (\alpha_1,\beta_{c,1}) $ with a simple zero at $ \alpha_1 $. Moreover, there exist constants $ \delta_{\alpha_1}>0 $ and $ A_{\alpha_1}>0 $, independent of $ c $, such that $ \zeta_{c,\alpha_1}(z) $ is conformal in $ \{|z-\alpha_1| \leq \delta_{\alpha_1}c \} $ and satisfies $ 4cA_{\alpha_1} \leq |\zeta^\prime_{c,\alpha_1}(\alpha_1)| $.
\end{Lemma}

As already apparent from \eqref{4.1.1}, the function $ h_c(\z) $ plays the central role in this subsection. Recall the special point $ \z_c $, see \eqref{hc1}, and its relation to $ \boldsymbol\beta_{c,1} $, see \eqref{cstars} and \eqref{point_zc}. Hence, while constructing conformal maps at $ \beta_{c,1} $, we need to consider several cases.

Given $ c\in (0,c^*] $, in which case $ \boldsymbol\beta_{c,1}=\z_c $ and $ h_c(\boldsymbol\beta_{c,1})$ is finite, define
\begin{equation}
\label{4.2.1}
\zeta_{\beta_{c,1}}(z) := \left(-\frac34\int_{\beta_{c,1}}^z\bigl(h_c^{(0)}-h_c^{(1)}\bigr)(s) {\rm d}s\right)^{2/3}, \qquad \alpha_1<\re z<\alpha_2,
\end{equation}
where the choice of the root is made so that $ \zeta_{\beta_{c,1}}(z) $ is positive for $ x>\beta_{c,1} $. Then the following lemma holds, see \cite[Lemma~7.5]{ApDenYa21}.

\begin{Lemma}
\label{lem:4.2}
For each $ c\in(0,c^*] $, the map $ \zeta_{\beta_{c,1}}(z) $ is positive on $ (\beta_{c,1},\alpha_2) $ and negative on $ (\alpha_1,\beta_{c,1}) $ with a simple zero at $ \beta_{c,1} $. Moreover, there exist constants $ \delta_{\beta_1}>0 $ and $ A_{\beta_1}>0 $, independent of $ c $, such that $ \zeta_{\beta_{c,1}}(z) $ is conformal in $ \{|z-\beta_{c,1}| \leq \delta_{\beta_1}c \} $ and satisfies $ 4c^{-1/3} A_{\beta_1} \leq \zeta_{\beta_{c,1}}^\prime(\beta_{c,1}) $.
\end{Lemma}

When $ c>c^* $, we can and do define a conformal around $ \beta_1 $ similarly to \eqref{4.1.1}, see \eqref{4.4.1}. However, the radii of conformality of these maps necessarily shrink as $ c\to c^{*+} $ since $ h_c(\z_c)=0 $ and $ \z_c $ approaches $ \boldsymbol \beta_1 $ in this situation. Hence, we use a special construction when $ c $ close to and larger than $ c^* $. The next lemma was shown in \cite[Lemma~7.6]{ApDenYa21}.

\begin{Lemma}
\label{lem:4.3}
The constant $ \delta_{\beta_1} $ from Lemma~{\rm\ref{lem:4.2}} can be adjusted so that there exist $ c^\prime>c^* $ and functions \smash{$ \hat\zeta_{c,\beta_1}(z) $}, $ c\in[c^*,c^\prime] $, conformal in $ \{|z-\beta_1| \leq \delta_{\beta_1}c\} $, satisfying
\begin{equation}
\label{4.3.4}
-\frac34\int_{\beta_1}^z\bigl(h_c^{(0)}-h_c^{(1)}\bigr)(s) {\rm d}s = \hat\zeta_{c,\beta_1}^{3/2}(z) - \hat\zeta_{c,\beta_1}(\beta_1+\epsilon_c)\hat\zeta_{c,\beta_1}^{1/2}(z)
\end{equation}
for some $ \epsilon_c>0 $. Each conformal map $ \hat\zeta_{c,\beta_1}(z) $ is positive on $ (\beta_1,\alpha_2) $ and negative on $ (\alpha_1,\beta_1) $ with a simple zero at $ \beta_1 $. Moreover, they form a continuous family in parameter $ c\in[c^*,c^\prime] $ and~$ \hat\zeta_{c^*,\beta_1}(z) = \zeta_{\beta_1}(z) $ $($recall that $ \beta_{c^*,1}=\beta_1)$.
\end{Lemma}

Similarly to \eqref{4.1.1}, given $ c\in( c^*,1) $, define
\begin{equation}
\label{4.4.1}
\zeta_{c,\beta_1}(z) := \left(\frac14\int_{\beta_1}^z\bigl(h_c^{(0)}-h_c^{(1)}\bigr)(s) {\rm d}s\right)^2, \qquad \alpha_1<\re z<\alpha_2.
\end{equation}
Then the following lemma holds, see \cite[Lemma~7.7]{ApDenYa21}.

\begin{Lemma}
\label{lem:4.4}
There exists a continuous and non-vanishing function $ \delta_{\beta_1}(c)>0 $ on $ (c^*,1] $ such that $ \zeta_{c,\beta_1}(z) $ is conformal in $ \{z : |z-\beta_1|\leq\delta_{\beta_1}(c) \} $, has a simple zero at $ \beta_1 $, is positive on~$ (\beta_1,\alpha_2) $ and negative on~$ (\alpha_1,\beta_1) $. The constant $ A_{\beta_1} $ in Lemma~{\rm\ref{lem:4.2}} can be adjusted so that $ 4A_{\beta_1}(z_c-\beta_1) \leq |\zeta_{c,\beta_1}^\prime(\beta_1)| $.
\end{Lemma}

For any $ c\in(0,1) $, define
\begin{equation}
\label{calH}
\mathcal H_c(z) := \re \left( \int_{\beta_{c,1}}^z \bigl( h_c^{(0)}-h_c^{(1)}\bigr)(s){\rm d}s \right), \qquad \alpha_1<\re z<\alpha_2,
\end{equation}
(please, note the change in notation as compared to \cite[Lemma~7.8]{ApDenYa21}, see \cite[equation~(7.30)]{ApDenYa21}).

\begin{Lemma}
\label{lem:4.6}
The constant $ \delta_{\beta_1} $ can be made smaller, if necessary, so that for any $ \delta\in(0,\delta_{\beta_1}] $ it holds that
\[
\mathcal H_c(x+\ic y) \leq -B\delta^{3/2}c, \qquad x\in[\beta_{c,1}+\delta c,\alpha_2-\delta c], \quad y\in[-\delta c/2,\delta c/2],
\]
for any $ c\in(0,c^*) $ and some constant $ B>0 $ independent of $ \delta $ and $ c $. Moreover, for any fixed~${ r>0} $ small enough there exist $ c_r>0 $ and $ \epsilon(r)>0 $ such that
\[
\mathcal H_c(x+\ic y) \leq -\epsilon(r), \qquad x\in[\alpha_1+r,\alpha_2-r],\quad y\in[-r/2,r/2],
\]
for all $ c\in(0,c_r) $. Furthermore, it holds that
\[
\mathcal H_c(x \pm \ic\delta c) \geq B\delta^{5/2} c, \qquad x\in[\alpha_1,\beta_{c,1}], \quad c\in(0,1).
\]
\end{Lemma}

\subsection{Additional properties}

Two questions were left unanswered in \cite[Section~7.4]{ApDenYa21}, namely, the behavior of the constants $ \epsilon_c $ in Lemma~\ref{lem:4.3} and of the function $ \delta_{\beta_1}(c) $ in Lemma~\ref{lem:4.4} as $ c\to c^{*+} $.

\begin{Lemma}
\label{lem:4.5}
The limit $ \lim_{c\to c^{*+}}\epsilon_c/\sqrt{z_c-\beta_1} $ exists and is non-zero. Moreover, $ \delta_{\beta_1}(c) \geq C\allowbreak\times\sqrt{z_c-\beta_1} $, $ c>c^* $, for some constant $ C $ independent of $ c $.
\end{Lemma}
\begin{proof}
Because $ h_c(\z) $ has a simple pole at $ \boldsymbol\beta_1 $ when $ c>c^* $, we can write
\[
h_c^{(0)}(x) - h_c^{(1)}(x) = 2u_c (x - \beta_1)^{-1/2} + \sum_{k=0}^\infty 2u_{k,c}(x - \beta_1)^{k+1/2},
\]
where $ (x-\beta_1)^{1/2} $ is a branch positive for $ x>\beta_1 $. Then it holds that
\[
\frac1{(x-\beta_1)^{1/2}} \left( -\frac34\int_{\beta_1}^x \bigl(h_c^{(0)}-h_c^{(1)}\bigr)(s) {\rm d}s \right) = -3u_c - \sum_{k=0}^\infty \frac{3u_{k,c}}{2k+3} (x-\beta_1)^{k+1}.
\]
Thus, we get from \eqref{4.3.4} and the last claim of Lemma~\ref{lem:4.3} that
\[
-3u_c = -\hat\zeta_{c,\beta_1}(\beta_1+\epsilon_c)\hat\zeta_{c,\beta_1}^\prime(\beta_1)^{1/2} = -(1+o(1)) \bigl(\zeta_{\beta_1}^\prime(\beta_1)^{3/2}\epsilon_c + \mathcal O\bigl(\epsilon_c^2\bigr) \bigr)
\]
as $ c\to c^{*+} $, where \smash{$ \zeta_{\beta_1}^\prime(\beta_1) > 0 $} by Lemma~\ref{lem:4.2}. The first claim of the lemma now follows from the fact that \smash{$ u_c/\sqrt{z_c-\beta_1} $} has a limit as $ c\to c^{*+} $, which has been shown in the proof of \cite[Lemma~7.7]{ApDenYa21}.

The function $ \delta_{\beta_1}(c) $ in Lemma~\ref{lem:4.4} was defined simply as the largest radius of conformality of~$ \zeta_{c,\beta_1}(z) $ from \eqref{4.4.1} (times a fixed constant less than $ 1 $ to ensure conformality in the closed disk). Combining \eqref{4.3.4} and \eqref{4.4.1}, we can see that
\[
9\zeta_{c,\beta_1}(z) = \hat\zeta_{c,\beta_1}(z)\bigl(\hat\zeta_{c,\beta_1}(z) - \hat\zeta_{c,\beta_1}(\beta_1+\epsilon_c) \bigr)^2 =: (z-\beta_1)(z-\beta_1-\epsilon_c)^2 F_c(z),
\]
where the functions $ F_c(z) $ are analytic and non-vanishing in~$ \{|z-\beta_1| \leq \delta_{\beta_1}c^* \} $, continuously depend on the parameter $ c $ (this is a property of the functions $ h_c(\z) $), and converge there uniformly to $ F_{c^*}(z) := (\zeta_{\beta_1}(z)/(z-\beta_1) )^3 $ as $ c\to c^{*+} $. Observe that $ F_{c^*}(\beta_1)>0 $ by Lemma~\ref{lem:4.2} and therefore the values $ F_c(\beta_1) $ are uniformly separated away from zero. Similarly, there exists a constant~$ K $, independent of $ c $, such that
$
|F_c(z)|,~|F_c^\prime(z)|,~|((z-\beta_1)F_c(z))^\prime| \leq K$
for $ z $ in~$ {\{|z-\beta_1| \leq \delta_{\beta_1}c^* \}} $ and any $ c\in[c^*,c^\prime] $. To simplify the notation slightly, let $ G_c(z) = F_c(z+\beta_1) $. Let $ A $ be such that $ 0< A\leq \min_{c^*\leq c\leq c^\prime}G_c(0)/(20K)$. For further simplicity, assume that $ c $ is sufficiently close to $ c^* $ so that $ \epsilon_c\leq 1 $ and $ A\epsilon_c<\delta_{\beta_1}c^* $. Trivially, it holds that
\begin{align*}
D(z_1,z_2):={} &\frac{z_1(z_1-\epsilon_c)^2G_c(z_1)-z_2(z_2-\epsilon_c)^2G_c(z_2)}{z_1-z_2}\\
 ={} & \epsilon_c^2 G_c(0) + z_1 \frac{(z_1-\epsilon_c)^2 G_c(z_1)-\epsilon_c^2 G_c(0)}{z_1}\\
 & + z_2 \frac{(z_1-\epsilon_c)^2 G_c(z_1)-(z_2 - \epsilon_c)^2 G_c(z_2)}{z_1-z_2}.
\end{align*}
Let $ z_1$, $z_2 $ be in $ \{ |z| \leq A\epsilon_c\} $. Since $ \epsilon_c$, $A\leq 1 $, we get that
\[
\left|\frac{(z_1-\epsilon_c)^2 G_c(z_1)-\epsilon_c^2 G_c(0)}{z_1}\right| \leq |z_1G_c(z_1)| + 2\epsilon_c|G_c(z_1)| + \epsilon_c^2\left|\frac{G_c(z_1)-G_c(0)}{z_1}\right| \leq 4K\epsilon_c.
\]
Similarly, we obtain that
\begin{gather*}
\left| \frac{(z_1-\epsilon_c)^2 G_c(z_1)-(z_2 - \epsilon_c)^2 G_c(z_2)}{z_1-z_2} \right| \\
\qquad\leq |z_1+z_2||G_c(z_1)| + |z_2|^2\left|\frac{G_c(z_1)-G_c(z_2)}{z_1-z_2}\right| \\
\qquad\phantom{\leq}{} + 2\epsilon_c \left| \frac{z_1G_c(z_1)-z_2G_c(z_2)}{z_1-z_2} \right| + \epsilon_c^2\left|\frac{G_c(z_1)-G_c(z_2)}{z_1-z_2}\right| \leq 6K\epsilon_c.
\end{gather*}
Hence, we can conclude that
$
|D(z_1,z_2)| \geq \epsilon_c^2 G_c(0) - 10AK\epsilon_c^2 \geq \epsilon_c^2 G_c(0)/2>0$.
That is, we have shown that $ z(z-\epsilon_c)^2G_c(z) $ is conformal in~$ {\{ |z|<A\epsilon_c\}} $. Thus, $ \zeta_{c,\beta_1}(z) $ is conformal in~$ {\{|z-\beta_1|<A\epsilon_c\}} $ and therefore $ \delta_{\beta_1}(c)\geq A\epsilon_c $. In view of the first claim of the lemma, the second one follows.
\end{proof}

Notice that we can assume that $ \delta_{\beta_1}(c) $ is an increasing function of $ c\in(c^*,1] $. Indeed, this always can be achieved by replacing $ \delta_{\beta_1}(c) $ with $ \min_{x\in [c,1]} \delta_{\beta_1}(x) $. The corresponding bound of Lemma~\ref{lem:4.5} will not change since
\[
 \min_{x\in [c,1]} \delta_{\beta_1}(x) = \delta_{\beta_1}(x_c) \geq C \sqrt{z_{x_c}-\beta_1}\geq C\sqrt{z_c-\beta_1}
\]
for some $ x_c \in [c,1] $, where the last inequality holds since $ z_c $ is an increasing function of $ c $.

\subsection{Polygons of conformality}
\label{ss:pc}

In this subsection, we describe domains of conformality that we shall use to define the ``lens''. The construction is not deep but somewhat technical.

Let constants $ \delta_{\alpha_1}$, $\delta_{\beta_1} $ be as in the preceding lemmas of this section and constants $ \delta_{\alpha_2}$, $\delta_{\beta_2} $ be defined similarly but with respect to $ \Delta_2 $. It follows from \eqref{c-rate} that
\begin{equation}
\label{dD}
\delta_\Delta := \min\left\{ 1,\delta_{\alpha_1},\delta_{\beta_1}, \delta_{\alpha_2},\delta_{\beta_2}, \frac{\alpha_2-\beta_1}3, \inf_{c\in(0,1)}\left\{\frac{|\Delta_{c,1}|}{3c},\frac{|\Delta_{c,2}|}{3(1-c)}\right\}\right\} > 0.
\end{equation}

Because the conformal maps $ \hat\zeta_{c,\beta_1}(z) $ from Lemma~\ref{lem:4.3} form a continuous family in the parameter $ c\in[c^*,c^\prime] $, it follows from Lemma~\ref{lem:4.5} that there exists a~constant $ K_1>0 $ such that~$
\hat\zeta_{c,\beta_1}(\beta_1+\epsilon_c) \leq K_1\sqrt{z_c-\beta_1}$.
On the other hand, the same continuity in the parameter $ c $ and \eqref{koebe} imply that there exists a~constant $ K_2>0 $ such that for every $ \delta\in(0,\delta_\Delta) $ and~$ c\in[c^*,c^\prime] $ it holds that
\begin{equation}
\label{zhat-b}
\min\bigl\{ \big|\hat\zeta_{c,\beta_1}(s)\big| : |s-\beta_1|=\delta c^*/\sqrt2\bigr\} \geq K_2 \delta.
\end{equation}
Since $ z_c $ is an increasing function of $ c $, given $ \delta\in(0,\delta_\Delta) $, there exists a unique $ c^\prime(\delta)>c^* $ such that $ z_{c^\prime(\delta)} - \beta_1 = (K_2/K_1)^2\delta^2 $, where we adjust the constants $ K_1$, $K_2 $ so that $ c^\prime(\delta_\Delta)\leq c^\prime $. Then
\begin{equation}
\label{min-con}
\hat\zeta_{c,\beta_1}(\beta_1+\epsilon_c) \leq \min\bigl\{ \big|\hat\zeta_{c,\beta_1}(s)\big| : |s-\beta_1|=\delta c/\sqrt2\bigr\}, \qquad c\in[c^*,c^\prime(\delta)],
\end{equation}
a technical inequality that will be important to us later. Another consequence of this definition of $ c^\prime(\delta) $ and Lemma~\ref{lem:4.5} is that
\smash{$
K^\prime := \min\bigl\{ 1, \inf_{0<\delta<\delta_\Delta} \frac{\delta_{\beta_1}(c^\prime(\delta))}{\delta}\bigr\} >0$}.

Denote by $ U(z,r) $ the interior of the square with vertices $ z\pm r$, $z\pm\ic r $. For any $ \delta\in(0,\delta_\Delta) $, set
\begin{equation}
\label{Kdeltac}
K = K(\delta,c) := \begin{cases}
1/3, & c<c^\prime(\delta), \beta_1-\beta_{c,1}>2\delta c/3, \\
1, & c<c^\prime(\delta), \beta_1-\beta_{c,1} \leq 2\delta c/3, \\
K^\prime, & c\geq c^\prime(\delta).
\end{cases}
\end{equation}
Since $ \delta_{\beta_1}(c) $ is an increasing function of $ \delta $, it therefore holds that $ \delta_{\beta_1}(c) \geq \delta_{\beta_1}(c^\prime(\delta)) \geq K\delta $ when~$ c \in [c^\prime(\delta),1) $. Define
\begin{equation}
\label{Uc1}
U_{c,e} := U(e,K\delta c), \qquad e\in \{\alpha_1,\beta_{c,1} \}.
\end{equation}
This definition achieves the following:
\begin{itemize}\itemsep=0pt
\item the map $ \zeta_{c,\alpha_1}(z) $ from Lemma~\ref{lem:4.1} is conformal in $ U_{c,\alpha_1} $;
\item since $ \beta_1-\beta_{c,1} $ is a decreasing function of $ c $ while $ 2\delta c/3 $ is clearly increasing, the squares $ U_{c,\beta_{c,1}} $ start out (as $ c $ increases from $0$) with $ K = 1/3 $ and in these cases the point $ \beta_1 $ does not belong to the interior of the squares and lies distance at least $ \delta c/3 $ from their boundary;
\item when the parameter $ c $ reaches the value for which $ \beta_1-\beta_{c,1}=2\delta c/3 $, the value of $ K $ changes to $ 1 $ and from that point on $ \beta_1 $ belongs to the interior of the squares $ U_{c,\beta_{c,1}} $ and lies distance at least $ \delta c/3\sqrt2 $ from their boundary;
\item when $ \beta_1 $ belongs to the interior of the square $ U_{c,\beta_{c,1}} $ it is not at its center unless $ c\geq c^* $;
\item the map $ \zeta_{\beta_{c,1}}(z) $ from Lemma~\ref{lem:4.2} is conformal within $ U_{c,\beta_{c,1}} $ for each $ c\leq c^* $;
\item when $ c\in(c^*,c^\prime(\delta)) $ the point $ \beta_1+\epsilon_c $ belongs to the interior of the square $ U_{c,\beta_1} $ by \eqref{min-con} and the map $ \hat\zeta_{c,\beta_1}(z) $ from Lemma~\ref{lem:4.3} is conformal in this square;
\item when $ c\in[c^\prime(\delta),1) $, the map $ \zeta_{c,\beta_1}(z) $ from Lemma~\ref{lem:4.4} is conformal within the square $ U_{c,\beta_1} $ whose size is proportional to $ \delta c $ as in all the other cases, which was the motivating reason behind the definition of $ c^\prime(\delta) $.
\end{itemize}

\begin{figure}[!ht]\centering

\subfigure[]{\includegraphics[scale=1.5]{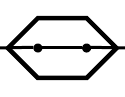}
\begin{picture}(0,0)
\put(-74,41){$\beta_{c,1}+\delta c/2$}
\put(-35,23){$\beta_1$}
\end{picture}} \qquad\qquad
\subfigure[]{\includegraphics[scale=1.5]{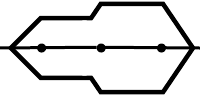}
\begin{picture}(0,0)
\put(-125,23){$\beta_{c,1}+\delta c/2$}
\put(-90,41){$\alpha_1+3r/2$}
\put(-35,23){$\beta_1$}
\end{picture}
}
\caption{Polygons $ \partial U_{c,\beta_1} $; panel (a) $ c\in[c_r,c^*) $; panel (b) $ c\in(0,c_r) $.}
\label{fig:polygon}
\end{figure}
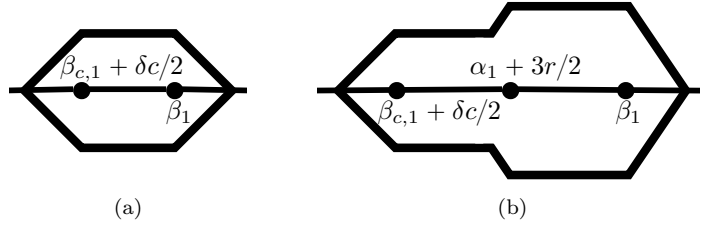

When $ \beta_1 $ does not belong to the interior of the square $ U_{c,\beta_{c,1}} $, we need to define a polygon~$ U_{c,\beta_1} $ containing $ \beta_1 $ in its interior. For reasons that have to do with future asymptotic analysis, the definition of this polygon is rather technical as well. Let $ r>0 $ be small enough (in particular, $ 3r/2<\beta_1-\alpha_1 $) and $ c_r $ be as in Lemma~\ref{lem:4.6}. When $ c\in(0,c^*) $ and $ \beta_1\not\in U_{c,\beta_{c,1}} $, we set
\begin{equation}
\label{Uc2}
U_{c,\beta_1} := \bigcup_{x\in[\beta_{c,1}+\delta c/2,\beta_1]}U(x,\delta c/6) \cup \bigcup_{x\in[\alpha_1+3r/2,\beta_1]} U(x,r/2),
\end{equation}
where the second union is present only if $ c<c_r $ and always stays within the rectangle of the second estimate of Lemma~\ref{lem:4.6}, see Figure~\ref{fig:polygon} (we always can decrease $ c_r $ if necessary so that $ \beta_{c_r,1}+c_r/2 < \alpha_1 + 3r/2 $). Notice that the domains $ U_{c,\alpha_1} $, $ U_{c,\beta_{c,1}} $, and $ U_{c,\beta_1} $ are disjoint, however, $ \partial U_{c,\beta_{c,1}} $ and $ \partial U_{c,\beta_1} $ share a common point when both sets are distinct and non-empty.

The domains $ U_{c,\alpha_2} $, $ U_{c,\alpha_{c,2}} $, and $ U_{c,\beta_2} $ can be defined similarly. We use domains with polygonal boundary rather then disks for a not very deep reason that in this case it is easier to explain uniform boundedness of Cauchy operators on our variable lenses, see \cite[Lemma~7.9]{ApDenYa21}.

\section{Orthogonal polynomials and Riemann--Hilbert problems}

To slightly simplify the notation we agree that from now on all the quantities that depend on the parameter $ c $ will simply be labeled by $ \n $ when referred to with $ c=c(\n)=n_1/|\n| $.

 We let $ [\boldsymbol A]_{i,j} $ stand for $ (i,j) $-th entry of a matrix $ \boldsymbol A $ and $ \boldsymbol E_{i,j} $ be the matrix with all zero entries except for $ [\boldsymbol E_{i,j}]_{i,j} =1 $. Also, we set $\sigma(\n) := \diag\left(|\n|, -n_1,-n_2\right)$, $ \n=(n_1,n_2) $.

\subsection{Initial RH problem}

Consider the following Riemann--Hilbert problem (\rhy): find a $ 3\times3 $ matrix function $ \boldsymbol Y(z) $ such that
\begin{itemize}\itemsep=0pt
\label{rhy}
\item[(a)] $\boldsymbol Y(z)$ is analytic in $\C\setminus(\Delta_1\cup \Delta_2)$ and $\displaystyle \lim_{z\to\infty} {\boldsymbol Y}(z)z^{-\sigma(\n)} = {\boldsymbol I}$;
\item[(b)] $\boldsymbol Y(z)$ has continuous boundary values on each $\Delta_i^\circ$ that satisfy
\[
{\boldsymbol Y}_+(x) = {\boldsymbol Y}_-(x)( \boldsymbol I + \rho_i(x) \boldsymbol E_{1,i+1}), \qquad i\in\{1,2\};
\]
\item[(c)] the entries of the $(i+1)$-st column of $\boldsymbol Y(z)$ behave like $\mathcal{O}\left(\log|z-e|\right)$ while the remaining entries stay bounded as $z\to e\in\{\alpha_i,\beta_i\}$, $ i\in\{1,2\} $.
\end{itemize}

Let $ R_\n^{(i)}(z) $, $ i\in\{1,2\} $, be the $ i $-th function of the second kind associated with $ P_\n(x) $. That~is,%
\begin{equation}
\label{Rni}
R_\n^{(i)}(z) := \frac1{2\pi\ic}\int\frac{P_\n(x)\rho_i(x){\rm d}x}{x-z} = \frac{h_{\n,i}}{z^{n_i+1}} + \mathcal O\left(\frac1{z^{n_i+2}}\right), \qquad z\in\overline\C\setminus\Delta_i,
\end{equation}
where the estimate follows from \eqref{typeII} and \eqref{connection} and holds as $ z\to\infty $. Assume that the multi-index $ \n $ is such that
\begin{equation}
\label{normality_con}
\deg(P_\n) = |\n| \qquad \text{and}\qquad  h_{\n-\vec e_1,1}h_{\n-\vec e_2,2} \neq 0.
\end{equation}
The following lemma holds, see \cite[Proposition~3.1]{Ya16}.

\begin{Lemma}
\label{lem:rhy}
If $ P_\n(z) $ satisfying \eqref{typeII} and $ R_\n^{(i)}(z) $, $ i\in\{1,2\} $, given by \eqref{Rni}, satisfy \eqref{normality_con}, then~\hyperref[rhy]{\rhy} is solved by
\begin{equation}
\label{Y}
\boldsymbol Y(z) = \begin{pmatrix}
P_\n(z) & R_\n^{(1)}(z) & R_\n^{(2)}(z) \\[1mm]
h_{\n-\vec e_1,1}^{-1}P_{\n-\vec e_1}(z) & h_{\n-\vec e_1,1}^{-1}R_{\n-\vec e_1}^{(1)}(z) & h_{\n-\vec e_1,1}^{-1}R_{\n-\vec e_1}^{(2)}(z) \\[1mm]
h_{\n-\vec e_2,2}^{-1}P_{\n-\vec e_2}(z) & h_{\n-\vec e_2,2}^{-1}R_{\n-\vec e_2}^{(1)}(z) & h_{\n-\vec e_2,2}^{-1}R_{\n-\vec e_2}^{(2)}(z)
\end{pmatrix}.
\end{equation}
Conversely, if a solution of \hyperref[rhy]{\rhy} exists, then it is unique and is given by \eqref{Y} with $ P_\n(z) $ and \smash{$ R_\n^{(i)}(z) $}, $ i\in\{1,2\} $, necessarily satisfying \eqref{normality_con}.
\end{Lemma}

\subsection{Opening of the lenses}
\label{ss:OL}

The next step in the Riemann--Hilbert analysis consists in factorizing the jump matrix and moving some of the jump relations into the complex plane, the so-called ``opening of the lenses''. In constructing this lens we rely heavily on the material of Section~\ref{ss:pc}.

\begin{figure}[!ht]\centering
\includegraphics[scale=1.1]{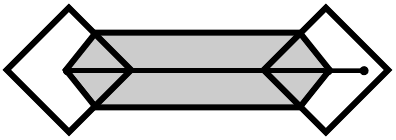}
\begin{picture}(0,0)
\put(-190,32){$\alpha_1$}
\put(-38,27){$\beta_{\n,1}$}
\put(-228,55){$\partial U_{\n,\alpha_1}$}
\put(-22,55){$\partial U_{\n,\beta_{\n,1}}$}
\put(-115,65){$\Gamma^+_{\n,1}$}
\put(-115,5){$\Gamma^-_{\n,1}$}
\put(-115,44){$\Omega_{\n,1}$}
\end{picture}
\caption{The boundaries $ \partial U_{\n,\alpha_1} $ and $ \partial U_{\n,\beta_{\n,1}} $, arcs $\Gamma_{\n,1}^\pm$, and domains $ \Omega_{\n,1}^\pm $ (shaded).}
\label{fig:lens}
\end{figure}
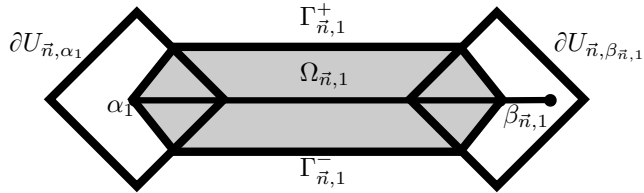

Given $ \delta\in(0,\delta_\Delta] $, see \eqref{dD}, let $ U_{\n,e} $, $ e\in \{\alpha_1,\beta_{\n,1},\beta_1\} $, be the squares defined in \eqref{Uc1} and~\eqref{Uc2} (and via similar formulae at $ \{\alpha_2,\alpha_{\n,2},\beta_2\} $). Recall the definition of the rays $ I_\pm $ after~\eqref{K} and the conformal maps constructed in \eqref{4.1.1}--\eqref{4.4.1}. Let $ \tau_* $ be as in Section~\ref{ss:sse}, $ c^* $ as in \eqref{cstars}, and $ c^\prime(\delta) $ as defined before \eqref{min-con}. Set $ \zeta_{\n,\alpha_1}(z) := \zeta_{c(\n),\alpha_1}(z) $, see \eqref{4.1.1},
\begin{equation}
\zeta_{\n,\beta_1}(z) :=
\begin{cases}
\zeta_{\beta_{\n,1}}(z), & c(\n) \leq c^*, \zeta_{\beta_{\n,1}}(\beta_1)>\tau_*|\n|^{-2/3}, \\
\zeta_{\beta_{\n,1}}(z)-\zeta_{\beta_{\n,1}}(\beta_1), & c(\n)\leq c^*, \zeta_{\beta_{\n,1}}(\beta_1)\leq \tau_*|\n|^{-2/3},
\end{cases}
\label{b-maps1}
\end{equation}
where $ \zeta_{\beta_{\n,1}}(z) $ is given by \eqref{4.2.1} (notice that $ \zeta_{\beta_{c,1}}(\beta_1)\geq0 $ when $ c\leq c^*$), and
\begin{equation}
\zeta_{\n,\beta_1}(z) :=
\begin{cases}
\hat\zeta_{c(\n),\beta_1}(z), & c^*<c(\n)<c^\prime(\delta), \quad \text{see} \ \eqref{4.3.4}, \\
\zeta_{c(\n),\beta_1}(z), & c^\prime(\delta) \leq c(\n), \quad \text{see} \ \eqref{4.4.1}
\end{cases}
\label{b-maps2}
\end{equation}
(let us emphasize that the choice whether \eqref{4.3.4} or \eqref{4.4.1} is used does depend on the value of $ \delta $; in the first line of \eqref{b-maps2} we also slightly departed from our labeling convention, hopefully without much ambiguity). The maps $ \zeta_{\n,\alpha_2}(z) $ and $ \zeta_{\n,\beta_2}(z) $ can be defined similarly.

We now select open Jordan arcs $ \Gamma_{\n,i}^\pm $ connecting $ \alpha_{\n,i} $ to $ \beta_{\n,i} $ so that
\begin{equation}
\label{Ipm}
\zeta_{\n,\beta_i}\bigl(\Gamma_{\n,i}^\pm\cap U_{\n,\beta_{\n,i}}\bigr) \subset I_\pm, \qquad \zeta_{\alpha_{\n,i}}\bigl(\Gamma_{\n,i}^\pm\cap U_{\n,\alpha_{\n,i}}\bigr) \subset I_\mp,
\end{equation}
and that consist of straight line segments outside of $ U_{\n,e} $, $ e\in E_\n $, see Figure~\ref{fig:lens}. These arcs are oriented from $ \alpha_{\n,i} $ to $ \beta_{\n,i} $. We denote by $ \Omega_{\n,i}^\pm $ the domains delimited by $ \Gamma_{\n,i}^\pm $ and $ \Delta_{\n,i} $.

As $ \rho_i(x) $ is a restriction of an analytic function, that we keep denoting by $ \rho_i(z) $, we can decrease the constant $ \delta_\Delta $, if necessary, so that $ \rho_i(z) $ is analytic and non-vanishing in a simply connected neighborhood of the connected component intersecting $ \Delta_i $ of the above constructed lens for any value of $ \delta\leq\delta_\Delta $.

\subsection{Factorized Riemann--Hilbert problem}

For compactness of notation, we introduce transformations $\mathrm{T}_i$, $i\in\{1,2\}$, that act on $2\times2$ matrices in the following way
\[
\mathrm{T}_1 \begin{pmatrix} e_{11} & e_{12} \\ e_{21} & e_{22} \end{pmatrix} = \begin{pmatrix} e_{11} & e_{12} & 0 \\ e_{21} & e_{22} & 0 \\ 0 & 0 & 1\end{pmatrix} \qquad \text{and}\qquad  \mathrm{T}_2 \begin{pmatrix} e_{11} & e_{12} \\ e_{21} & e_{22} \end{pmatrix} = \begin{pmatrix} e_{11} & 0 & e_{12} \\ 0 & 1 & 0 \\ e_{21} & 0 & e_{22}\end{pmatrix}.
\]
Given $\boldsymbol Y(z)$, the solution of \hyperref[rhy]{\rhy}, set
\begin{equation}
\label{X}
\boldsymbol X(z):= \boldsymbol Y(z)
\begin{cases}
\mathrm{T}_i \begin{pmatrix} 1 & 0 \\ \mp1/\rho_i(z) & 1 \end{pmatrix}, & z\in \Omega_{\n,i}^{\pm}, \\
\boldsymbol I, & \text{otherwise}.
\end{cases}
\end{equation}
Then $\boldsymbol X(z)$ solves the following Riemann--Hilbert problem (\rhx):
\begin{itemize}\itemsep=0pt
\item[(a)] $\boldsymbol X(z)$ is analytic in $\C\setminus\bigcup_{i=1}^2\bigl(\Delta_i\cup\Gamma_{\n,i}^+\cup\Gamma_{\n,i}^-\bigr)$ and $ \lim_{z\to\infty} {\boldsymbol X}(z)z^{-\sigma(\n)} = {\boldsymbol I}$;
\item[(b)] $\boldsymbol X(z)$ has continuous boundary values on $\bigcup_{i=1}^2\bigl(\Delta_i^\circ\cup\Gamma_{\n,i}^+\cup\Gamma_{\n,i}^-\bigr)$ that satisfy
\[
{\boldsymbol X}_+(s)={\boldsymbol X}_-(s)
\begin{cases}
\displaystyle \mathrm{T}_i \begin{pmatrix} 0 & \rho_i(s) \\ -1/\rho_i(s) & 0 \end{pmatrix}, & s\in \Delta_{\n,i}, \vspace{0.5mm}\\
\displaystyle \mathrm{T}_i \begin{pmatrix} 1 & 0 \\ 1/\rho_i(s) & 1 \end{pmatrix}, & s\in \Gamma_{\n,i}^+\cup\Gamma_{\n,i}^-, \vspace{0.5mm}\\
\displaystyle \mathrm{T}_i \begin{pmatrix} 1 & \rho_i(s) \\ 0 & 1 \end{pmatrix}, & s\in \Delta_i^\circ\setminus\Delta_{\n,i};
\end{cases}
\]
\item[(c)] the entries of the first and $(i+1)$-st columns of $\boldsymbol X(z)$ behave like $\mathcal{O}\left(\log|z-e|\right)$ while the remaining entries stay bounded as $z\to e\in\{\alpha_i,\beta_i\}$.
\end{itemize}

The following lemma is a combination of \cite[Lemma~8.1]{Ya16} and \cite[Lemma~6.4]{BYa10}.

\begin{Lemma}
\label{rhx}
\hyperref[rhx]{\rhx} and \hyperref[rhy]{\rhy} are simultaneously solvable and the solutions are unique and connected by \eqref{X}.
\end{Lemma}

\section{Global parametrix}

As will become apparent later, away from the intervals $ \Delta_{\n,1} $ and $ \Delta_{\n,2} $ we expect $ \boldsymbol X(z) $ to behave like the solution of the following Riemann--Hilbert problem (\rhn)
\begin{itemize}\itemsep=0pt
\label{rhn}
\item[(a)] $\boldsymbol N(z)$ is analytic in $\C\setminus(\Delta_{\n,1}\cup \Delta_{\n,2})$ and $ \lim_{z\to\infty} {\boldsymbol N}(z)z^{-\sigma(\n)} = {\boldsymbol I}$;
\item[(b)] $\boldsymbol N(z)$ has continuous boundary values on each $\Delta_{\n,i}^\circ$ that satisfy
\[
\boldsymbol N_+(s) = \boldsymbol N_-(s) \mathrm{T}_i \begin{pmatrix} 0 & \rho_i(s) \\ -1/\rho_i(s) & 0 \end{pmatrix};
\]
\item[(c)] it holds that $ \boldsymbol N(z) = \boldsymbol{\mathcal O}\bigl(|z-e|^{-1/4}\bigr) $ as $ z\to e \in E_\n $.
\end{itemize}

Let $\Phi_\n(\z)$, $ w_{\n,i}(z) $, and $S_\n(\z)$ be the functions given by \eqref{Phin}, \eqref{wci}, and \eqref{Szego}, respectively. Recall also the definition of $ \chi_\n(\z) $ in \eqref{chi} as well as \eqref{AngPar1}. Set
\begin{equation}
\label{Ups_n}
\Upsilon_{\n,i}(\z) := A_{\n,i} \bigl(\chi_\n(\z) - B_{\n,i} \bigr)^{-1}, \qquad i\in\{1,2\}.
\end{equation}
It follows from \eqref{AngPar1} that $ \Upsilon_{\n,i}(\z) $ is a conformal map of $ \RS_\n $ onto $ \overline\C $ that maps~$ \infty^{(i)} $ into $ \infty $ and~$ \infty^{(0)} $ into $ 0 $. More precisely, it holds that
\[
\Upsilon_{\n,i}^{(i)}(z) = z+ \mathcal O(1) \qquad \text{and}\qquad  \Upsilon_{\n,i}^{(0)}(z) = A_{\n,i}z^{-1} + \mathcal O\bigl(z^{-2}\bigr)
\]
as $ z\to\infty $. Put $\boldsymbol S(z):=\diag\bigl(S_\n^{(0)}(z),S_\n^{(1)}(z),S_\n^{(2)}(z)\bigr)$ and define
\begin{equation}
\label{matrix-M}
\boldsymbol M(z):=\boldsymbol S^{-1}(\infty) \begin{pmatrix}
1 & 1/w_{\n,1}(z) & 1/w_{\n,2}(z) \\[1mm]
\Upsilon_{\n,1}^{(0)}(z) & \Upsilon_{\n,1}^{(1)}(z)/w_{\n,1}(z) & \Upsilon_{\n,1}^{(2)}(z)/w_{\n,2}(z) \\[1mm]
\Upsilon_{\n,2}^{(0)}(z) & \Upsilon_{\n,2}^{(1)}(z)/w_{\n,1}(z) & \Upsilon_{\n,2}^{(2)}(z)/w_{\n,2}(z)
\end{pmatrix} \boldsymbol S(z).
\end{equation}
Then it can be readily verified using \eqref{szego-pts2} that \hyperref[rhn]{\rhn} is solved by $ \boldsymbol N(z) := \boldsymbol C (\boldsymbol{MD})(z)$, see~\cite[Section~7.3]{ApDenYa21}, where $\boldsymbol C$ is a diagonal matrix of constants such that
\begin{equation}
\label{matrix-CD}
\lim_{z\to\infty}\boldsymbol C\boldsymbol D(z)z^{-\sigma(\n)} = \boldsymbol I \qquad \text{and}\qquad  \boldsymbol D(z):=\diag\bigl(\Phi_\n^{(0)}(z),\Phi_\n^{(1)}(z),\Phi_\n^{(2)}(z)\bigr).
\end{equation}
Since the jump matrices in \hyperref[rhn]{\rhn}(b) have determinant $ 1 $, it follows from the second identity in \eqref{szego-pts2} and the normalization at infinity that $ \det \boldsymbol N(z) $ is holomorphic in the entire extended complex plane except for at most square root singularities at the elements of $ E_\n $. As those singularities are isolated, they are removable and $ \det \boldsymbol N(z)\equiv1 $. In fact, it holds that $ \det \boldsymbol M(z)\equiv\det \boldsymbol D(z)\equiv\det \boldsymbol C=1 $.

\begin{Lemma}
\label{lem:Mbounds1}
It holds that\footnote{We write $ |A(z)| \sim |B(z)| $ if $ C^{-1}|A(z)| \leq |B(z)| \leq C|A(z)| $ for some $ C>1 $ and $ |\boldsymbol A(z)| \sim |\boldsymbol B(z)| $ if $ |[\boldsymbol A]_{i,j}(z)| \sim |[\boldsymbol B]_{i,j}(z)| $ for each pair $ i,j\in\{1,2,3\} $.} $ \boldsymbol M(z) = \boldsymbol{\mathcal O}\bigl(\delta^{-1/2}\bigr) $ uniformly for $ z $ such that $ \delta c(\n) \leq \dist(z,E_{\n,1} ) $ and $ \delta(1-c(\n)) \leq \dist(z,E_{\n,2}) $, where the constants in $ \boldsymbol{\mathcal O}(\cdot) $ are independent of $ c(\n) $ and $ \delta $. Moreover, it holds that
\[
|\boldsymbol M(z)| \sim \begin{pmatrix} \delta^{-1/4} & \delta^{-1/4} & 1-c(\n) \\ \delta^{-1/4} & \delta^{-1/4} & c(\n)(1-c(\n)) \\ (1-c(\n))\delta^{-1/4} & (1-c(\n))\delta^{-1/4} & 1 \end{pmatrix}
\]
uniformly on $ |z-\alpha_1| = \delta c(\n) $ and $ |z-\beta_{\n,1}| = \delta c(\n) $ and a similar formula holds on $ |z-\alpha_{\n,2}| = \delta (1-c(\n)) $ and $ |z-\beta_2| = \delta (1-c(\n)) $, where the constants of proportionality are independent of~$ c(\n) $ and $ \delta $.
\end{Lemma}
\begin{proof}
This is \cite[Lemma~7.3]{ApDenYa21}. There it was stated that $ \boldsymbol M(z) = \boldsymbol{\mathcal O}_\delta(1) $. However, the actual proof shows that $ \boldsymbol M(z) = \boldsymbol{\mathcal O}\bigl(\delta^{-1/2}\bigr) $, see also the forthcoming Lemma~\ref{lem:Mbounds2}.
\end{proof}

Let $ \Pi_\n(\z) $ be the rational function defined in \eqref{Pin}. Furthermore, let $ \Pi_{\n,i}(\z) $, $ i\in\{1,2\} $, be the rational functions on $ \RS_\n $ with the zero/pole divisors and the normalizations given by
\begin{equation}
\label{Pini}
\infty^{(0)} + \infty^{(i)} +2 \infty^{(3-i)} - \boldsymbol\alpha_1 - \boldsymbol\beta_{\n,1} - \boldsymbol\alpha_{\n,2} - \boldsymbol\beta_2 \qquad \text{and}\qquad  \Pi_{\n,i}^{(i)}(z) = \frac1z + \mathcal O\left(\frac1{z^2}\right)
\end{equation}
as $ z\to\infty $. It was shown in \cite[equation~(7.7)]{ApDenYa21} that
\begin{equation}
\label{M-inverse}
\boldsymbol M^{-1}(z) = \boldsymbol S^{-1}(z) \begin{pmatrix} \Pi_\n^{(0)}(z) & \Pi_{\n,1}^{(0)}(z) & \Pi_{\n,2}^{(0)}(z) \\[1mm] w_{\n,1}(z)\Pi_\n^{(1)}(z) & w_{\n,1}(z)\Pi_{\n,1}^{(1)}(z) & w_{\n,1}(z)\Pi_{\n,2}^{(1)}(z) \\[1mm] w_{\n,2}(z)\Pi_\n^{(2)}(z) & w_{\n,2}(z)\Pi_{\n,1}^{(2)}(z) & w_{\n,2}(z)\Pi_{\n,2}^{(2)}(z) \end{pmatrix} \boldsymbol S(\infty).
\end{equation}
Then the following lemma takes place.

\begin{Lemma}
\label{lem:Mbounds2}
It holds that $ \boldsymbol M^{-1}(z) = \boldsymbol{\mathcal O}\bigl(\delta^{-1/2}\bigr) $ uniformly for $ z $ satisfying $ \delta c(\n) \leq \dist(z,E_{\n,1} ) $ and $ \delta(1-c(\n)) \leq \dist(z,E_{\n,2}) $, where the constants in $ \boldsymbol{\mathcal O}(\cdot) $ are independent of $ c(\n) $ and $ \delta $. Moreover, it holds that
\[
\boldsymbol M^{-1}(z) = \boldsymbol{\mathcal O} \begin{pmatrix} \delta^{-1/4} & \delta^{-1/4} & (1-c(\n))\delta^{-1/4} \\ \delta^{-1/4} & \delta^{-1/4} & (1-c(\n))\delta^{-1/4} \\ 1-c(\n) & 1-c(\n) & 1 \end{pmatrix}
\]
uniformly on $ |z-\alpha_1| = \delta c(\n) $ and $ |z-\beta_{\n,1}| = \delta c(\n) $ and a similar formula holds on $ |z-\alpha_{\n,2}| = \delta (1-c(\n)) $ and $ |z-\beta_2| = \delta (1-c(\n)) $, where $\boldsymbol{\mathcal O}(\cdot) $ is independent of $ c(\n) $ and $ \delta $.
\end{Lemma}
\begin{proof}
It was shown in \cite[Lemma~7.3]{ApDenYa21} that
\begin{gather}
\boldsymbol S(\infty) \sim \diag\bigl(c(\n)^{-1/3}(1-c(\n))^{-1/3},c(\n)^{2/3}(1-c(\n))^{-1/3},c(\n)^{-1/3}(1-c(\n))^{2/3}\bigr),\nonumber \\
|\boldsymbol S(z)| \sim \boldsymbol S(\infty)\diag\bigl(\delta^{-1/4},\delta^{1/4},1\bigr),\label{6.2.1}
\end{gather}
uniformly on $ |z-\alpha_1| = \delta c(\n) $ and $ |z-\beta_{\n,1}| = \delta c(\n) $, where the constants of proportionality are independent of $ c(\n) $ and $ \delta $. It was further shown in \cite[Lemma~5.3]{ApDenYa21} that
\begin{equation}
\label{6.2.3}
(-1)^{3-i}(w_{\n,1}w_{\n,2})(z)\Pi_{\n,3-i}(\z) =
\begin{cases}
\bigl(\Upsilon_{\n,i}^{(2)}-\Upsilon_{\n,i}^{(1)}\bigr)(z), & \z\in\RS_\n^{(0)}, \\[1mm]
\bigl(\Upsilon_{\n,i}^{(0)}-\Upsilon_{\n,i}^{(2)}\bigr)(z), & \z\in\RS_\n^{(1)}, \\[1mm]
\bigl(\Upsilon_{\n,i}^{(1)}-\Upsilon_{\n,i}^{(0)}\bigr)(z), & \z\in\RS_\n^{(2)},
\end{cases}
\end{equation}
for $ i\in\{1,2\} $ and
\begin{equation*}
(w_{\n,1}w_{\n,2})(z)\Pi_\n(\z) =
\begin{cases}
\bigl(\Upsilon_{\n,2}^{(2)}\Upsilon_{\n,1}^{(1)}-\Upsilon_{\n,2}^{(1)}\Upsilon_{\n,1}^{(2)}\bigr)(z), & \z\in\RS_\n^{(0)}, \\[1mm]
\bigl(\Upsilon_{\n,2}^{(0)}\Upsilon_{\n,1}^{(2)}-\Upsilon_{\n,2}^{(2)}\Upsilon_{\n,1}^{(0)}\bigr)(z), & \z\in\RS_\n^{(1)}, \\[1mm]
\bigl(\Upsilon_{\n,2}^{(1)}\Upsilon_{\n,1}^{(0)}-\Upsilon_{\n,2}^{(0)}\Upsilon_{\n,1}^{(1)}\bigr)(z), & \z\in\RS_\n^{(2)}.
\end{cases}
\end{equation*}
Next, it was proven, see \cite[equation~(5.23)]{ApDenYa21}, that
\[
\sqrt{3(\alpha_2-\beta_1)} < |w_{\n,1}(z)|/ (c(\n)\sqrt\delta) < 3\sqrt{\beta_2-\alpha_1}
\]
on $ |z-\alpha_1| = \delta c(\n) $ and $ |z-\beta_{\n,1}| = \delta c(\n) $. Finally, it was deduced in \cite[Lemma~5.2]{ApDenYa21} that
\begin{gather}
c(\n)^{-1}\big|\Upsilon_{\n,1}^{(0)}(z)\big|, c(\n)^{-1}\big|\Upsilon_{\n,1}^{(1)}(z)\big|, c(\n)^{-2}\big|\Upsilon_{\n,1}^{(2)}(z)\big| \sim 1, \nonumber \\
(1-c(\n))^{-2}\big|\Upsilon_{\n,2}^{(0)}(z)\big|, (1-c(\n))^{-2}\big|\Upsilon_{\n,2}^{(1)}(z)\big|, \big|\Upsilon_{\n,2}^{(2)}(z)\big| \sim 1,\label{6.2.6}
\end{gather}
on $ \{z : \dist(z,\Delta_{c,1})\leq c\delta_*\} $ for all $ 0<\delta_*\leq (\alpha_2-\beta_1)/2 $, where the constants of proportionality depend only on $ \delta_* $. Then it follows from \eqref{M-inverse} and \eqref{6.2.3}--\eqref{6.2.6} that
\begin{equation}
\label{6.2.7}
\boldsymbol M^{-1}(z) = \boldsymbol S^{-1}(z) \boldsymbol{\mathcal O} \begin{pmatrix} \delta^{-1/2} & c(\n)^{-1}\delta^{-1/2} & \delta^{-1/2} \\ c(\n) & 1 & c(\n) \\ (1-c(\n))^2 & (1-c(\n))^2 c(\n)^{-1} & 1 \end{pmatrix} \boldsymbol S(\infty)
\end{equation}
uniformly on $ |z-\alpha_1| = \delta c(\n) $ and $ |z-\beta_{\n,1}| = \delta c(\n) $, where $\boldsymbol{\mathcal O}(\cdot) $ is independent of $ c(\n) $ and~$ \delta $ and one needs to observe that the functions in the third row of the middle matrix in \eqref{M-inverse} are holomorphic at $ \alpha_1 $ and $ \beta_{\n,1} $ and therefore their estimates, stated in \eqref{6.2.7}, can be obtained via~\eqref{6.2.3}--\eqref{6.2.6} and the maximum modulus principle for holomorphic functions. The second claim of the lemma now follows from \eqref{6.2.1} and \eqref{6.2.7}. To see the validity of the first claim, one needs to combine \eqref{6.2.6}--\eqref{6.2.7} with the maximum modulus principle as well as the estimate
\begin{gather*}
|\boldsymbol S_\pm(x)|^{-1} \lesssim \boldsymbol S^{-1}(\infty) \diag\bigl( 1,\delta^{-1/4},\delta^{-1/4} \bigr), \\
 x\in (\alpha_1+\delta c(\n),\beta_{\n,1}-\delta c(\n)) \cup (\alpha_{\n,2}+\delta c(\n),\beta_2-\delta c(\n)),
\end{gather*}
that follows from \cite[Lemma~5.4]{ApDenYa21}.
\end{proof}

\section{Local parametrices}

To describe the behavior of $ \boldsymbol X(z) $ within the domains $ U_{\n,e} $, we are seeking solutions of the following Riemann--Hilbert problems (\rhp$_e$):
\begin{itemize}\itemsep=0pt\setlength{\leftskip}{0.1cm}\label{rhp}
\item[(a--c)] $\boldsymbol P_e(z)$ satisfies \hyperref[rhx]{\rhx}(a--c) within $U_{\n,e}$;
\item[(d)] $\boldsymbol P_e(s)=\boldsymbol M(s)\bigl(\boldsymbol I+\boldsymbol{\mathcal O} \bigl( \delta^{-3/2} \varepsilon_\n^{1/3} \bigr)\bigr)\boldsymbol D(s)$ uniformly on $\partial U_{\n,e} $.
\end{itemize}

The asymptotic formula in \hyperref[rhp]{\rhp$_e$}(d) will hold as long as $ \delta^{-3/2}\varepsilon_\n \leq \mathfrak b $ for some $ {\mathfrak b>0 }$ fixed and small enough, which is of course the only asymptotically interesting case, and constants in $ \boldsymbol{\mathcal O}(\cdot) $ will depend on $ \mathfrak b $, but will be independent of $ \delta $ and $ \n $. We solve \hyperref[rhp]{\rhp$_e$} only for~$ e\in\{\alpha_1,\beta_{\n,1},\beta_1\} $ with the understanding that solutions for $ e\in \{ \alpha_2,\alpha_{\n,2},\beta_2 \} $ are constructed similarly.

\subsection[Local parametrix at alpha\_1]{Local parametrix at $ \boldsymbol{\alpha_1 }$}

Observe that the principal branch \smash{$ \zeta_{\n,\alpha_1}^{1/2}(z) $} is positive on $ (-\infty,\alpha_1) $, see Lemma~\ref{lem:4.1}. Since $ \Phi_\n(\z) $ has a pole at $ \infty^{(0)} $ and a zero at $ \infty^{(1)} $, it follows from \eqref{Phin2} and \eqref{4.1.1} that
\begin{equation}
\label{rhpsi-1}
\exp\bigl\{4|\n|\zeta_{\n,\alpha_1}^{1/2}(z) \bigr\} = \Phi_\n^{(0)}(z)/\Phi_\n^{(1)}(z), \qquad z\in U_{\n,\alpha_1}.
\end{equation}
According to \eqref{Uc1}, the square $ U_{\n,\alpha_1} $ contains a disk of radius at least $ K\delta c(\n)/\sqrt 2 $ centered at~$ \alpha_1 $. Then Lemma~\ref{lem:4.1} and \eqref{koebe} yield that
\begin{equation}
\label{rhpsi-2}
 A_{\alpha_1}^* \bigl(\delta n_1^2\bigr) \leq |\n|^2 \min_{s\in\partial U_{\n,\alpha_1}}\big|\zeta_{\n,\alpha_1}(s)\big|,
\end{equation}
where $A_{\alpha_1}^*:=A_{\alpha_1}K/\sqrt 2$. Then \hyperref[rhwpsi]{\rhwpsi}(a--c) and \eqref{rhpsi-1} imply that
\begin{equation}
\label{rhpsi-3}
\boldsymbol P_{\alpha_1}(z) := \boldsymbol E_{\alpha_1}(z) \mathrm{T}_1 \big[ \boldsymbol\Psi_* \bigl(|\n|^2\zeta_{\n,\alpha_1}(z)\bigr) \exp\bigl\{-2|\n|\zeta_{\n,\alpha_1}^{1/2}(z) \sigma_3\bigr\} \rho_1^{-\sigma_3/2}(z) \big] \boldsymbol D(z)
\end{equation}
satisfies \hyperref[rhp]{\rhp$_{\alpha_1}$}(a--c) for any holomorphic prefactor $\boldsymbol E_{\alpha_1}(z)$. Using \hyperref[rhn]{\rhn}(b) and the definition of $ \boldsymbol N(z) $ in terms of $ \boldsymbol M(z) $ after \eqref{matrix-M}, which implies that these matrices obey exactly the same jump relations, one can readily check that
\begin{equation}
\label{rhpsi-4}
\boldsymbol E_{\alpha_1}(z) := \boldsymbol M(z)\mathrm{T}_1\big[ \bigl(\sigma_3\boldsymbol K\sigma_3\bigr) \bigl(|\n|^2\zeta_{\n,\alpha_1}(z)\bigr) \rho_1^{-\sigma_3/2}(z)\big]^{-1}
\end{equation}
is holomorphic in $U_{\n,\alpha_1}\setminus\{\alpha_1\}$. Since the first and second columns of $\boldsymbol M(z)$ have at most quarter root singularities at $ \alpha_1 $ and the third one is bounded, see Lemma~\ref{lem:Mbounds1}, $\boldsymbol E_{\alpha_1}(z)$ is in fact holomorphic in $U_{\n,\alpha_1}$. Therefore, it follows from \hyperref[rhwpsi]{\rhwpsi}(d) and \eqref{rhpsi-2} that
\[
\boldsymbol P_{\alpha_1}(s)= \boldsymbol M(s) \bigl( \boldsymbol I + \boldsymbol{\mathcal O}\bigl(\delta^{-1/2} n_1^{-1} \bigr) \bigr) \boldsymbol D(s), \qquad s\in \partial U_{\n,\alpha_1},
\]
where $\boldsymbol{\mathcal O}(\cdot) $ is independent of $ \n $ and $ \delta $, but does depend on $ \mathfrak b $, which needs to be small enough so that \hyperref[rhwpsi]{\rhwpsi}(d) is applicable. Recall also that $ \det \boldsymbol M(z) \equiv \det \boldsymbol D(z) \equiv 1 $ as explained between~\eqref{matrix-CD} and \eqref{M-inverse}. Hence, it holds that $ \det \boldsymbol E_{\alpha_1}(z) \equiv 1/\sqrt2 $ and respectively $ {\det \boldsymbol P_{\alpha_1}(z)\equiv 1 }$.

\subsection[Local parametrix at beta\_1 when c]{Local parametrix at $\boldsymbol{\beta_1} $ when $\boldsymbol{c^\prime(\delta)\leq c(\n)} $}

In this case, $ \zeta_{\n,\beta_1}(z) $ is given by the second line in \eqref{b-maps2}, i.e., by \eqref{4.4.1}. Lemma~\ref{lem:4.4} yields that the principal square root branch of this map is positive for $ x>\beta_1 $ (within the domain of conformality). Thus, we can deduce from continuity in the parameter $ c $ as well as \eqref{4.3.4} that this branch is given by the expression within the parenthesis in \eqref{4.4.1}. Hence, it follows from \eqref{Phin2} that
\[
\exp\bigl\{4|\n|\zeta_{\n,\beta_1}^{1/2}(z) \bigr\} = \Phi_\n^{(0)}(z)/\Phi_\n^{(1)}(z), \qquad z\in U_{\n,\beta_1}.
\]
Recall that $ \sqrt{z_c-\beta_1} \geq (K_2/K_1)\delta $ for $ c\geq c^\prime(\delta) $ by the very definition of $ c^\prime(\delta) $ just before \eqref{min-con}. Since $ c^*<c(\n) $, the square $ U_{\n,\beta_1} $ contains a disk of radius at least $ K\delta c^*/\sqrt 2 $ centered at $ \beta_1=\beta_{\n,1} $, see \eqref{Uc1}. Then it follows from Lemma~\ref{lem:4.4} and \eqref{koebe} that
\begin{equation}
\label{rhpsi-6}
A_{\beta_1}^* \bigl( \delta^3 |\n|^2\bigr) \leq |\n|^2 \min_{s\in\partial U_{\n,\beta_1}} \big|\zeta_{\n,\beta_1}(s)\big|,
\end{equation}
where $ A_{\beta_1}^* := c^*A_{\beta_1} K (K_2/K_1)^2/\sqrt2 $. Similarly to \eqref{rhpsi-3}--\eqref{rhpsi-4}, a solution of \hyperref[rhp]{\rhp$_{\beta_1}$} is given~by
\begin{gather*}
\boldsymbol P_{\beta_1}(z) := \boldsymbol E_{\beta_1}(z) \mathrm{T}_1\big[ \boldsymbol\Psi\bigl(|\n|^2\zeta_{\n,\beta_1}(z)\bigr) \exp\bigl\{-2|\n|\zeta_{\n,\beta_1}^{1/2}(z)\sigma_3 \bigr\} \rho_1^{-\sigma_3/2}(z) \big] \boldsymbol D(z), \\
\boldsymbol E_{\beta_1}(z) := \boldsymbol M(z)\mathrm{T}_1\big[ \boldsymbol K \bigl(|\n|^2\zeta_{\n,\beta_1}(z)\bigr) \rho_1^{-\sigma_3/2}(z)\big]^{-1}.
\end{gather*}
It again holds that $ \det \boldsymbol P_{\beta_1}(z) \equiv 1 $. We also get from \hyperref[rhpsi]{\rhpsi}(d) and \eqref{rhpsi-6} that the error term in \hyperref[rhp]{\rhp$_{\beta_1}$}(d) is of order $ \boldsymbol{\mathcal O}\left(\delta^{-3/2} |\n|^{-1} \right) $ with constants independent of $ \delta $ and $ \n $ but dependent on $ \mathfrak b $, which needs to be small enough so that \hyperref[rhpsi]{\rhpsi}(d) is applicable.

\subsection[Local parametrix at beta\_1 when c\^* < c(n) < c\^prime(delta)]{Local parametrix at $\boldsymbol{\beta_1} $ when $ \boldsymbol{c^* < c(\n) < c^\prime(\delta) }$}

In this case $ \zeta_{\n,\beta_1}(z) $ is defined by the first line of \eqref{b-maps2}, i.e., \eqref{4.3.4}. Hence, it follows from \eqref{Phin2} and \eqref{4.3.4} that
\begin{equation}
\label{rhphi-1}
\exp\left\{ -\frac43 \bigl( \zeta_{\n,\beta_1}^{3/2}(z) - \zeta_{\n,\beta_1}(\beta_1+\epsilon_\n) \zeta_{\n,\beta_1}^{1/2}(z) \bigr) \right\} = \Phi_{\n}^{(0)}(z)/\Phi_{\n}^{(1)}(z), \qquad z\in U_{\n,\beta_1}.
\end{equation}
In the considered region of the parameter $ c(\n) $ each square $ U_{\n,\beta_1} $ is defined in \eqref{Uc1} with $ K=1 $. Thus, it follows directly from \eqref{zhat-b} and \eqref{min-con} that
\begin{equation}
\label{rhphi-2}
|\n|^{2/3}\min_{s\in\partial U_{\n,\beta_1}}|\zeta_{\n,\beta_1}(s)| \geq
\begin{cases}
|\n|^{2/3}\zeta_{\n,\beta_1}(\beta_1+\epsilon_\n) =: -s_\n, \\
K_2 \delta|\n|^{2/3}.
\end{cases}
\end{equation}
Then, as in the previous two subsections, one can verify using \eqref{rhphi-1}--\eqref{rhphi-2}, \hyperref[rhwphi]{\rhwphi}, and \hyperref[rhn]{\rhn}(b) that \hyperref[rhp]{\rhp$_{\beta_1}$} is solved by
\begin{gather*}
\boldsymbol P_{\beta_1}(z) := \boldsymbol E_{\beta_1}(z) \mathrm{T}_1\big[ \widetilde{\boldsymbol\Phi} \bigl(|\n|^{2/3}\zeta_{\n,\beta_1}(z);s_\n\bigr) \bigl(\Phi^{(0)}_\n/\Phi^{(1)}_\n\bigr)^{-\sigma_3/2}(z)\rho_1^{-\sigma_3/2}(z) \big] \boldsymbol D(z), \\
\boldsymbol E_{\beta_1}(z) := \boldsymbol M(z)\mathrm{T}_1\big[ \boldsymbol K \bigl(|\n|^{2/3}\zeta_{\n,\beta_1}(z)\bigr) \rho_1^{-\sigma_3/2}(z)\big]^{-1},
\end{gather*}
where $\boldsymbol E_{\beta_1}\!(z)$ is holomorphic in $ U_{\n,\beta_1} $ and the error term in \hyperref[rhp]{\rhp}(d) is of order $ \boldsymbol{\mathcal O}\bigl(\delta^{-1/2} |\n|^{-1/3} \bigr) $ with constants independent of $ \delta $ and $ \n $ but dependent on $ \mathfrak b $, which needs to be small enough so that \hyperref[rhwphi]{\rhwphi}(d) is applicable. Observe that if $ c(\n) $'s are separated away from $ c^* $, then $ s_\n \sim -|\n|^{2/3} $ and the error term can be improved to $ \boldsymbol{\mathcal O}\bigl(\delta^{-1/2} |\n|^{-1} \bigr) $. As before, $ \det \boldsymbol P_{\beta_1}(z) \equiv 1 $.

\subsection[Local parametrix at beta\_n,1 when c(n)leq c\^* and zeta\_beta]{Local parametrix at $ \boldsymbol{\beta_{\n,1} }$ when $ \boldsymbol{c(\n)\leq c^* }$ and $ \boldsymbol{\zeta_{\beta_{\n,1}}(\beta_1) \leq \tau_* |\n|^{-2/3}} $}
\label{ss:64}

Recall also that we set $ \zeta_{\n,\beta_1}(z) = \zeta_{\beta_{\n,1}}(z)-\zeta_{\beta_{\n,1}}(\beta_1) $ in \eqref{b-maps1}. We get from Lemma~\ref{lem:4.2} that the principal branch \smash{$ \zeta_{\beta_{\n,1}}^{3/2}(z) $} is positive on $ (\beta_{c,1},\alpha_2) $. It can be readily inferred from \eqref{calH} and the first estimate of Lemma~\ref{lem:4.6} that this branch is equal to the expression within the parenthesis in \eqref{4.2.1}. Thus, it follows from \eqref{Phin2} that
\begin{equation}
\label{rhphi-4}
\exp\left\{-\frac43|\n|\zeta_{\beta_{\n,1}}^{3/2}(z) \right\} = \Phi_\n^{(0)}(z)/\Phi_\n^{(1)}(z), \qquad z\in U_{\n,\beta_{\n,1}}.
\end{equation}
Define $ s_\n := |\n|^{2/3} \zeta_{\beta_{\n,1}}(\beta_1) \in [0,\tau_*] $, where the last conclusion is the restriction placed on $ s_\n $ in this subsection. This restriction, Lemma~\ref{lem:4.2}, and \eqref{koebe} imply that
\[
\big|\beta_{\n,1} - \beta_1 \big| \leq \bigl( c^*)^{1/3} A_{\beta_1}^{-1} \tau_* |\n|^{-2/3}.
\]
In particular, as $ \beta_{c,1} $ is a continuous increasing function of $ c\in[0,c^*] $ with $ \beta_{c^*,1}=\beta_1 $, recall~\eqref{cstars}, it must necessarily hold that $ c(\n)\to c^* $. Thus, by recalling \eqref{Uc1}, we see that uniformly for all~$ |\n| $ large enough the square $ U_{\n,\beta_{\n,1}} $ contains a disk of radius at least $ \delta c^*/2\sqrt 2 $ centered at $ \beta_1 $ (the factor $ 1/2 $ is there to move the center from $ \beta_{\n,1} $ to $ \beta_1 $). Then Lemma~\ref{lem:4.2} and \eqref{koebe} yield that
\begin{equation}
\label{rhphi-6}
\tilde A_{\beta_1} \bigl( \delta |\n|^{2/3} \bigr) \leq |\n|^{2/3} \min_{s\in\partial U_{\n,\beta_{\n,1}}} \big|\zeta_{\n,\beta_1}(s)\big|,
\end{equation}
where $ \tilde A_{\beta_1} := (c^*)^{2/3}A_{\beta_1}/2\sqrt 2 $. As before, one can check using \hyperref[rhphi]{\rhphi}(a--c) and \hyperref[rhn]{\rhn}(b) that
\begin{gather}
\boldsymbol P_{\beta_{\n,1}}(z) := \boldsymbol E_{\beta_{\n,1}}(z) \mathrm{T}_1\big[ \boldsymbol\Phi\bigl(|\n|^{2/3} \zeta_{\n,\beta_1}(z);s_\n\bigr) \exp\left\{\frac23|\n|\zeta_{\beta_{\n,1}}^{3/2}(z)\sigma_3 \right\} \rho_1^{-\sigma_3/2}(z) \big] \boldsymbol D(z), \nonumber\\
\boldsymbol E_{\beta_{\n,1}}(z) := \boldsymbol M(z)\mathrm{T}_1\big[ \boldsymbol K \bigl(|\n|^{2/3}\zeta_{\beta_{\n,1}}(z)\bigr) \rho_1^{-\sigma_3}(z)\big]^{-1},\label{rhphi-5}
\end{gather}
satisfies \hyperref[rhp]{\rhp$_{\beta_{\n,1}}$}(a--c) and that the prefactor $\boldsymbol E_{\beta_{\n,1}}(z)$ is holomorphic in $ U_{\n,\beta_{\n,1}} $ (it is by design that $ \zeta_{\n,\beta_1}(z) $ is used as the argument of $ \boldsymbol \Phi(\zeta) $ and $ \zeta_{\beta_{\n,1}}(z) $ is used everywhere else). Since
\[
\left|\frac{\zeta_{\beta_{\n,1}}(s)}{\zeta_{\n,\beta_1}(s)}\right| \leq 1+ \left|\frac{\zeta_{\beta_{\n,1}}(\beta_1)}{\zeta_{\n,\beta_1}(s)}\right| \leq 1 + \frac{\bigl(\tau_*/\tilde A_{\beta_1}\bigr)}{\delta|\n|^{2/3}},
\]
it holds that
\[
\boldsymbol K^{-1} \bigl(|\n|^{2/3}\zeta_{\beta_{\n,1}}(z)\bigr)\boldsymbol K \bigl(|\n|^{2/3}\zeta_{\n,\beta_1}(z)\bigr) = \boldsymbol I + \boldsymbol{\mathcal O}\bigl( \delta^{-1}|\n|^{-2/3} \bigr),
\]
where $ \boldsymbol{\mathcal O}(\cdot) $ is independent of $ \delta $ and $ \n $ as long as $ \delta |\n|^{2/3} \geq \tau_*/2\tilde A_{\beta_1} $. Assume further that $ \delta |\n|^{2/3} $ is large enough so that \hyperref[rhphi]{\rhphi}(d) takes place. Then \eqref{rhphi-4}, \eqref{rhphi-6}, and \hyperref[rhphi]{\rhphi}(d) imply that \hyperref[rhp]{\rhp$_{\beta_{\n,1}}$}(d) holds with the error term of order $ \boldsymbol{\mathcal O}\left(\delta^{-1} |\n|^{-1/3} \right) $. As in the previous subsections, we point out that $ \det \boldsymbol P_{\beta_{\n,1}}(z) \equiv 1 $.

\subsection[Local parametrix at beta\_n,1 when c(n)leq c\^{}* and zeta\_beta\_n,1(beta\_1)]{Local Parametrix at $ \boldsymbol{\beta_{\n,1}} $ when $ \boldsymbol{c(\n)\leq c^* }$ and $ \boldsymbol{\zeta_{\beta_{\n,1}}(\beta_1) > \tau_* |\n|^{-2/3}} $}
\label{ss:65}

It follows from \eqref{4.2.1} that \eqref{rhphi-4} still holds with $ \zeta_{\n,\beta_1}(z) $ ``replaced'' by $ \zeta_{\beta_{\n,1}}(z) $ as these symbols denote the same function in the considered case, see \eqref{b-maps1}. Definition of $ U_{\n,\beta_{\n,1}} $ in \eqref{Uc1} implies that this square contains a disk of radius at least \smash{$ \delta c(\n)/3\sqrt 2 $}. Therefore, Lemma~\ref{lem:4.2} and \eqref{koebe} yield that
\begin{equation}
\label{rhth-1}
A_{\beta_1}^* \bigl(\delta n_1^{2/3}\bigr) \leq |\n|^{2/3} \min_{s\in\partial U_{\n,\beta_{\n,1}}} \big|\zeta_{\n,\beta_1}(s) \big|,
\end{equation}
where $ A_{\beta_1}^*:=A_{\beta_1}/3\sqrt2 $. Let $ \boldsymbol E_{\beta_{\n,1}}(z) $ be the same as in \eqref{rhphi-5}, where once again we must keep in mind the relabeling of $ \zeta_{\beta_{\n,1}}(z) $ as $ \zeta_{\n,\beta_1}(z) $. Define
\[
\tau_\n := \begin{cases}
 \infty, & \beta_1 \not\in U_{\n,\beta_{\n,1}}, \\
 |\n|^{2/3} \zeta_{\n,\beta_1}(\beta_1), & \text{otherwise}.
 \end{cases}
\]
Under the conditions considered in this subsection, it holds that $ \tau_\n>\tau_* $.

The following paragraph is applicable only if $ \beta_1 \in U_{\n,\beta_{\n,1}} $, i.e., when $ \tau_\n $ is finite. It follows from \eqref{Kdeltac} that $ \beta_1-\beta_{\n,1}\leq 2\delta c(\n)/3\leq 2c(\n)/3$ in this case. Since $ \beta_1-\beta_{c,1} $ is a strictly decreasing function of $ c $ while $ 2c/3 $ is obviously strictly increasing, there exists a unique $ \tilde c<c^* $ such that $ \beta_1-\beta_{\tilde c,1} = 2\tilde c/3$ (so, $ c(\n)\geq \tilde c $ when $ \tau_\n $ is finite). Each map $ \zeta_{\beta_{c,1}}(z) $ is conformal in~$ U(\beta_{c,1},\delta c) $ for $ c\in[\tilde c,c^*] $ and as a family they continuously depend on the parameter $ c $. As~$ \beta_1 $ is separated from $ \partial U_{\n,\beta_{\n,1}} $ by a distance of at least $ \delta \tilde c/3\sqrt2 $ when $ \beta_1\in U_{\n,\beta_{\n,1}} $, it follows from a~compactness argument that
\[
|\zeta_{\n,\beta_1}(\beta_1)-\zeta_{\n,\beta_1}(s)| \geq d>0, \qquad s\in \partial U_{\n,\beta_{\n,1}},
\]
where $ d $ is independent of $ n $ (but does depend on $ \delta $). Since $ |\n|^{2/3}d>1 $ for all $ |\n| $ large enough, we get that $ |\n|^{2/3}\zeta_{\n,\beta_1}(U_{\n,\beta_{\n,1}} ) $ contains $ \{|\zeta-\tau_\n|\leq 1\} $ in its interior.

We now get from \eqref{rhth-1} and \hyperref[rhth]{\rhth} that \hyperref[rhp]{\rhp$_{\beta_{\n,1}}$} is solved by
\begin{equation}
\label{rhth-2}
\boldsymbol P_{\beta_{\n,1}}(z) := \boldsymbol E_{\beta_{\n,1}}(z) \mathrm{T}_1\left[ \boldsymbol\Theta\bigl(|\n|^{2/3}\zeta_{\n,\beta_1}(z);\tau_\n\bigr) \exp\left\{\frac23|\n|\zeta_{\n,\beta_1}^{3/2}(z)\sigma_3 \right\} \rho_1^{-\sigma_3/2}(z) \right] \boldsymbol D(z),
\end{equation}
where the error term in \hyperref[rhp]{\rhp$_{\beta_{\n,1}}$}(d) is of order \smash{$ \boldsymbol{\mathcal O}\bigl(\delta^{-1/2} n_1^{-1/3} \bigr) $} with constants independent of~$ \delta $ and $ \n $, which can be improved to \smash{$ \boldsymbol{\mathcal O}\bigl(\delta^{-1/2} n_1^{-1} \bigr) $} if $ c(\n) $'s are separated from $ c^* $. As in all the previous subsections, we have that $ \det \boldsymbol P_{\beta_{\n,1}}(z) \equiv 1 $.

\subsection[Local parametrix at beta\_1 when]{Local parametrix at $ \boldsymbol{\beta_1} $ when $\boldsymbol{\beta_1 \notin U_{\n,\beta_{\n,1}} }$}

Observe that it necessarily holds in the considered case that $ c(\n) < c^* $. Let
\[
I_1(z) := \int_{\Delta_1}\frac{\rho_1(x)}{x-z}\frac{{\rm d}x}{2\pi\ic} - \frac{\rho_1(\alpha_1)}{2\pi\ic}\log(z-\alpha_1), \qquad z\notin(-\infty,\beta_1],
\]
where we use the principal branch of the logarithm ($\log(x-\alpha_1)>0 $ for $ x>\alpha_1 $). $ I_1(z) $ is analytic in the domain of its definition, has continuous traces on $ \Delta_1^\circ $ that satisfy $ (I_{1+} - I_{1-})(x) = \rho_1(x) $ according to Plemelj--Sokhotski formulae, see \cite[Section~I.4.2]{Gakhov}, it has a logarithmic singularity at $ \beta_1 $ and is bounded in the vicinity of $ \alpha_1 $, see \cite[Section~I.8.2]{Gakhov}. Due to the construction of~$ U_{\n,\beta_1} $ in \eqref{Uc2}, $ \partial U_{\n,\beta_1} $ never approaches $ \beta_1 $ and hence, $ |I_1(s)| $ is uniformly bounded on $ \partial U_{\n,\beta_1} $ with the bound dependent on the chosen value of $ r $.

Now, it follows from \eqref{calH}, the first and second claims of Lemma~\ref{lem:4.6}, and the definition of~$ U_{\n,\beta_1} $ in \eqref{Uc2} that
\begin{equation}
\label{rhc-1}
\big|\Phi_\n^{(0)}(z)/\Phi_\n^{(1)}(z)\big| = {\rm e}^{|\n| \mathcal H_\n(z)} \leq {\rm e}^{-B(\delta/3)^{3/2}n_1} = \mathcal O\bigl( \delta^{-3/2} n_1^{-1}\bigr), \qquad z\in \overline U_{\n,\beta_1}
\end{equation}
(here we are using the facts that the second estimate of Lemma~\ref{lem:4.6} is needed only when $ c(\n) < c_r $ and we always can make $ c_r $ small enough so that $ Bc_r < 3^{3/2}\epsilon(r) $). Hence, a solution of \hyperref[rhp]{\rhp$_{\beta_1}$} in this case is given by
\begin{gather}
\boldsymbol P_{\beta_1}(z) := \boldsymbol M(z) \bigl( \boldsymbol I + I_{\n,1}(z) \boldsymbol E_{1,2} \bigr) \boldsymbol D(z),\nonumber \\
I_{\n,1} (z) := I_1(z)\Phi_\n^{(0)}(z)/\Phi_\n^{(1)}(z),
\qquad z\in U_{\n,\beta_1}\setminus\Delta_1,\label{rhc-2}
\end{gather}
Indeed, as matrices $ \boldsymbol M(z) $ and $ \boldsymbol D(z) $ are holomorphic in $ U_{\n,\beta_1} $ and $ I_1(z) $ is holomorphic in $ U_{\n,\beta_1}\setminus\Delta_1 $, the above matrix satisfies \hyperref[rhp]{\rhp$_{\beta_1}$}(a). Requirement \hyperref[rhp]{\rhp$_{\beta_1}$}(b) easily follows from the form of the additive jump of $ I_1(z) $. Since $ I_1(z) $ appears only in the second column of~$ \boldsymbol P_{\beta_1}(z) $ and has a logarithmic singularity at $ \beta_1 $, \hyperref[rhp]{\rhp$_{\beta_1}$}(c) is fulfilled. Finally, \hyperref[rhp]{\rhp$_{\beta_1}$}(d) is a consequence of~\eqref{rhc-1}.

\section{Small norm problem}

\begin{figure}[!ht]
\centering
\includegraphics[scale=.9]{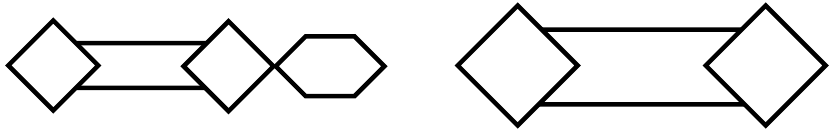}
\begin{picture}(370,0)
\put(326,40){$\partial U_{\n,\beta_2}$}
\put(210,40){$\partial U_{\n,\alpha_2}$}
\put(-5,20){$\partial U_{\n,\alpha_1}$}
\put(110,20){$\partial U_{\n,\beta_{\n,1}}$}
\put(130,60){$\partial U_{\n,\beta_1}$}
\put(48,58){$\Gamma_{\n,1}^+\setminus \overline U_\n$}
\put(48,19){$\Gamma_{\n,1}^-\setminus \overline U_\n$}
\put(265,64){$\Gamma_{\n,2}^+\setminus \overline U_\n$}
\put(265,12){$\Gamma_{\n,2}^-\setminus \overline U_\n$}
\end{picture}
\caption{\small Lens $ \Sigma_{\n,\delta} $ consisting of two connected components $ \Sigma_{\n,\delta,1} $ (the left one) and $ \Sigma_{\n,\delta,2} $ (the right one).}
\label{fig:sigma-lens}
\end{figure}
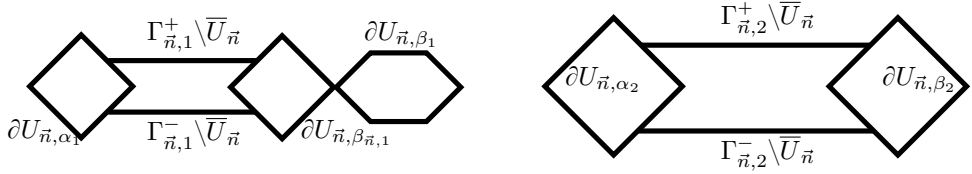

In this section, we make the last preparatory step before solving \hyperref[rhx]{\rhx}. Recall \eqref{Ipm}. Set
\[
\Sigma_{\n,\delta} := \partial U_\n \cup \bigl( \Gamma_\n \setminus \overline U_\n \bigr), \qquad U_\n := \bigcup_e U_{\n,e}, \qquad \Gamma_\n := \bigcup_{i=1}^2 \bigl(\Gamma_{\n,i}^+ \cup \Gamma_{\n,i}^-\bigr),
\]
see Figure~\ref{fig:sigma-lens}. The parts of $ \Sigma_{\n,\delta} $ that belong to $ \Gamma_\n $ inherit their orientation from the original arcs and the individual polygons in $ \partial U_\n $ are oriented clockwise. We shall further denote by $ \Sigma_{\n,\delta,1} $ and $ \Sigma_{\n,\delta,2} $ the left and right, respectively, connected components of $ \Sigma_{\n,\delta} $ and by $\Sigma_{\n,\delta}^\circ$ the subset of points around which $ \Sigma_{\n,\delta} $ is locally a Jordan arc.

Given the global parametrix $\boldsymbol N(z)=\boldsymbol C(\boldsymbol{MD})(z)$ solving \hyperref[rhn]{\rhn}, see \eqref{matrix-M} and \eqref{matrix-CD}, and local parametrices $\boldsymbol P_e(z)$ solving \hyperref[rhp]{\rhp$_e$} and constructed in the previous section, consider the following Riemann--Hilbert Problem (\rhz):
\begin{itemize}\itemsep=0pt
\label{rhz}
\item[(a)] $\boldsymbol{Z}(z)$ is a holomorphic matrix function
 in $\overline\C\setminus\Sigma_{\n,\delta}$ and
 $\boldsymbol{Z}(\infty)=\boldsymbol{I}$;
\item[(b)] $\boldsymbol{Z}(z)$ has continuous and bounded boundary values on $\Sigma_{\n,\delta}^\circ$ that satisfy
\[
\boldsymbol{Z}_+(s) = \boldsymbol{Z}_-(s)
\begin{cases}
(\boldsymbol{MD})(s) (\boldsymbol I+\rho_i^{-1}(s)\boldsymbol E_{1,i+1}) (\boldsymbol{MD})^{-1}(s), & s\in\Gamma_\n\setminus \overline U_\n, \\
\boldsymbol P_e(s) (\boldsymbol{MD})^{-1}(s), & s\in\partial U_{\n,e},
\end{cases}
\]
where $ e \in \{ \alpha_1,\beta_{\n,1},\beta_1,\alpha_2,\alpha_{\n,2},\beta_2 \} $.
\end{itemize}

Then the following lemma takes place.

\begin{Lemma}
\label{lem:rhz}
For each $ \delta\in(0,\delta_\Delta/2) $, there exists $ \varepsilon_\delta>0 $ such that \hyperref[rhz]{\rhz} is solvable for all~${ \varepsilon_\n\leq \varepsilon_\delta }$. Moreover, for each $ r>0 $ small enough there exists a constant $ C_{2,r} $, independent of~$ \delta $ and $ \n $, such that\footnote{In particular, the estimate holds in $ \overline\C $ when $ c^*\leq c(\n) \leq c^{**} $.}
\begin{equation}
\label{rhz-1}
| Z_{i,k}(z)| \leq C_{2,r}\frac{\varepsilon_\n^{1/3}}{\delta^5}, \qquad
\begin{cases}
|z-\alpha_1|\geq 2r, & c(\n)\leq c^{\ast\ast}, \\
|z-\beta_2| \geq 2r, & c(\n)\geq c^\ast,
\end{cases}
\end{equation}
$ i,k\in\{0,1,2\} $, where, for the ease of the future use, we let $ Z_{i,k}(z) := [\boldsymbol Z(z)]_{i+1,k+1} - \delta_{ik} $, $ \delta_{ik} $ is the Kronecker symbol, and $ Z_{i,k}(z) $ needs to be replaced by $ Z_{i,k\pm}(z) $ for $ z\in\Sigma_{\n,\delta} $. The exponent~$ 1/3 $ of $ \varepsilon_\n $ can be replaced by $ 1 $ if we additionally require that $ |c(\n)-c^\ast|,|c(\n)-c^{\ast\ast}|\geq\epsilon>0 $, where the constants $ \varepsilon_\delta $ and $ C_{2,r} $ will depend on $ \epsilon $ in this case.
\end{Lemma}
\begin{proof}
Let us more generally consider \hyperref[rhz]{\rhz} on $ \Sigma_{\n,\nu\delta} $, where $ \nu\in[1/2,2] $ and we also scale the parameter $ r $ by $ \nu $, see \eqref{Uc2}. Put
\[
\boldsymbol I + \boldsymbol V(s) := \bigl(\boldsymbol Z_-^{-1}\boldsymbol Z_+\bigr)(s), \qquad s\in \Sigma_{\n,\nu\delta}^\circ,
\]
to be the jump of $ \boldsymbol Z(z) $ on $ \Sigma_{\n,\nu\delta} $. It can be readily seen from analyticity of $ \boldsymbol M(z) $, $ \boldsymbol D(z) $, $ \boldsymbol P_e(z) $, and \hyperref[rhz]{\rhz}(b) that $ \boldsymbol V(s) $ can be analytically continued off each Jordan subarc of $ \Sigma_{\n,\nu\delta} $. Thus, the solutions of \hyperref[rhz]{\rhz} for different values of $ \nu $, if exist, are analytic continuations of each other.

Let us now estimate the size of $ \boldsymbol V(s) $ in the supremum norm. We shall do it only on $ \Sigma_{\n,\nu\delta,1} $, understanding that the estimates on $ \Sigma_{\n,\nu\delta,2} $ can be carried out in a similar fashion. For $ s\in\partial U_{\n,e} $, $ e\in \{\alpha_1,\beta_{\n,1},\beta_1\} $, it holds that
\[
\boldsymbol V(s) = \boldsymbol P_e(s) (\boldsymbol{MD})^{-1}(s) - \boldsymbol I = \boldsymbol M(s) \boldsymbol{\mathcal O} \bigl( \delta^{-3/2} \varepsilon_\n^{1/3} \bigr) \boldsymbol M^{-1}(s) = \boldsymbol{\mathcal O}\bigl( \delta^{-2} \varepsilon_\n^{1/3}\bigr)
\]
by \hyperref[rhp]{\rhp$_e$}(d) and Lemmas~\ref{lem:Mbounds1}--\ref{lem:Mbounds2} (since $ \nu\geq1/2 $, there is no need to explicitly include it in the estimates). Notice also that the power $ 1/3 $ of $ \varepsilon_\n $ comes only from the local problems discussed in Sections~\ref{ss:64}--\ref{ss:65}. When $ |c(\n)-c^\ast|,|c(\n)-c^{\ast\ast}|\geq\epsilon>0 $, the material of Section~\ref{ss:64} is no longer relevant and the estimate in Section~\ref{ss:65} is of order $ \varepsilon_\n $ as remarked after \eqref{rhth-2}. For~$ s\in \Gamma_{\n,1}^\pm\setminus U_\n $, it follows from the third inequality in Lemma~\ref{lem:4.6} and Lemmas~\ref{lem:Mbounds1}--\ref{lem:Mbounds2} that
\begin{align*}
\boldsymbol V(s) &= (\boldsymbol {MD})(s) \bigl(\boldsymbol I+\rho_1^{-1}(s)\boldsymbol E_{1,2}\bigr) (\boldsymbol{MD})^{-1}(s) - \boldsymbol I \\
& = \rho_1^{-1}(s) \bigl( \Phi_\n^{(1)}(s)/\Phi_\n^{(0)}(s) \bigr) \boldsymbol M(s)\boldsymbol E_{2,1}\boldsymbol M^{-1}(s) = \boldsymbol{\mathcal O}\bigl( \delta^{-7/2}\varepsilon_\n \bigr).
\end{align*}
Altogether, we have shown that $ \|\boldsymbol V\|_{\Sigma_{\n,\nu\delta}} = \mathcal O\bigl( \delta^{-4}\varepsilon_\n ^{1/3} \bigr) $, where $ \mathcal O(\cdot) $ is independent of $ \delta $ and $ \n $, and that $ \|\boldsymbol V\|_{\Sigma_{\n,\nu\delta}} = \mathcal O_\epsilon\bigl( \delta^{-4}\varepsilon_\n \bigr) $ when $ |c(\n)-c^\ast|,|c(\n)-c^{\ast\ast}|\geq\epsilon $.

It was explained in \cite[Lemma~7.9]{ApDenYa21} that the norms of Cauchy operators (functions are mapped into traces of their Cauchy integrals) as operators from $ L^2(\Sigma_{\n,\nu\delta}) $ into itself are uniformly bounded above independently of $ \n $ and $ \nu\delta $. Then, as in Section~\ref{sec:MLP}, \cite[Theorem~8.1]{FokasItsKapaevNovokshenov} allows us to conclude that \hyperref[rhz]{\rhz} is solvable for all $ \varepsilon_\n \leq \varepsilon_*\delta^{12} $ and
\begin{equation}
\label{rhz-4}
\|\boldsymbol Z_\pm\|_{L^2(\Sigma_{\n,\nu\delta})} \leq C^\prime \delta^{-4} \varepsilon_\n^{1/3} \qquad \text{or} \qquad \|\boldsymbol Z_\pm\|_{L^2(\Sigma_{\n,\nu\delta})} \leq C^\prime_\epsilon \delta^{-4}\varepsilon_\n,
\end{equation}
where $ C^\prime$, $C^\prime_\epsilon $ are independent of $ \n $ and $ \delta $ and the second estimate holds when $ |c(\n)-c^\ast|,|c(\n)-c^{\ast\ast}|\geq\epsilon $.

Recall now that the squares $ U_{\n,\alpha_1} $ and $ U_{\n,\beta_{\n,1}} $ have diameters that are at most $ 2\nu\delta c(\n) $ and at least $ \min\{1/3,K^\prime\}\delta c(\n)$ long, see \eqref{Uc1}, while the narrowest part of $ U_{\n,\beta_1} $ is similarly proportioned, see \eqref{Uc2}. Moreover, analogous claims hold for $ U_{\n,\alpha_2} $ and $ U_{\n,\alpha_{\n,2}} $, and $ U_{\n,\beta_2} $ with $ c(\n) $ replaced by $ 1-c(\n) $. Hence, one can choose a finite collection of values for the parameter $ \nu $ (the values $1/2,1,2$ should do the job) to make sure that there exists a constant~$ \gamma\in(0,1)$ such that every $ z $ is at least distance $ \gamma\delta \varsigma_\n $ away from one of the lenses $ \Sigma_{\n,\nu\delta} $, where~$ {\varsigma_\n = \min\{ c(\n),1-c(\n) \} }$. Since the arclengths of $ \Sigma_{\n,\nu\delta} $ are uniformly bounded above, it readily follows from the Cauchy integral formula, a straightforward estimate, and the H\"older inequality that
\begin{equation}
\label{rhz-2}
|Z_{i,k}^\nu(z)| \leq C^{\prime\prime} \bigl( \|\boldsymbol Z_+\|_{L^2(\Sigma_{\n,\nu\delta})} + \|\boldsymbol Z_-\|_{L^2(\Sigma_{\n,\nu\delta})} \bigr) / (\gamma \varsigma_\n\delta ),
\end{equation}
for all $ \dist(z,\Sigma_{\n,\nu\delta}) \geq \gamma\varsigma_\n\delta $ and $ i,k\in\{0,1,2\} $, and some constant $ C^{\prime\prime} $, independent of $ \n $ and $ \delta $, where the superscript $ \nu $ signifies that these functions come from $ \boldsymbol Z(z) $ with jump on $ \Sigma_{\n,\nu\delta} $. As matrices $ \boldsymbol{Z}(z) $ for different values of $ \nu $ are analytic continuations of each other, inequality \eqref{rhz-2} can be improved to
\begin{equation}
\label{rhz-3}
|Z_{i,k}(z)| \leq C^{\prime\prime} \max_\nu \bigl( \|\boldsymbol Z_+\|_{L^2(\Sigma_{\n,\nu\delta})} + \|\boldsymbol Z_-\|_{L^2(\Sigma_{\n,\nu\delta})} \bigr) / (\gamma \varsigma_\n\delta)
\end{equation}
for all $ z\in\C $ and $ i,k\in\{0,1,2\} $, where $ Z_{i,k}(z) $ needs to be replaced by $ Z_{i,k\pm}(z) $ for $ z\in\Sigma_{\n,\delta} $.

To prove \eqref{rhz-1}, assume that $ c(\n)\leq c^{\ast\ast} $ (the cases $ c^\ast\leq c(\n)$ can be handled similarly). Recall the definition of the polygon $ \partial U_{\n,\beta_1} $ in \eqref{Uc2}, see also Figure~\ref{fig:polygon}\,(b). If we only consider indices $ \n $ such that $ c(\n) \geq c_r $, then \eqref{rhz-4} and \eqref{rhz-3} clearly yield \eqref{rhz-1}. When $ c(\n)<c_r $, the part of $ \Sigma_{\n,\nu\delta,1} $ that lies in $ \{|z-\alpha_1|\geq2\nu r\} $ only depends on $ \nu r $ (this must be a part of the boundary of the second union in \eqref{Uc2}) and every point of $ \{|z-\alpha_1|\geq2r\} $ lies distance at least $ \gamma \delta $ from one of the sets $ \Sigma_{\n,\nu\delta,2} $, $ \nu\in[1/2,1] $ (perhaps at expense of decreasing $ \gamma $). Thus, similarly to \eqref{rhz-2}, it holds that
\[
|Z_{i,k}^\nu(z)| \leq C^{\prime\prime}_r ( \|\boldsymbol Z_+\|_{L^2(\Sigma_{\n,\nu\delta})} + \|\boldsymbol Z_-\|_{L^2(\Sigma_{\n,\nu\delta})} ) /(\gamma \delta )
\]
for $ \dist(z,\Sigma_{\n,\nu\delta}) \geq \gamma\delta $ and $ |z-\alpha_1|\geq 2\nu r $, and $ i,k\in\{0,1,2\} $. We now get \eqref{rhz-1} from \eqref{rhz-4} and the appropriate analog of \eqref{rhz-3} that is obtained from the above estimate by varying $ \nu\in[1/2,1] $ over a finite set of values.
\end{proof}

\section[Proofs of Theorems 1.1 and 1.6]{Proofs of Theorems~\ref{thm:1} and \ref{thm:2}}

The proofs of all the main results are based on the following lemma, which is an immediate consequence of Lemma~\ref{lem:rhz}.

\begin{Lemma}
\label{lem:rha}
A solution of \hyperref[rhx]{\rhx} exists for all $ \varepsilon_\n $ small enough and is given by
\[
 \boldsymbol X(z) := \boldsymbol C\boldsymbol Z(z)
 \begin{cases}
(\boldsymbol{MD})(z), &z\in \overline\C\setminus\overline U_\n, \\
\boldsymbol P_e(z), & z\in U_{\n,e}, e\in\{\alpha_1,\beta_{\n,1},\beta_1,\alpha_2,\alpha_{\n,2},\beta_2\},
\end{cases}
 \]
 where $ \boldsymbol Z(z) $ solves \hyperref[rhz]{\rhz}, $ \boldsymbol N(z) = \boldsymbol C(\boldsymbol{MD})(z) $ solves \hyperref[rhn]{\rhn}, and $ \boldsymbol P_e(z) $ solve \hyperref[rhp]{\rhp$_e$}.
\end{Lemma}

\subsection[Proof of Theorem 1.1]{Proof of Theorem~\ref{thm:1}}

Let $ K_{\n,d} := \{z\in\overline\C : \dist(z,\Delta_\n)\geq d\} $. We can choose parameter $ \delta $ in the definition of the contour $ \Sigma_{\n,\delta} $ so that $ K_{\n,d} $ does not intersect the closures of those connected components of the complement of $ \Sigma_{\n,\delta} $ that intersect each $ \Omega_{\n,i}^\pm $, $ i\in\{1,2\} $, see Figures~\ref{fig:lens} and~\ref{fig:sigma-lens}, and recall \eqref{Uc2}. However, $ K_{\n,d} $ can intersect $ U_{\n,\beta_1} $ or $ U_{\n,\alpha_2 }$ when $ \beta_1\not\in U_{\n,\beta_{\n,1}} $ or $ \alpha_2\not\in U_{\n,\alpha_{\n,2}} $, respectively (these two things cannot happen simultaneously). We also assume that $ 2r<d $ in \eqref{rhz-1}.

Recall the definition of $ I_{\n,1}(z) $ in \eqref{rhc-2}. We define these functions to be non-zero only in~$ U_{\n,\beta_1} $ and only when $ \beta_1\not\in U_{\n,\beta_{\n,1}} $. The functions $ I_{\n,2}(z) $ are defined similarly in $ U_{\n,\alpha_2 } $ and only when~$ \alpha_2\not\in U_{\n,\alpha_{\n,2}} $. Then
\begin{equation}
\label{Y1}
\boldsymbol Y(z) = \boldsymbol C (\boldsymbol{ZM})(z)\bigl(\boldsymbol I+ I_{\n,1}(z)\boldsymbol E_{1,2}+ I_{\n,2}(z)\boldsymbol E_{1,3}\bigr)\boldsymbol D(z), \qquad z\in K_{\n,d},
\end{equation}
by \eqref{X}, \eqref{rhc-2}, and Lemma~\ref{lem:rha}. To be more precise, we need to write $ \boldsymbol Y_\pm(z) $ and $ I_{\n,i\pm}(z) $ for~${ z\in \Delta \setminus \Delta_\n }$ when this set is non-empty as well as $ \boldsymbol Z_\pm(z) $ for $ z\in K_{\n,d}\cap U_{\n,\beta_1} $ or $ z\in K_{\n,d}\cap U_{\n,\alpha_2} $ when $ \beta_1\notin U_{\n,\beta_{\n,1}} $ or $ \alpha_2\notin U_{\n,\alpha_{\n,2}} $, respectively. With the notation of Lemma~\ref{lem:rhz}, we get from~\eqref{Y}, \eqref{Y1}, the definition of $ \boldsymbol M(z) $ in \eqref{matrix-M} and of $ \boldsymbol C $, $ \boldsymbol D(z) $ in \eqref{matrix-CD} that
\begin{align}
P_\n(z) &= [\boldsymbol Y(z)]_{1,1} = [\boldsymbol C]_{1,1} [(\boldsymbol{ZM})(z)]_{1,1}[\boldsymbol D(z)]_{1,1} \nonumber \\
& = \bigl( 1 + Z_{0,0}(z) + s_{\n,1}\Upsilon_{\n,1}^{(0)}(z)Z_{0,1}(z) + s_{\n,2}\Upsilon_{\n,2}^{(0)}(z)Z_{0,2}(z) \bigr) \mathcal P_\n(z),
\label{Pn}
\end{align}
where one needs to recall \eqref{Fn} and we set \smash{$ s_{\n,i}:=S_\n^{(0)}(\infty)/S_\n^{(i)}(\infty) $}, $ i\in\{1,2\} $, as well as observe that $ [\boldsymbol C]_{1,1} = 1/\tau_\n $, see \eqref{Phin}. Estimates \eqref{6.2.6} imply that
\begin{equation}
\label{Up-bounds1}
c(\n)^{-1}\big|\Upsilon_{\n,1\pm}^{(0)}(x)\big| \lesssim 1 \qquad \text{and}\qquad  (1-c(\n))^{-2}\big|\Upsilon_{\n,2\pm}^{(0)}(x)\big| \lesssim 1
\end{equation}
for $ x\in \Delta_{\n,1} $, where the constants in $ \lesssim $ are independent of $ \n $. It essentially follows from symmetry and was rigorously shown in \cite[Lemma~5.2]{ApDenYa21} that
\begin{equation}
\label{Up-bounds2}
c(\n)^{-2}\big|\Upsilon_{\n,1\pm}^{(0)}(x)\big| \lesssim 1 \qquad \text{and}\qquad  (1-c(\n))^{-1}\big|\Upsilon_{\n,2\pm}^{(0)}(x)\big| \lesssim 1
\end{equation}
for $ x\in \Delta_{\n,2} $. Moreover, we get from the first estimate in \eqref{6.2.1} that
\begin{equation}
\label{est_sni}
s_{\n,1} \sim c(\n)^{-1} \qquad \text{and}\qquad  s_{\n,2} \sim (1-c(\n))^{-1},
\end{equation}
where the constants in $ \sim $ are independent of $ \n $. Hence, we deduce from the maximum modulus principle that
\begin{equation}
\label{Up-bounds}
\big|s_{\n,1}\Upsilon_{\n,1}^{(0)}(z)\big|, \big|s_{\n,2}\Upsilon_{\n,2}^{(0)}(z)\big| \lesssim 1
\end{equation}
uniformly for all $ \n $ and $ z $ in the extended complex plane, including the traces on $ \Delta_\n $. Plugging~\eqref{Up-bounds} and \eqref{rhz-1} into \eqref{Pn} yields
$
P_\n(z) = (1+o(1)) \mathcal P_\n(z)
$
uniformly on $ K_{\n,d} $, where the error terms are exactly as described in the statement of the theorem. Recall now that $ \mathcal P_\n(z) = \mathcal P_{\n,1}(z)\mathcal P_{\n,2}(z) $. Let $ \Gamma_{\n,d,i}:= \{ z : \dist(z,\Delta_{\n,i}) = d \} $, ${ i\in\{1,2\} }$. Assume that $ d $ is small enough so that these two curves have disjoint interiors. It is easy to see from their very definitions that~$ \mathcal P_{\n,i}(z) $ has winding number $ n_i $ on $ \Gamma_{\n,d,i} $ (it is analytic and non-vanishing in the exterior of~$ \Gamma_{\n,d,i} $ with a pole of order $ n_i $ at infinity) while $ \mathcal P_{\n,3-i}(z) $ has winding number $ 0 $ (it is analytic and non-vanishing in the interior of~$ \Gamma_{\n,d,i} $). Hence, $ \mathcal P_\n(z) $ has winding number $ n_i $ on $ \Gamma_{\n,d,i} $ and so does $ P_\n(z) $ for all $ \varepsilon_\n $ sufficiently small. For all such~$ \n $, let~$ P_{\n,i}(z) $ be the monic subfactor of $ P_\n(z) $ of degree $ n_i $ that has zeros only in the interior of~$ \Gamma_{\n,d,i} $. Then~${ P_{\n,i}(z)/\mathcal P_{\n,i}(z) }$ is a holomorphic and non-vanishing function in the exterior of~$ \Gamma_{\n,d,i} $ that assumes value $ 1 $ at infinity. As this exterior is simply connected, we have that
$
f_{\n,i}(z) := \log \bigl( P_{\n,i}(z)/\mathcal P_{\n,i}(z) \bigr)
$
admits a holomorphic branch in the exterior of $ \Gamma_{\n,d,i} $ that vanishes at infinity. Moreover, $ f_{\n,1}(z) + f_{\n,2}(z) = \log(1+o(1)) = o(1) $ uniformly in $ K_{\n,d} $. Of course, this is true for any $ d>0 $. Thus, we deduce from the Cauchy integral formula that
\[
f_{\n,i}(z) = \int_{\Gamma_{\n,d/2,i}} \frac{f_{\n,i}(s)}{z-s}\frac{{\rm d}s}{2\pi\ic} = \int_{\Gamma_{\n,d/2,i}} \frac{o(1) - f_{\n,3-i}(s)}{z-s}\frac{{\rm d}s}{2\pi\ic} = \int_{\Gamma_{\n,d/2,i}} \frac{o(1)}{z-s}\frac{{\rm d}s}{2\pi\ic}
\]
for $ z $ in the exterior of $ \Gamma_{\n,d/2,i} $, where we used analyticity of $ f_{\n,3-i}(z) $ in the interior of $ \Gamma_{\n,d/2,i} $ on the last step. A trivial estimate now yields that $ f_{\n,i}(z) = o(1) $ in the exterior of $ \Gamma_{\n,d,i} $, which is equivalent to the first claim of \eqref{main1}.

Let now $ F_{\n,d} \subset\Delta_\n^\circ $ be such that $ \dist(F_{\n,d},E_\n)\geq d$. Again, we can adjust $ \delta $ so that $ F_{\n,d} $ does not intersect the closure of $ U_\n $. Hence,
\begin{equation}
\label{Y2}
\boldsymbol Y_{\pm}(x) = \boldsymbol C (\boldsymbol{ZM}_{\pm}\boldsymbol D_{\pm})(x)\bigl(\boldsymbol I \pm \rho_i^{-1}(x) \boldsymbol E_{i+1,1} \bigr), \qquad x\in F_{\n,d}\cap\Delta_{\n,i},
\end{equation}
for $ i\in\{1,2\} $, again by \eqref{X} and Lemma~\ref{lem:rha}. Then we get for $ x\in F_{\n,d} \cap \Delta_{\n,i} $ that
\begin{align*}
P_\n(x) = {}&\frac{(S_\n\Phi_\n)_\pm^{(0)}(x)}{S_\n^{(0)}(\infty) \tau_\n} \bigl(1+Z_{0,0}(x) + s_{\n,1} \bigl(\Upsilon_{\n,1\pm}^{(0)}Z_{0,1}\bigr)(x) + s_{\n,2}\bigl(\Upsilon_{\n,2\pm}^{(0)}Z_{0,2}\bigr)(x) \bigr) \\
 &\pm\frac1{S_\n^{(0)}(\infty) \tau_\n} \\
 &\times\frac{(S_\n\Phi_\n)_\pm^{(i)}(x)}{(\rho_iw_{\n,i\pm})(x)}\bigl(1+Z_{0,0}(x) + s_{\n,1} \bigl(\Upsilon_{\n,1\pm}^{(i)}Z_{0,1}\bigr)(x) + s_{\n,2}\bigl(\Upsilon_{\n,2\pm}^{(i)}Z_{0,2}\bigr)(x) \bigr).
\end{align*}
Since $ G^{(0)}_\pm(x) = G^{(i)}_\mp(x) $ on $ \Delta_{\n,i} $ for any rational function $ G(\z) $ on $ \RS_\n $, it follows from the definition of $ F_\n(z) $, \eqref{szego-pts2}, \eqref{rhz-1}, and \eqref{Up-bounds} that
\[
P_\n(x) = (1+ o(1)) \mathcal P_{\n+}(x) + (1+ o(1)) \mathcal P_{\n-}(x)
\]
uniformly on $ F_{\n,d} $. From this we immediately deduce that
\[
P_{\n,i}(x) = \left[ (1+ o(1)) \mathcal P_{\n,i+}(x) + (1+ o(1)) \mathcal P_{\n,i-}(x) \right] \frac{\mathcal P_{\n,3-i}(x)}{P_{\n,3-i}(x)}
\]
uniformly on $ F_{\n,d}\cap\Delta_{\n,i} $, $ i\in\{1,2\} $. The second asymptotic formula of \eqref{main1} now easily follows from the first.

\subsection[Proof of Theorem 1.6]{Proof of Theorem~\ref{thm:2}}

We retain the notation introduced in the previous subsection. Similarly to the matrix $ \boldsymbol Y(z) $ defined in \eqref{Y}, set
\begin{equation}
\label{hatY}
\widehat{\boldsymbol Y}(z) := \left(\begin{matrix} L_\n(z) & - A_\n^{(1)}(z) & -A_\n^{(2)}(z) \\ -h_{\n,1} L_{\n+\vec e_1}(z) & h_{\n,1} A_{\n+\vec e_1}^{(1)}(z) & h_{\n,1}A_{\n+\vec e_1}^{(2)}(z) \\ -h_{\n,2} L_{\n+\vec e_2}(z) & h_{\n,2} A_{\n+\vec e_2}^{(1)}(z) & h_{\n,2}A_{\n+\vec e_2}^{(2)}(z) \end{matrix}\right),
\end{equation}
where $ L_\n(z) := \int(z-x)^{-1}Q_\n(x) $. It was shown in \cite[Theorem~4.1]{GerKVA01} that
\begin{equation}
\label{hatY1}
\widehat{\boldsymbol Y}(z) = \bigl(\boldsymbol Y^\mathrm{t}(z)\bigr)^{-1},
\end{equation}
where $ \cdot^\mathrm{t} $ denotes the transpose. As before, let $ K_{\n,d} := \{z\in\overline\C : \dist(z,\Delta_\n)\geq d\} $. We keep all the restrictions on $ \delta $ and $ r $ in the definition of $ \Sigma_{\n,\delta} $ and in \eqref{rhz-1}, respectively, specified at the beginning of the previous subsection. It follows from \eqref{Y1} and \eqref{hatY1} that
\[
\widehat{\boldsymbol Y}(z) = \boldsymbol C^{-1}\bigl(\boldsymbol Z^{-1}\bigr)^{\mathrm t}(z)\bigl(\boldsymbol M^{-1}\bigr)^\mathrm{t}(z) \bigl(\boldsymbol I - I_{\n,1}(z)\boldsymbol E_{2,1} - I_{\n,2}(z)\boldsymbol E_{3,1}\bigr) \boldsymbol D^{-1}(z),
\]
$ z\in K_{\n,d} $, where one needs to remember that the functions $ I_{\n,1}(z) $ and $ I_{\n,2}(z) $ are never non-zero simultaneously and we need to replace matrix functions by their boundary values when appropriate as described after \eqref{Y1}. The above equation and \eqref{hatY} yield that
\begin{equation}
\label{Y3}
- A_\n^{(i)}(z) = \tau_\n\big[\bigl(\boldsymbol Z^{-1}\bigr)^{\mathrm t}(z) \bigl(\boldsymbol M^{-1}\bigr)^\mathrm{t}(z) \big]_{1,i+1}/\Phi_\n^{(i)}(z), \qquad z\in K_{\n,d},
\end{equation}
where one needs to remember that $ [\boldsymbol C]_{1,1} = 1/\tau_\n $, see \eqref{Phin} and \eqref{matrix-CD}. Let us rewrite \eqref{M-inverse} as
\[
\boldsymbol M^{-1}(z) = \diag\left(\frac1{S_\n^{(0)}(z)},\frac{w_{\n,1}(z)}{S_\n^{(1)}(z)},\frac{w_{\n,2}(z)}{S_\n^{(2)}(z)}\right)\boldsymbol \Pi(z) \boldsymbol S(\infty),
\]
which also serves as a definition of the matrix $ \boldsymbol\Pi(z) $. Then it follows from \eqref{Y3} that
\begin{equation}
\label{Y4}
-A_\n^{(i)}(z) = \big[\bigl(\boldsymbol Z^{-1}\bigr)^{\mathrm t}(z) \boldsymbol S(\infty) \boldsymbol \Pi^\mathrm{t}(z)\big]_{1,i+1} \frac{\tau_\n w_{\n,i}(z)}{\bigl(S_\n\Phi_\n\bigr)^{(i)}(z)}, \qquad z\in K_{\n,d}.
\end{equation}
Observe that all the jump matrices in \hyperref[rhz]{\rhz}(b) have determinant one. Since $ \boldsymbol Z(\infty)=\boldsymbol I $, we therefore get that $ \det(\boldsymbol Z(z))\equiv 1 $ in the entire extended complex plane. Hence,
\begin{equation}
\label{defhatB}
\big| \widehat Z_{j,k}(z)\big| \leq C_{2,r}\frac{\varepsilon_\n^{1/3}}{\delta^{10}}, \qquad
\begin{cases}
|z-\alpha_1|\geq 2r, & c(\n)\leq c^{\ast\ast}, \\
|z-\beta_2| \geq 2r, & c(\n)\geq c^\ast,
\end{cases} \qquad j,k\in\{0,1,2\},
\end{equation}
by \eqref{rhz-1} with perhaps modified constant $ C_{2,r} $, where \smash{$ \widehat Z_{j,k}(z) := \big[(\boldsymbol Z^{-1})^\mathrm{t}(z)\big]_{j+1,k+1} - \delta_{jk} $}. Moreover, as in \eqref{rhz-1}, \smash{$ \varepsilon_\n^{1/3} $} can be replaced by $ \varepsilon_\n $ when the parameters $ c(\n) $ are uniformly separated from $ c^*$, $c^{**} $. Notice also that $ \widehat Z_{j,k}(\infty) =0 $. Thus, we can write
\begin{align*}
\big[\bigl(\boldsymbol Z^{-1}\bigr)^{\mathrm t}(z)\boldsymbol S(\infty)\boldsymbol \Pi^\mathrm{t}(z)\big]_{1,i+1} ={}&
S_\n^{(0)}(\infty)\bigl(\Pi_\n^{(i)}(z) + \Pi_\n^{(i)}(z)\widehat Z_{0,0}(z) \\
& + s_{\n,1}^{-1}\Pi_{\n,1}^{(i)}(z)\widehat Z_{0,1}(z) + s_{\n,2}^{-1}\Pi_{\n,2}^{(i)}(z)\widehat Z_{0,2}(z) \bigr),
\end{align*}
$ z\in K_{\n,d} $, where, as before, $ s_{\n,l} = S_\n^{(0)}(\infty)/S_\n^{(l)}(\infty) $. Now, observe that
\[
\Pi_{\n,l}(\z) / \Pi_\n(\z) = -A_{\n,l}^{-1}\Upsilon_{\n,l}(\z), \qquad l\in\{1,2\},
\]
which follows by comparing zero/pole divisors in \eqref{Pin}, \eqref{Pini}, and the sentence after \eqref{Ups_n} as well as the normalizations at $ \infty^{(l)} $, see the sentence after \eqref{Pin}, \eqref{Pini}, and the display after~\eqref{Ups_n}. Therefore, \eqref{Y4} can be rewritten as
\begin{equation}
\label{Y5}
A_\n^{(i)}(z) = \Bigg(1 + \widehat Z_{0,0}(z) - \frac{\Upsilon_{\n,1}^{(i)}(z)}{s_{\n,1}A_{\n,1}} \widehat Z_{0,1}(z) - \frac{\Upsilon_{\n,2}^{(i)}(z)}{s_{\n,2}A_{\n,2}} \widehat Z_{0,2}(z) \Bigg) \mathcal A_{\n,i}(z).
\end{equation}
Recall that $ A_{c,i} $ are continuous and non-vanishing functions of the parameter $ c\in[0,1] $ except for~$ {A_{0,1}=0 }$ and $ A_{1,2} =0 $, see \eqref{AngPar2}, which satisfy \eqref{c_vanish}. Since \smash{$ \Upsilon_{\n,j\pm}^{(0)}(x) = \Upsilon_{\n,j\mp}^{(i)}(x) $} for~${ x\in\Delta_{\n,i} }$, we get from \eqref{Up-bounds1}, \eqref{Up-bounds2}, and \eqref{est_sni} that
\begin{equation}
\label{Y6}
\left| \frac{\Upsilon_{\n,3-i}^{(i)}(z)}{s_{\n,3-i}A_{\n,3-i}} \right| \lesssim \big| s_{\n,3-i}^{-1} \big| \lesssim 1
\end{equation}
uniformly for all $ \n $ and $ z $ in the whole extended complex plane by the maximum modulus principle for holomorphic functions). On the other hand, because $ \Upsilon_{\n,i}^{(i)}(z) $ has a simple pole at infinity, the same line of arguments yields that
\begin{equation}
\label{Y7}
\left| \frac{\Upsilon_{\n,i}^{(i)}(z)}{s_{\n,i}A_{\n,i}} \right| \lesssim \frac2{\beta_{\n,i}-\alpha_{\n,i}} \left| z - \frac{\beta_{\n,i}+\alpha_{\n,i}}2 + w_{\n,i}(z) \right|
\end{equation}
uniformly for all $ \n $ and $ z $ in the whole extended complex plane (the right-hand side above is simply the absolute value of a conformal map that takes the complement of $ \Delta_{\n,i} $ onto the complement of the unit disk normalized to take infinity into itself). Since the functions \smash{$ \widehat Z_{0,i}(z) $} vanish at infinity, it then follows from \eqref{c-rate}, \eqref{est_sni}, \eqref{defhatB}, \eqref{Y6}, and \eqref{Y7} that
\[
A_\n^{(i)}(z) / \mathcal A_{\n,i}(z) = 1 + o(s_{\n,i})
\]
uniformly on $ K_{\n,d} $, where $ o(1) $ has the same meaning as in Theorem~\ref{thm:1}. As the left-hand sides above is analytic off $ \Delta_{\n,i} $, this estimate, in fact, holds uniformly for $ \dist(z,\Delta_{\n,i})\geq d $ by the maximum modulus principle. This finishes the proof of the top formulae in \eqref{main2}.

Let now $ F_{\n,d,i} \subset \Delta_{\n,i} $ such that $ \dist(F_{\n,d,i}, E_{\n,i})\geq d $. As usual, we can adjust $ \delta $ so that $ F_{\n,d,i} $ does not intersect the closure of $ U_\n $. Relations \eqref{Y2} and \eqref{hatY1} give us
\[
\widehat{\boldsymbol Y}_{\pm}(x) = \boldsymbol C^{-1}\bigl(\boldsymbol Z^{-1}\bigr)^{\mathrm t}(x)\bigl(\boldsymbol M_{\pm}^{-1}\bigr)^\mathrm{t}(x)\boldsymbol D_{\pm}^{-1}(x)\bigl(\boldsymbol I\mp \rho_i^{-1}(x)\boldsymbol E_{1,i+1}\bigr), \qquad x\in F_{\n,d,i}.
\]
Similarly to the proof of the second formula in \eqref{main1}, observe that
\[
\mp\rho_i^{-1}(x)\Pi_{\n\pm}^{(0)}(x)/(S_\n\Phi_\n)_\pm^{(0)}(x) = \bigl(w_{\n,i\mp}\Pi_{\n\mp}^{(i)}\bigr)(x)/(S_\n\Phi_\n)_\mp^{(i)}(x), \qquad x\in\Delta_{\n,i},
\]
by \eqref{szego-pts2}. Hence, analogously to \eqref{Y5}, the above two formulae yield that
\begin{align*}
A_\n^{(i)}(x) ={}& \Bigg(1 + \widehat Z_{0,0}(x) - \frac{\Upsilon_{\n,1\pm}^{(i)}(x)}{s_{\n,1}A_{\n,1}} \widehat Z_{0,1}(x) - \frac{\Upsilon_{\n,2\pm}^{(i)}(x)}{s_{\n,2}A_{\n,2}}\widehat Z_{0,2}(x) \Bigg) \mathcal A_{\n,i\pm}(x) \\
& + \Bigg(1 + \widehat Z_{0,0}(x) - \frac{\Upsilon_{\n,1\mp}^{(i)}(x)}{s_{\n,1}A_{\n,1}} \widehat Z_{0,1}(x) - \frac{\Upsilon_{\n,2\mp}^{(i)}(x)}{s_{\n,2}A_{\n,2}} \widehat Z_{0,2}(x) \Bigg) \mathcal A_{\n,i\mp}(x).
\end{align*}
Since the estimates in \eqref{Y6} and \eqref{Y7} do hold on $ F_{\n,d,i} $, the bottom formula in \eqref{main2} follows.

\section[Proof of Theorems 1.10 and 1.11]{Proof of Theorems~\ref{thm:3} and~\ref{thm:4}}

We retain the notation from the previous section.

\subsection[Asymptotics of a\_n,i]{Asymptotics of $\boldsymbol{a_{\n,i} }$}

In this subsection, we prove the first two equalities of the top line of \eqref{main4} and the first asymptotic formula in \eqref{main3}. To this end, we need the formula
\begin{equation}
\label{9.1.1}
h_{\n,i} = a_{\n,i} h_{\n-\vec e_i,i},
\end{equation}
which is obtained by multiplying the first recurrence relation in \eqref{rec_rel} by $ x^{n_i-1} $, integrating it against $ {\rm d}\mu_i(x) $, and recalling \eqref{connection}.

To extract asymptotics of $ h_{\n,i} $ we use \eqref{Rni}. Formulae \eqref{Y}, \eqref{matrix-M}, and \eqref{Y1} yield that
\begin{align*}
R_\n^{(i)}(z) &= [\boldsymbol Y(z)]_{1,i+1} =[(\boldsymbol{ZM})(z)]_{1,i+1}[\boldsymbol D(z)]_{i+1,i+1} \\
& = \bigl( 1+ Z_{0,0}(z) + s_{\n,1}\Upsilon_{\n,1}^{(i)}(z)Z_{0,1}(z) + s_{\n,2}\Upsilon_{\n,2}^{(i)}(z)Z_{0,2}(z) \bigr) \frac{ [\boldsymbol C]_{1,1}}{S_\n^{(0)}(\infty)} \frac{(S_\n\Phi_\n)^{(i)}(z)}{w_{\n,i}(z) }
\end{align*}
in a neighborhood of infinity. We get from \eqref{c_vanish}, \eqref{c-rate}, \eqref{est_sni}, \eqref{Y6}, and \eqref{Y7} that
\[
R_\n^{(i)}(z) = \bigl( 1 + o(s_{\n,i}) \bigr) \frac{ [\boldsymbol C]_{1,1}}{S_\n^{(0)}(\infty)} \frac{(S_\n\Phi_\n)^{(i)}(z)}{w_{\n,i}(z) }
\]
in the vicinity of infinity, where $ o(1) $ is exactly the same as described in Theorem~\ref{thm:1}. Hence, it follows from the second equality in \eqref{Rni} and the very definition of the matrix $ \boldsymbol C $ in \eqref{matrix-CD} that
\begin{equation}
\label{9.1.2}
h_{\n,i} = \frac{1 + o(s_{\n,i})}{s_{\n,i}} \frac{[\boldsymbol C]_{1,1}}{[\boldsymbol C]_{i+1,i+1}}.
\end{equation}
On the other hand, since $ 1/h_{\n-\vec e_i,i} $ is the leading coefficient of $ A_\n^{(i)}(z) $, we get from the first two formulae of \eqref{main2}, the limits stated right after \eqref{Pin}, the definition of $ \boldsymbol C $, and \eqref{est_sni} that
\begin{equation}
\label{9.1.3}
\frac1{h_{n-\vec e_i,i}} = \bigl( 1 + o(s_{\n,i}) \bigr) s_{\n,i}A_{\n,i} \frac{[\boldsymbol C]_{i+1,i+1}}{[\boldsymbol C]_{1,1}}.
\end{equation}
Since $ s_{\n,i}^2A_{\n,i} \sim 1$, the first claim in \eqref{main3} follows by plugging \eqref{9.1.2} and \eqref{9.1.3} into \eqref{9.1.1}.

The above formulae also allow us to prove the validity of the first two equations in the top line of \eqref{main4}. First of all, let us observe that the differentiability of $ A_{c,i} $ and $ B_{c,i} $ as functions of $ c $ on $ (0,c^*) \cup (c^{**},1) $ was established in \cite[Section~4]{ApKoz20} based on an explicit parametrization of the Riemann surfaces $ \RS_c $ established in \cite{ApKalLysTul09,LysTul17,LysTul18}.

Fix $ c \in (0,c^*) \cup (c^{**},1) $ and let $ \mathcal N_c $ be a sequence of multi-indices such that $ c(\n) \to c $ as~$ {|\n| \to\infty }$, $ \n\in \mathcal N_c $. Clearly, we also have $ c(\n-\vec e_i) \to c $ as $ |\n| \to\infty $, $ \n\in \mathcal N_c $. Notice that the numbers $ s_{\n,i} $ and $ s_{\n-\vec e_i,i} $ are uniformly bounded along $ \mathcal N_c $ and $ s_{\n,i}/s_{\n-\vec e_i,i}\to 1 $ as $ |\n| \to\infty $, $ \n\in \mathcal N_c $, due to the continuous dependence of the Szeg\H{o} functions on the parameter. Hence, we get from \eqref{9.1.1}--\eqref{9.1.3} that
\[
(1+o(1))A_{\n,i} = a_{\n,i} = (1+o(1))\frac{[\boldsymbol C_\n]_{1,1}}{[\boldsymbol C_\n]_{i+1,i+1}}\frac{[\boldsymbol C_{\n-\vec e_i}]_{i+1,i+1}}{[\boldsymbol C_{\n-\vec e_i}]_{1,1}}
\]
as $ |\n| \to\infty $, $ \n\in \mathcal N_c $, where we explicitly indicate the dependence of the matrices $ \boldsymbol C $ on the multi-index $ \n $. Using \eqref{Phin}, we then get that
\[
1+o(1) = \left( \frac{A_{\n,i}}{A_{\n-\vec e_i,i}} \right)^{n_i-1} \left( \frac{B_\n}{B_{\n-\vec e_i}} \right)^{n_{3-i}}
\]
as $ |\n| \to\infty $, $ \n\in \mathcal N_c $, where $ B_\n = B_{\n,2} - B_{\n,1} $. By taking logarithms of both sides and using the mean-value theorem, we get that
\[
o(1) = (c(\n) - c(\n-\vec e_i)) \left((n_i-1) \frac{A_{\xi(\n),i}^\prime}{A_{\xi(\n),i}} + n_{3-i} \frac{B_{\eta(\n)}^\prime}{B_{\eta(\n)}} \right)
\]
as $ |\n| \to\infty $, $ \n\in \mathcal N_c $, for some $ \xi(\n)$, $\eta(\n) $ that lie between $ c(\n) $ and $ c(\n-\vec e_i) $. Notice that
\[
c(\n) - c(\n-\vec e_i) = \frac{(-1)^{i-1} n_{3-i}}{|\n|(|\n|-1)}.
\]
Since $ n_{3-i}/|\n| $ approaches either $ c $ or $ 1-c $ along $ \mathcal N_c $, it therefore holds that
\[
o(1) = \frac{n_i-1}{|\n|}\frac{A_{\xi(\n),i}^\prime}{A_{\xi(\n),i}} + \frac{n_{3-i}}{|\n|} \frac{B_{\eta(\n)}^\prime}{B_{\eta(\n)}}
\]
as $ |\n| \to\infty $, $ \n\in \mathcal N_c $. Because $ \xi(\n),\eta(\n) \to c $ as $ |\n| \to\infty $, $ \n\in \mathcal N_c $, the first two differential equation in the top line of \eqref{main4} follow by taking the limit in the above equality.

\subsection[Proof of (1.35)]{Proof of (\ref{main5})}

In this subsection, we prove \eqref{main5} and the last equality in the top line of \eqref{main4}. Let us set~${ R_1(c) := c^{-2}A_{c,1}} $ and $ R_2(c) := (1-c)^{-2}A_{c,2} $, which are continuous non-vanishing functions of the parameter $ c\in[0,1] $, see \eqref{c_vanish}, that are continuously differentiable on $ (0,1)\setminus\{c^*,c^{**}\}$. It follows from \eqref{AcBc} that
\begin{equation}
\label{9.3.1}
R_1^\prime(c) + R_2^\prime(c) = 2B_cB_c^\prime.
\end{equation}
On the other hand, expressing $ A_{c,i}^\prime $ through $ A_{c,i} $, $ B_c $, and $B_c^\prime $ by using the first two equalities in the top line of \eqref{main4} together with adding these expressions up gives
\begin{equation}
\label{9.3.2}
A_{c,1}^\prime + A_{c,2}^\prime = -\frac{B_c^\prime}{B_c}\left( \frac{1-c}c A_{c,1} + \frac c{1-c} A_{c,2} \right) = -c(1-c) B_c B_c^\prime,
\end{equation}
where we also used \eqref{AcBc} for the last step. Combining equations \eqref{9.3.1} and \eqref{9.3.2} to eliminate~$ B_cB_c^\prime $ as well as dividing by $ R_1(c) $ yields
\begin{equation}
\label{9.3.3}
(1-c)(2-c)R(c)\frac{R_2^\prime(c)}{R_2(c)} + c(1+c)\frac{R_1^\prime(c)}{R_1(c)} = -4c + 4(1-c)R(c).
\end{equation}
On the other hand, the first two relations in top line of \eqref{main4} can be rewritten as
\begin{equation}
\label{9.3.4}
(1-c)^2\frac{R_2^\prime(c)}{R_2(c)} - c^2\frac{R_1^\prime(c)}{R_1(c)} = 2.
\end{equation}
Equations \eqref{9.3.3} and \eqref{9.3.4} form a two by two linear system whose solution is given by
\begin{equation}
\label{9.3.5}
\frac{R_2^\prime(c)}{R_2(c)} = 2\frac{2c+1+2cR(c)}{1-c^2+c(2-c)R(c)}
\end{equation}
and
\begin{equation}
\label{9.3.6}
\frac{R_1^\prime(c)}{R_1(c)} = -2\frac{2(1-c)+(3-2c)R(c)}{1-c^2+c(2-c)R(c)}.
\end{equation}
Equation \eqref{main5} now follows by taking the difference between \eqref{9.3.5} and \eqref{9.3.6}. Finally, the last equality in the top line of \eqref{main4} comes from \eqref{9.3.6} and the identity
\[
\frac{B_c^\prime}{B_c} = -\frac c{1-c} \frac{A_{c,1}^\prime}{A_{c,1}} = -\frac2{1-c} - \frac c{1-c}\frac{R_1^\prime(c)}{R_1(c)}.
\]

\subsection[Differentiability of chi\_c and beta\_c,1]{Differentiability of $ \boldsymbol{\chi_c }$ and $\boldsymbol{\beta_{c,1}} $}

In this subsection, we establish several facts that will be needed in the remainder of the proof of Theorems~\ref{thm:3} and~\ref{thm:4}. We investigate only the case $ c\in(0,c^*) $ as the behavior for $ c\in(c^{**},1) $ can be deduced similarly.

Since the symmetric functions of the branches of $ \chi_c(\z) $ are necessarily polynomials, we can use \eqref{chi} and \eqref{AngPar1} to derive the cubic equation satisfied by $ \chi_c $. This equation, after some straightforward algebraic simplifications, can be written as
\begin{equation}
\label{alg eq for chi}
z = \chi_c(\z) + \frac{A_{c,1}}{\chi_c(\z)-B_{c,1}} + \frac{A_{c,2}}{\chi_c(\z)-B_{c,2}}.
\end{equation}
Given $ s\in \overline\C $, the above equation can be interpreted as
\begin{equation}
\label{derivative zc}
 z(c,s) := \pi\bigl( \chi_c^{-1}(s)\bigr) = s + \frac{A_{c,1}}{s-B_{c,1}} + \frac{A_{c,2}}{s-B_{c,2}}.
\end{equation}
In particular, we see from \eqref{main5}, the top line of \eqref{main4}, and the remark made after Theorem~\ref{thm:4} that $ \partial_c z(c,s) $ exists and is locally uniformly bounded in $ s\in \C\setminus\{B_{c,1},B_{c,2}\} $ (this estimate is also uniform in $ c \in (0,c^*) $ if $ s $ stays away from $ \{B_{c,i} : c\in [0,c^*] \}$).

As pointed out right before \eqref{AcBc}, we have that $ \chi_c(\beta_{c,1}) = (1-c) B_{c,1} + cB_{c,2} $ when $ c\in(0,c^*) $. Thus, we get from \eqref{alg eq for chi} that
\[
\beta_{c,1} = (1-c) B_{c,1} + cB_{c,2} + \frac1{B_c}\left(\frac{A_{c,1}}{c}- \frac{A_{c,2}}{1-c}\right).
\]
This, of course, immediately shows that $ \beta_{c,1} $ is a differentiable function of $ c $ on $ (0,c^*) $. The above two observations will be sufficient for us to finish the proof of Theorem~\ref{thm:4}, that is, to prove the bottom line of \eqref{main4}.

Now we deduce from \eqref{AcBc} and the bottom line of \eqref{main4} (the use of the bottom line of \eqref{main4} is not really necessary but allows us to get a more compact formula) that
\[
\beta_{c,1}^\prime = (2c-1)B_c^\prime + \frac 1{B_c}\left( \frac{A_{c,1}^\prime}c - \frac{A_{c,2}^\prime}{1-c} \right) - \frac{B_c^\prime}{B_c^2}(cR_1(c)-(1-c)R_2(c)).
\]
Using the top line of \eqref{main4}, we obtain that
\[
\frac{A_{c,1}^\prime}c - \frac{A_{c,2}^\prime}{1-c} = \frac{B_c^\prime}{B_c} ( c R_2(c) - (1-c)R_1(c) ).
\]
Thus, using \eqref{AcBc}, \eqref{main5}, and \eqref{main4} once more, we arrive at
\begin{align}
\beta_{c,1}^\prime & = (2c-1)B_c^\prime + ( R_2(c) - R_1(c))\frac{B_c^\prime}{B_c^2} \nonumber \\
& = 2 ( cR_2(c) - (1-c)R_1(c))\frac{B_c^\prime}{B_c^2} = 6\frac{R(c)}{R^\prime(c)} \left( \frac{B_c^\prime}{B_c} \right)^2.
\label{derivative beta c}
\end{align}
Formula \eqref{derivative beta c} clearly shows that $ \beta_{c,1}^\prime $ is a bounded continuous function of $ c\in(0,c^*) $ (it also shows that $ \beta_{c,1} $ is an increasing function of $ c $).

\subsection[Differentiability of sigma(c)]{Differentiability of $\boldsymbol{\sigma(c)} $}

Before we go back to the proof of Theorems~\ref{thm:3} and~\ref{thm:4}, we need to analyze the behavior of one more auxiliary quantity, namely $ \sigma(c) $, defined by
\begin{equation}
\label{def sigma c}
S_c^{(0)}(z) = S_c^{(0)}(\infty) \left( 1 + \frac{\sigma(c)}z + \mathcal O\left(\frac1{z^2}\right) \right), \qquad z\to\infty.
\end{equation}
More precisely, we need to establish its differentiability with respect to $ c $ and the boundedness of $ \sigma^\prime(c) $. Again, we only consider the case $ c\in(0,c^*) $ as the behavior for $ c\in(c^{**},1) $ can be deduced similarly.

We shall utilize a representation of $ S_c(\z) $ different from \eqref{Szego}. To this end, let $ S_{\rho_2}(z) $ be given by \eqref{classical_szego}. Set $ \hat\rho_1(x) := \rho_1(x)S_{\rho_2}(x) $ and define
\[
S_{c,1}(z) := {\rm e}^{F_{c,1}(z)}, \qquad F_{c,1}(z) := \frac{w_{c,1}(z)}{2\pi\ic}\int_{\Delta_{c,1}}\frac{\log(\hat\rho_1w_{c,1+})(x)}{z-x}\frac{{\rm d}x}{w_{c,1+}(x)},
\]
which is the standard Szeg\H{o} function of $ (\hat\rho_1w_{c,1+})(x)$ considered as a weight on $ \Delta_{c,1} $. That is, it is a non-vanishing holomorphic function off $ \Delta_{c,1} $ (including at infinity) with continuous traces on both sides of $ (\alpha_1,\beta_{c,1}) $ that satisfy
$
S_{c,1+}(x)S_{c,1-}(x) (\hat\rho_{1|\Delta_{c,1}}w_{c,1+})(x) \equiv 1
$
and quarter-root singularities at $ \alpha_1$, $\beta_{c,1} $. Next, recall that $ \chi_c(\boldsymbol\Delta_2) $ is a Jordan curve that approaches the circle~${ \chi_0(\boldsymbol\Delta_2) := \bigl\{|z-B_{0,2}|^2=A_{0,2}\bigr\}} $ as $ c\to0 $, see \eqref{AngPar2} and \eqref{chi_curves}. We orient $ \chi_c(\boldsymbol\Delta_2) $ clockwise. The conformal map $ \chi_c(\z) $ maps \smash{$ \RS_c^{(2)}\setminus\boldsymbol\Delta_2 $} and \smash{$ \RS_c\setminus \RS_c^{(2)} $} onto the interior and exterior domains of $ \chi_c(\boldsymbol\Delta_2) $ with \smash{$ \infty^{(0)} $} mapped into $ \infty $, see \eqref{chi}. Recall \eqref{derivative zc}. Set
\[
\mathcal S_{c,2}(\z) := \exp\left\{\oint_{\chi_c(\boldsymbol\Delta_2)}\frac{F_{c,1}(z_c(s))}{s-\chi_c(\z)}\frac{{\rm d}s}{2\pi\ic}\right\},
\]
which is a sectionally analytic and non-vanishing function in $ \RS_c\setminus\boldsymbol\Delta_2 $, $ \mathcal S_{c,2}\bigl(\infty^{(0)}\bigr) = 1 $, and $ \mathcal S_{c,2-}(\x) = \mathcal S_{c,2+}(\x) S_{c,1}(x) $ on $ \boldsymbol\Delta_2 $ (according to our chosen orientation of $ \boldsymbol\Delta_2 $, the positive approach to $ \boldsymbol\Delta_2 $ is from \smash{$ \RS_c^{(0)} $}). The properties of $ S_{c,1}(z) $ and $ \mathcal S_{c,2}(\z) $ readily yield that
\[
S_c(\z) = k_c\begin{cases}
S_{c,1}(z)\mathcal S_{c,2}(\z), & \z\in\RS_c^{(0)}, \\
\mathcal S_{c,2}(\z)/S_{c,1}(z), & \z\in\RS_c^{(1)}, \\
\mathcal S_{c,2}(\z), & \z\in\RS_c^{(2)},
\end{cases}
\]
where $ k_c $ is a normalizing constant \big(one can readily check that the right-hand side above satisfies~\eqref{szego-pts2} and $ k_c $ is there to achieve the normalization \smash{$ \bigl(S_c^{(0)}S_c^{(1)} S_c^{(2)} \bigr)(z)\equiv 1 $}\big). Therefore,
\begin{equation}
\label{expression sigma c}
\sigma(c) = F_{c,1}^\prime(\infty) - \frac1{2\pi\ic}\oint_{\chi_0(\boldsymbol\Delta_2)}F_{c,1}(z(c,s)){\rm d}s,
\end{equation}
where we understand the derivative at infinity in local coordinates, that is, $ f(z) = f(\infty) + f^\prime(\infty)/z + \cdots $, and we moved the curve of integration from $ \chi_c(\boldsymbol\Delta_2) $ to the circle $ \chi_0(\boldsymbol\Delta_2) $ by the analyticity of $ F_{c,1}(z_c(s)) $ (if $ \chi_0(\boldsymbol\Delta_2) $ does not work for all $ c\in(0,c^*) $, we can partition $ (0,c^*) $ into finitely many overlapping intervals and on each of them similarly replace $ \chi_c(\boldsymbol\Delta_2) $ by a curve independent of $ c $; this will not alter the forthcoming computations).

Let $ \ell_c(y) = \frac{\beta_{c,1}-\alpha_1}2(y+1) + \alpha_1 $ be the linear function with positive leading coefficient that takes $ [-1,1] $ onto $ \Delta_{c,1} $. Then,
\[
F_{c,1}(\ell_c(z)) = \frac{\sqrt{z^2-1}}{2\pi}\int_{-1}^1 \frac{\log\bigl(\ic\hat\rho_1(\ell_c(y)) \frac{\beta_{c,1}-\alpha_1}2\sqrt{1-y^2} \bigr)}{y-z}\frac{{\rm d}y}{\sqrt{1-y^2}}.
\]
It can be readily verified that the above expression can be rewritten as
\[
\int_{-1}^1 \frac{\sqrt{z^2-1}}{y-z}\frac{\log \hat\rho_1(\ell_c(y)){\rm d}y}{2\pi\sqrt{1-y^2}} + \int_{-1}^1 \frac{\sqrt{z^2-1}}{y-z}\frac{\log\sqrt{1-y^2}{\rm d}y}{2\pi\sqrt{1-y^2}} - \frac12\log\left(\frac \ic2 (\beta_{c,1}-\alpha_1)\right).
\]
Notice that we do not need to worry about the logarithmic term as it is constant in $ z $ an will disappear after integration on the circle $ \chi_0(\boldsymbol\Delta_2 ) $. The same reason allows us to subtract $ 1 $ from the kernel in the above integrals without altering the expression in \eqref{expression sigma c}. Hence, we need to analyze integrals of the form
\[
I_{d,k}(z) := \int_{-1}^1 \left(\frac{\sqrt{z^2-1}}{y-z} - 1 \right) \frac{d(\ell_c(y))k(y){\rm d}y}{2\pi\sqrt{1-y^2}},
\]
where either $ d(x) = \log\hat\rho_1(x) $ and $ k(y)= 1 $ or $d(x) = 1 $ and $ k(y)= \log\sqrt{1-y^2} $. It is not hard to see that $ I_{d,k}(z) = \mathcal O\bigl(z^{-1}\bigr) $ and $ I_{d,k}^\prime(z) = \mathcal O\bigl(z^{-2}\bigr) $, where $ \mathcal O $ terms are uniform in $ c $. Moreover, it follows from \eqref{derivative beta c} that the above expression is differentiable with respect to $ c $ and it holds that
$
(\partial_cI_{d,k})(z) = \beta_{c,1}^\prime I_{d^\prime,(\cdot+1)k/2}(z) = \mathcal O\left(z^{-1}\right)$,
where $ \mathcal O $ term is again uniform in $ c $. Hence, we have that
\[
\partial_cI_{d,k}\left(\ell_c^{-1}(z)\right) = (\partial_cI_{d,k})\left(\ell_c^{-1}(z)\right) + I_{d,k}^\prime\left(\ell_c^{-1}(z)\right) \frac{-2\beta_{c,1}^\prime(z-\alpha_1)}{(\beta_{c,1}-\alpha_1)^2} = \mathcal O(1),
\]
where $ \ell_c^{-1}(z) $ is the inverse (not reciprocal) function of $ \ell_c(z) $ and $ \mathcal O $ term is again uniform in $ c $ and locally uniform in $ z $. Hence, we have that
\[
\partial_c \oint_{\chi_0(\boldsymbol\Delta_2)}F_{c,1}(z(c,s)){\rm d}s = \mathcal O(\partial_c z(c,s)) = \mathcal O(1),
\]
uniformly in $ c $, where we used boundedness of $ \partial_c z(c,s) $ to reach the last conclusion, see the sentence after \eqref{derivative zc}. Similarly, since integrals of odd functions on $ (-1,1) $ are necessarily zero, we have that
\[
F_{c,1}^\prime(\infty) = -\frac{\beta_{c,1}-\alpha_1}{4\pi} \int_{-1}^1 \log \hat\rho_1(\ell_c(y)) \frac{y{\rm d}y}{\sqrt{1-y^2}}.
\]
Considerations virtually identical to the ones presented above show that $ F_{c,1}^\prime(\infty) $ is differentiable with respect to $ c $ and that $ \partial_c F_{c,1}^\prime(\infty) = \mathcal O(1) $ uniformly in $ c $. Hence, we now deduce from \eqref{expression sigma c} that $ \sigma(c) $ is differentiable and $ \sigma^\prime(c) $ is bounded for $ c\in(0,c^*) $.

\subsection[Asymptotics of b\_n,i]{Asymptotics of $\boldsymbol{b_{\n,i} }$}

In this subsection, we prove the bottom formula in \eqref{main4} and the second asymptotic formula in~\eqref{main3}, thus, finishing the proof of Theorems~\ref{thm:3} and~\ref{thm:4}. Both proofs rely on the formula
\begin{equation}
\label{9.2.1}
b_{\n,i} = \lim_{z\to\infty} \left( z - \frac{P_{\n+\vec e_i}(z)}{P_\n(z)} \right),
\end{equation}
which readily follows from the first relation in \eqref{rec_rel}.

Consider \eqref{Y} with $ \n $ replaced by $ \n+\vec e_i $. It holds that $ h_{\n,i}[\boldsymbol Y(z)]_{i+1,1} = P_\n(z) $. Hence, similarly to \eqref{Pn}, we have that
\begin{align*}
P_\n(z) ={}& h_{\n,i} \mathcal P_{\n+\vec e_i}(z)\frac{[\boldsymbol C]_{i+1,i+1}}{[\boldsymbol C]_{1,1}} \\
&\times\bigl( Z_{i,0}(z) + s_{\n+\vec e_i,i}\Upsilon_{\n+\vec e_i,i}^{(0)}(z)(1+Z_{i,i}(z)) + s_{\n+\vec e_i,3-i}\Upsilon_{\n+\vec e_i,3-i}^{(0)}(z)Z_{i,3-i}(z) \bigr)
\end{align*}
for $ z $ in a neighborhood of infinity. Recall that the functions $ Z_{i,k}(z) $ and $ \Upsilon_{\n+\vec e_i,k}(z) $ vanish at infinity and
\[
\Upsilon_{\n+\vec e_i,i}^{(0)}(z) = \frac{A_{\n+\vec e_i,i}}z\left(1 + \frac{B_{\n+\vec e_i,i}}z + \mathcal O\left(\frac1{z^2}\right) \right)
\]
as $ z\to\infty $, see \eqref{Ups_n}. Since $ P_\n(z) $ is a monic polynomial, it therefore follows from \eqref{rhz-1} and \eqref{Up-bounds} that
\begin{equation}
\label{9.2.2}
P_\n(z) = \left( \frac 1z + \frac{s_{\n+\vec e_i,i}A_{\n+\vec e_i,i}B_{\n+\vec e_i,i}+o(1)}{s_{\n+\vec e_i,i}A_{\n+\vec e_i,i} + o(1)} \frac1{z^2} + \mathcal O\left(\frac1{z^3}\right) \right) \mathcal P_{\n+\vec e_i}(z)
\end{equation}
in the vicinity of infinity, where $ o(1) $ term in the denominator above is the constant next $ 1/z $ term of $ Z_{i,0}(z) $. Of course, it also holds that
\begin{equation}
\label{9.2.3}
P_{\n+\vec e_i}(z) = (1 + o(1) ) \mathcal P_{\n+\vec e_i}(z)
\end{equation}
for $ z $ around infinity. Plugging \eqref{9.2.2} and \eqref{9.2.3} into \eqref{9.2.1} gives
\[
b_{\n,i} = \frac{s_{\n+\vec e_i,i}A_{\n+\vec e_i,i}B_{\n+\vec e_i,i}+o(1)}{s_{\n+\vec e_i,i}A_{\n+\vec e_i,i} + o(1)} + o(1).
\]
The above formula is insufficient to prove the second formula of \eqref{main3} as $ s_{\n+\vec e_i,i}A_{\n+\vec e_i,i} \sim c(\n) $ as $ c(\n) \to 0 $ by \eqref{c_vanish} and \eqref{est_sni}. However, if $ c\in(0,c^*) \cup (c^{**},1) $ and $ \mathcal N_c $ is as in the second paragraph after \eqref{9.1.3}, then
\begin{equation}
\label{9.2.4}
b_{\n,i} = B_{\n+\vec e_i,i}+o(1) \qquad \text{as} \ |\n| \to\infty, \quad \n\in \mathcal N_c.
\end{equation}

On the other hand, \eqref{9.2.1}, \eqref{9.2.3}, and \eqref{9.2.3} with $ \n+\vec e_i $ replaced by $ \n $, give us that
\begin{equation}
\label{9.2.5}
b_{\n,i} = (\mathcal P_\n)_1 - (\mathcal P_{\n+\vec e_i})_1 + o(1)
\end{equation}
for all $ \varepsilon_\n $ small enough, where the error term is as in Theorem~\ref{thm:1} and we write $ \mathcal P_\n(z) = z^{|\n|} + (\mathcal P_\n)_1 z^{|\n|-1} + \cdots $. According to \eqref{Phin} and \eqref{Fn}, it holds that
\[
\mathcal P_\n(z) = \bigl( \chi_\n^{(0)}(z) -B_{\n,1} \bigr)^{n_1} \bigl( \chi_\n^{(0)}(z)-B_{\n,2} \bigr)^{n_2} S_\n^{(0)}(z) / S_\n^{(0)}(\infty).
\]
Then we get from \eqref{chi}, \eqref{AngPar1}, and \eqref{9.2.5} that
\begin{equation}
\label{9.2.6}
b_{\n,i} = B_{\n+\vec e_i,i} + n_1( B_{\n+\vec e_i,1} - B_{\n,1} ) + n_2( B_{\n+\vec e_i,2} - B_{\n,2} ) + \sigma_\n - \sigma_{\n+\vec e_i} + o(1),
\end{equation}
where $ \sigma_\n=\sigma(c(\n))$, see \eqref{def sigma c}, and the error term is as in Theorem~\ref{thm:1}. When $ c(\n),c(\n+\vec e_i)\in[c^*,c^{**}] $, $ \RS_\n=\RS_{\n+\vec e_i} $ and therefore \eqref{9.2.6} simply reduces to the second formula in \eqref{main3}. To deal with the remaining cases, recall that $ \sigma(c) $ is a differentiable function of $ c $ on $ (0,c^*)\cup(c^{**},1) $ with a~bounded derivative there. Since
\[
c(\n+\vec e_i) - c(\n) = \frac{(-1)^{i-1}n_{3-i}}{|\n|(|\n|+1)},
\]
it follows from \eqref{9.2.6} and the mean-value theorem that
\begin{equation}
\label{9.2.7}
b_{\n,i} - B_{\n+\vec e_i,i} + o(1) = \frac{(-1)^{i-1}n_{3-i}}{|\n|+1}\left( c(\n)B^\prime_{\xi(\n),1} + (1-c(\n))B_{\eta(\n),2}^\prime + \frac{\sigma^\prime(\zeta(\n))}{|\n|}\right),
\end{equation}
where the error term is as in Theorem~\ref{thm:1} and $ \xi(\n)$, $\eta(\n)$, $\zeta(\n) $ lie between $ c(\n) $ and $ c(\n+\vec e_i) $.

To finish the proof of Theorem~\ref{thm:4}, let $ c $ and $ \mathcal N_c $ be as \eqref{9.2.4}. The last relation in \eqref{main4} now follows from \eqref{9.2.4} and \eqref{9.2.7} since $ \xi(\n),\eta(\n),\zeta(\n) \to c $ as $ |\n|\to\infty $, $ \n\in\mathcal N_c $.

As pointed out right after Theorem~\ref{thm:4}, \eqref{main5} and the top line of \eqref{main4} imply infinite differentiability and boundedness of all the derivatives. Thus, we get from the mean-value theorem that
\[
\big| B^\prime_{\xi(\n),1} - B^\prime_{c(\n),1} \big|,\big| B^\prime_{\eta(\n),2} - B^\prime_{c(\n),2} \big| \lesssim |\n|^{-1}
\]
with a constant independent of $ \n $. Hence, we get from the bottom relation in \eqref{main4} and \eqref{9.2.7} that
$
b_{\n,i} - B_{\n+\vec e_i,i} + o(1) = \mathcal O\bigl(|\n|^{-1}\bigr)
$
uniformly in $ \n $ when at least one of the numbers $ c(\n)$, $c(\n+\vec e_i) $ belongs to $ (0,c^*] \cup [c^{**},1) $ (notice that if for example $ c(\n+\vec e_1)<c^*<c(\n) $, then $ B_{\n,i}=B_{c^*,i} $, $ \sigma_\n=\sigma(c^*) $ and therefore \eqref{9.2.6} still yields \eqref{9.2.7} with $ \xi(\n),\eta(\n),\zeta(\n) \in (c(\n+\vec e_1),c^*) $). This finishes the proof of Theorem~\ref{thm:3}.

\subsection*{Acknowledgements}

The research was supported in part by a grant from the Simons Foundation, CGM-706591. The author is grateful to the anonymous referees for their careful reading of the manuscript and helpful suggestions.

\pdfbookmark[1]{References}{ref}
\LastPageEnding

\end{document}